\documentclass[12pt,a4paper]{amsart}

\usepackage{amsmath,amsthm,amsfonts,amssymb,array,amscd}
\usepackage{amsmath,geometry,amssymb,amsfonts,amsthm,graphicx,enumerate,latexsym,tabularx,amscd}

\usepackage{refcheck} 
\norefnames
\nocitenames
\theoremstyle{plain}
\usepackage{amsthm}
 \newtheorem{thm}{Theorem}[section]
 \newtheorem{prop}[thm]{Proposition}
 \newtheorem{lem}[thm]{Lemma}
 \newtheorem{cor}[thm]{Corollary}
 \newtheorem{conj}[thm]{Conjecture}
\theoremstyle{definition}
 \newtheorem{exm}[thm]{Example}
 \newtheorem{dfn}[thm]{Definition}
 \newtheorem{rem}[thm]{Remark}
 \numberwithin{equation}{section}
\theoremstyle{definition}
\theoremstyle{remark}
 \numberwithin{equation}{section}

\renewcommand{\le}{\leqslant}\renewcommand{\leq}{\leqslant}
\renewcommand{\ge}{\geqslant}\renewcommand{\geq}{\geqslant}
\renewcommand{\setminus}{\smallsetminus}
\usepackage{tikz}
\usepackage[pdftex]{hyperref}
\usetikzlibrary{matrix,arrows}

\newcommand{\bbA}{\mathbb{A}}
\newcommand{\bbC}{\mathbb{C}}
\newcommand{\bbF}{\mathbb{F}}

\newcommand{\bbQ}{\mathbb{Q}}
\newcommand{\bbR}{\mathbb{R}}

\newcommand{\bbZ}{\mathbb{Z}}   

\renewcommand{\and}{\quad \mbox{and} \quad}  
\renewcommand{\le}{\leqslant}\renewcommand{\leq}{\leqslant}
\renewcommand{\ge}{\geqslant}\renewcommand{\geq}{\geqslant}
\renewcommand{\setminus}{\smallsetminus}
\title{existence of non-abelian local constants, and their properties}


\author[Biswas]{\bfseries Sazzad Ali Biswas}

\address{
Department of Mathematics \\ 
SRM University AP\\ 
Guntur, 522502\\
India}
\email{sazzadali.b@srmap.edu.in, sazzad.jumath@gmail.com}

\begin{document}

\vspace{18mm}
\setcounter{page}{1}
\thispagestyle{empty}

\begin{abstract}

In his Ph.D. thesis, John Tate attached the (abelian) local constants to the characters of a non-Archimedean local field
of characteristic zero. Robert Langlands proved the existence theorem of a non-abelian local constant 
of a higher-dimensional complex local Galois representation. In 1990, Helmut Koch summarized Langlands' strategy for 
the existence of a non-abelian local constant (group theoretically).
The Brauer induction formula plays a crucial role 
in Langlands' proof. Robert Boltje gives a canonical version of the Brauer induction formula.
In this paper, we review Langlands' strategy using Boltje's canonical Brauer induction formula.
We then review various
properties of local constants, some applications, and open problems.

\end{abstract}

\maketitle

\section{\textbf{Introduction}}

Let $F$ be a non-Archimedean local field of characteristic zero, i.e., $F$ is a finite extension of $\bbQ_p$, where $p$
is a prime. Let $\overline{F}$ be an algebraic closure of $F$, and $G_F$ be the absolute Galois group of $F$. 
It can be proved that $G_F$ is a solvable group. One of the aims of this paper is to understand the answer to
the following question: \\

{\it How to attach local constants 
({\it or local root numbers, i.e., a special value (when $s=1/2$) of the epsilon factors})
with finite-dimensional complex representations of $G_F$?}\\

The local constant is an invariant of a local Galois representation that preserves under
the local Langlands correspondence, 
therefore, it is an important object in the Langlands
program. In the Langlands program, the local constant 
plays an important role as part of the detection machinery in 
the local Langlands conjecture.
Let $\chi:F^\times\to \bbC^\times$ be a nontrivial character of $F^\times$ 
and $\pi_F$ be a uniformizer of $F$. When the conductor of $\chi$ is {\bf zero}, we call
$\chi$ is {\bf unramified}, and otherwise $\chi$ is {\bf ramified}. With any character $\chi$ of $F^\times$, we can 
attach an $L$-function as follows: 
$$L(\chi)=\begin{cases}
           1 & \text{if $\chi$ is ramified},\\
           (1-\chi(\pi_F))^{-1} & \text{if $\chi$ is unramified}.
          \end{cases}
$$
John Tate showed that this $L(\chi)$ satisfies a functional equation (cf. Equation (\ref{eqn 4.4})), where 
the local constant 
$W(\chi,\psi)\in \bbC^\times$ appears, here $\psi$ is a nontrivial additive character of $F$.
Now the question is:\\

{\it Can we extend the notion of the local constant $W(\chi,\psi)$ for higher-dimensional local Galois 
representation such that the extended definition of local constants behaves well with 
respect to short exact sequences, change of multiplicative measure, induction (in dimension zero), and agree with Tate's 
definition in the one-dimensional case?} \\

The answer is {\bf YES} (see Theorem A of \cite{RL}). This was first proved by Robert Langlands in his unpublished article 
\cite{RL} by the local method, and the paper is about 300 pages. Later, P. Deligne gave a smaller proof of the existence of 
a non-abelian local constant using the global method (cf. \cite{D1}). In \cite{HK}, Helmut Koch described the strategy
of Langlands'
proof in the language of the extendible function. In Langlands' proof, the Brauer induction formula 
plays a crucial 
role. In \cite{RB}, Robert Boltje introduces a canonical Brauer induction formula. In Section 2, and Section 3, 
we review
Koch's paper but using {\bf Boltje's canonical induction formula} (see \cite{RB}, \cite{RB2}).\\

Before going to explain, here I mention
the main reasons for writing this paper:\\
{\it 1. Because of the {\bf uniqueness} of Boltje's canonical Brauer induction formula,
we can reduce many computations of Koch's\footnote{Recently, Helmut Koch and E.-W. Zink (cf. \cite{HKEWZ}) 
revisited Langlands' paper,
\cite{RL}.}
\cite{HK} to simple computations,\\
2.  To show the connections between the three different conventions of local constants (due to Langlands, Deligne, and
Bushnell $\&$ Henniart, see Subsection 4.4),\\
3. To review various properties of local constants, applications, and open problems.}\\

 

 Let $G$ be a finite group, and $R(G)$ be its Grothendieck group
(tensor product induces multiplication on it makes 
 it an abelian ring, and even it is called Grothendieck $\Lambda$-ring,
 that is, the free abelian group with basis the finite-dimensional complex representations of $G$, 
 quotient out by short exact sequences in the usual way).
 Let $R_{+}(G)$ denote the free abelian group with basis the set of isomorphism classes of pairs $(H,\varphi)$, where 
 $H$ is a subgroup of $G$, and $\varphi$ is a one-dimensional representation of $H$. Define a group homomorphism 
 $b_G:R_{+}(G)\to R(G)$ by sending $(H,\chi)\mapsto Ind_{H}^{G}(\chi)$. By the Brauer induction theorem
 this a {\bf surjective} group 
 homomorphism. 
 Because $b_G$ is surjective group homomorphism, then from the first 
 isomorphism theorem we have
 \begin{equation}
  R_{+}(G)/Ker(b_G)\cong R(G).
 \end{equation}

 Our aim is to define local constants for every element of $R(G)$. 
 For this, it is equivalent to studying the behavior of local root numbers on $Ker(b_G)$. 
Every element $\sigma\in Ker(b_G)$, we have
 $$\dim(b_G(\sigma))=\dim(0)=0.$$
 Therefore,  working on $Ker(b_G)$ is equivalent to 
working with virtual representations of dimension zero. Local constants for dimension zero virtual representations are 
invariant under induction, and it is easy to attach the non-abelian local constants for dimension zero virtual representations
(cf. the definition of extendible function \ref{Subsection extendible function}).
Therefore, studying $Ker(b_G)$ is important for the existence proof of 
non-abelian local constants (cf. 
 Remark \ref{Remark on Kernel of bG}).
In \cite{RL}, Langlands worked on the {\bf generator of $Ker(b_G)$} (an abelian group), and proved that
{\it the existence of a local constant for every element in $R(G)$ is 
 equivalent to checking that certain functions vanish on $Ker(b_G)$.}
However, P. Deligne's approach (cf. {\cite{D1}) is different from that of Langlands.
In Deligne's approach (global method), {\bf $Ker(b_G)$ is used but not its 
generators.}
 
 Deligne (following Langlands) wrote down three types of elements in $Ker(b_G)$ explicitly.
They are called {\bf elements of type I, type II, and type III.}
(cf. Definition \ref{Definition of type I, II, III}), and here
is the existence theorem of non-abelian local constants:
\begin{thm}[Langlands-Deligne, \cite{D1}]\label{Theorem 1.1}
 {\bf i.} If $G$ is abelian, the kernel $Ker(b_G)$ is generated as an abelian group by the relation of type I.\\
 {\bf ii.} If $G$ is nilpotent, the kernel $Ker(b_G)$ is generated as an abelian group by the relations of type 
 I, and type II.\\
 {\bf iii.} If $G$ is solvable, the kernel $Ker(b_G)$ is generated by relations of type I, type II, and type III.
\end{thm}

\begin{rem}
 Actually, Langlands did much more than this. He analyzed the generators of kernel $Ker(b_G)$ for an arbitrary finite 
 group. However, because our local Galois groups are solvable, for existence proof, we only need to study finite solvable
 groups. To analyze any functional properties of solvable groups, one should follow a standard pattern as follows:
 $$\text{{\bf first, abelian groups, then nilpotent groups, and finally solvable groups}}.$$
 In this paper, we also follow this pattern to prove Theorems \ref{Theorem 1.1}, and \ref{Theorem 1.2}.
\end{rem}

Article \cite{RL} is about the existence of the non-abelian local constants.
These local constants are certain functions that can be extended
from one-dimensional representations 
to higher-dimensional representations nicely. After proving Theorem \ref{Theorem 1.1},
he actually proved the following
theorem. This is the {\bf criterion of being an extendible function.}

\begin{thm}[Langlands, Theorem 3.1 of \cite{HK}]\label{Theorem 1.2}
 Let $G$ be a solvable group, and let $\mathcal{F}$ be a function on 
 $\mathcal{M}_G:=\{(H,\varphi)| H\leq G, \varphi\in \widehat{H}:=\{\chi:H\to\bbC^\times\}\}$
 with values in a multiplicative abelian group 
 $\mathcal{A}$ satisfying certain properties (cf. Equations (\ref{eqn 2.2.1}), and (\ref{eqn 2.2.2})).
  Let $K$ be a 
 normal subgroup of $B\le G$, and $A:=B/K$. Let $\chi$ be a character of $B$.
 Then $\mathcal{F}$ is (weakly) extendible to $\mathcal{M}_G^v:=\{(H,\rho)| H\le G,\quad\text{$\rho$ is a virtual representation
 of $H$}\}$ 
 if and only if the following three conditions are fulfilled for all subgroups $B$ of $G$. \\
 {\bf I.} In the situation of {\bf type I} (i.e., when $A$ is a cyclic group of order a prime $\ell$), we have
 \begin{equation}\label{eqn 5.15}
  \mathcal{F}((K,\chi))\prod_{\mu\in\widehat{A}}\mathcal{F}((A,\mu))=\prod_{\mu\in\widehat{A}}\mathcal{F}((B,\mu\chi)).
 \end{equation}
{\bf II.} In the situation of {\bf type II} (i.e., when $A$ is a central extension of an abelian group of order 
$\ell^2$: $\bbZ/\ell\bbZ\times\bbZ/\ell\bbZ\cong X_1\times X_2\twoheadleftarrow A\hookleftarrow Z$), we have 
\begin{equation}\label{eqn 5.16}
 \mathcal{F}((G_1,\chi_1\chi))\prod_{\mu\in\widehat{X_1}}\mathcal{F}((X_1,\mu))=
 \mathcal{F}((G_2,\chi_2\chi))\prod_{\mu\in\widehat{X_2}}\mathcal{F}((X_2,\mu)),
\end{equation}
where $G_i$ is the inverse image of $H_i$ in $B$ (cf. Definition \ref{Definition of type I, II, III} for {\bf type II}). \\
{\bf III.} In the situation of {\bf type III} (i.e., when $A=H\rtimes C$ is a semidirect product where $C\ne\{ 1\}$ is an 
abelian normal subgroup contained in all nontrivial abelian normal subgroups of $A$) 
(cf. Definition \ref{Definition of type I, II, III} for {\bf Type III}), we have
\begin{equation}\label{eqn 1.3}
 \mathcal{F}((G',\chi))\prod_{\mu\in T}\mathcal{F}((A_\mu',\mu))=\prod_{\mu\in T}\mathcal{F}((B_\mu',\mu\chi)),
\end{equation}
where $B_\mu'$ is the inverse image of $A_\mu'$ in $B$, and $G'$ is the inverse image of $H$ in $B$.
\end{thm}

\begin{rem}

For one-dimensional Galois representation of $G_F$, hence using class field theory for the character of $F^\times$, we have 
the notion of local constant. Properties of the local constants of one-dimensional representations 
are the same as the initial conditions
for defining extendible function (cf. Subsection \ref{Subsection extendible function}).
Therefore,  proving the above theorem is equivalent to proving the existence 
of non-abelian local constants, and Langlands explicitly did the same.

\end{rem}

{\bf Organization of the paper:} 
Including the Introduction section, there are six sections of this paper. In Section 2,
we first define Robert Boltje's canonical
Brauer induction formula, and mention all the necessary properties of canonical Brauer induction. Then, for any finite group, 
we give the definition of the {\it extendible function}.
Then, we study all the necessary conditions of a function to be extendible. 

In Section 3, we first study Kernel $Ker(b_G)$. Next, we prove the Theorem \ref{Theorem 1.2}.
In Section 4, we review all known properties of local constants. In Section 5,
we mention some applications of root numbers, and open 
problems.
In Section 6 (Appendix), we make some
remarks on the global constants (also known as Artin's root numbers or global epsilon factors).

\section{\textbf{Boltje's canonical Brauer induction formula and extendible functions}}

\subsection{\textbf{Brauer Canonical induction formula}}

Let $G$ be a finite group. In this subsection, we define a canonical Brauer induction
formula, which is due to Robert Boltje (see \cite{RB}, \cite{RB2}).

A finite-dimensional complex representation $\rho$ of $G$ is called {\bf monomial} if we have the following decomposition:
$$\rho=\rho_1\oplus\rho_2\oplus\cdots\oplus\rho_n,$$
where $\rho_1,\cdots,\rho_n$ are one-dimensional sub-representation of $\rho$ which are called the lines of $\rho$, and 
$\bbC$-linear action of $G$ on $\rho$, such that the group elements permute the lines of $\rho$.

We call a representation $\rho$ is {\bf simple} if this action of $G$ on the lines is transitive.
In this article, we need the following setting 
for $G$:
\begin{itemize}
\item \textbf{$\bbC[G]:=$} the group ring of $G$ over $\bbC$.
 \item \textbf{$R(G):=$} the character ring or equivalently the Grothendieck group of $\mathbb{C}[G]$-modules provided with 
 the tensor product as the multiplication, and unit representation as the unit element.
 \item $\mathcal{M}_{G}:=\{(H,\varphi)\mid H\leq G,  \varphi:H\rightarrow \mathbb{C}^\times,
  \text{linear character of $H$} \}$, the set of pairs consisting a 
 subgroup and a \textbf{linear} character of that subgroup. Group $G$ acts on $\mathcal{M}_G$ by 
  conjugation, i.e.,
  \begin{center}
   $G\times \mathcal{M}_G\rightarrow \mathcal{M}_G$ , $G$ acts from the left, $g\times(H,\varphi):=(H^g,\varphi^g)$, 
   where $H^g:=gHg^{-1}$,
  \end{center}
 and $\varphi^g:=\varphi(g^{-1}hg)$ for all $g\in G$, and $h\in H$. In particular, this covers the whole $G$-conjugacy classes
   of group $H$. Only for the \textbf{normal subgroups} $H$, $H$ will be fixed.
We denote the $G$-orbit of $(H,\varphi)$ 
by 
$[H,\varphi]_G=[H,\varphi](:=\overline{(H,\varphi)}^{G}$ Boltje's notation), and we denote the set of 
$G$-orbit by $\mathcal{M}_{G}/G$.

\textbf{Poset structure}.
We define natural poset (partially ordered set) structures on $\mathcal{M}_G$, and $\mathcal{M}_G/G$ by 
\begin{center}
 $(\Gamma,\psi)\leq(H,\varphi)\Longleftrightarrow \Gamma\leq H \quad\text{and}\quad \psi=\varphi|_\Gamma$, and \\
 $[\Gamma,\psi]\leq[H,\varphi]\Longleftrightarrow(\Gamma,\psi)\leq {^g}(H,\varphi)$ for some $g\in G$.
\end{center}
Infima exists in $\mathcal{M}_G$ but in general not in $\mathcal{M}_G/G$.
\item $R_{+}(G):=$ the Free abelian group whose basis consists of isomorphic classes of irreducible
monomial representations
  of $G$.
 In our setting,  $\mathcal{M}_{G}/G$ is the basis of $R_{+}(G)$. $R_{+}(G)$ is also the 
 Grothendieck group of monomial representations
 of $G$. 
 $R_{+}(G)$
 again can 
be provided with an abelian ring 
structure with the identity element $[G,1_G]$. We can write each element of $R_{+}(G)$ as an integral linear combination of the 
basis elements $[H,\varphi]$. 
Because $G$ acts trivially on the pairs $(G,\varphi), \varphi\in \hat{G}$, we can consider 
$\bbZ[\hat{G}]$ as a subgroup of $R_{+}(G)$.
\item $b_{G}:R_{+}(G)\to R(G)$, $[H,\varphi]\mapsto \mathrm{Ind}_{H}^{G}\varphi$ the induction map from
 equivalence classes of monomial pairs to the monomial representations of $G$. 
 This $b_G$ is a \textbf{surjective} ring homomorphism 
 by the Brauer induction theorem (cf. \cite{RB2}).

\item {\bf Multiplication in $R_{+}(G)$:}

 Natural {\bf multiplication} in $R_+(G)$ corresponds to the tensor product in $R(G)$:
 \begin{equation}\label{eqn 2.3}
  [H_1,\chi_1]\cdot [H_2,\chi_2]=\sum_{(g_1,g_2)\in H_1\backslash G\times H_2\backslash G}
  [H_1^{g_1}\cap H_2^{g_2},\chi_1^{g_1}\chi_2^{g_2}]=
  \sum_{g\in H_1\backslash G/H_2}[H_1^g\cap H_2,\chi_1^g\chi_2],
 \end{equation}
where the brackets on the right side indicate the pairs which will occur, and which have still put together 
into equivalence classes, such that 
$$[H_1^{g_1}\cap H_2^{g_2}, \chi_1^{g_1}\chi_2^{g_2}]=[H_1^{g_1g_2^{-1}}\cap H_2, \chi_1^{g_1g_2^{-1}}\chi_2]=
[H_1\cap H_2^{g_2g_1^{-1}},\chi_1\chi_2^{g_2g_1^{-1}}],$$
because $[.,.]$ denotes $G$-equivalence classes. Therefore, writing also the right side of equation (\ref{eqn 2.3}) 
as a sum of 
equivalence classes, it is then only over a subset of pairs $(g_1,g_2)$ representing the different equivalence classes, 
which means that
$$g_1g_2^{-1}\in H_1\backslash G/H_2\quad or\quad g_2g_1^{-1}\in H_2\backslash G/H_1$$
cover the different double cosets. Then 
$$b_G([H_1,\chi_1]\cdot [H_2,\chi_2])=b_G([H_1,\chi_1])\otimes b_G([H_2,\chi_2])=
Ind_{H_1}^{G}(\chi_1)\otimes Ind_{H_2}^{G}(\chi_2),$$
because 
\begin{align*}
Ind_{H_1}^{G}(\chi_1)\otimes Ind_{H_2}^{G}(\chi_2)=Ind_{H_1}^{G}(\chi_1\otimes Res_{H_1}^{G}(Ind_{H_2}^{G}(\chi_2)))\\
= Ind_{H_1}^{G}\big(\chi_1\otimes\sum_{g\in H_2\backslash G/H_1}Ind_{H_2^{g}\cap H_1}^{H_1}(\chi_2^g|_{H_2^g\cap H_1})\big)\\
=\sum_{g\in H_2\backslash G/H_1} Ind_{H_2^g\cap H_1}^{G}(\chi_1|_{H_2^g\cap H_1}\cdot\chi_2^g|_{H_2^g\cap H_1}).
\end{align*}
Similarly, 
\begin{align*}
Ind_{H_1}^{G}(\chi_1)\otimes Ind_{H_2}^{G}(\chi_2)=Ind_{H_2}^{G}(Res_{H_2}^{G}(Ind_{H_1}^{G}(\chi_1))\otimes\chi_1)\\
=\sum_{g\in H_1\backslash G/H_2} Ind_{H_1^g\cap H_2}^{G}(\chi_1^g|_{H_1^g\cap H_2}\cdot\chi_2|_{H_1^g\cap H_2}).
\end{align*}
\begin{rem}
1. \underline{{\bf Twisting in $R_+(G)$:}} If $H_2=G$, hence $[H_2,\chi_2]=(G,\chi_2)$, then (\ref{eqn 2.3}) comes down to 
$$[H_1,\chi_1]\cdot (G,\chi_2)=[H_1,\chi_1\cdot Res_{H_1}^{G}(\chi_2)]$$
corresponding to 
$$Ind_{H_1}^{G}(\chi_1)\otimes\chi_2=Ind_{H_1}^{G}(\chi_1\otimes Res_{H_1}^{G}(\chi_2)).$$
We call this {\bf twist} of $[H_1,\chi_1]$ by $\chi_2$.

The twisting formula in $R_{+}(G)$ shows that the pairs $[G,\varphi]$ are multiplied with each other as in $\bbZ[\hat{G}]$.\\

2. We also see that $R_{+}(G)$ contains $\bbZ[\hat{G}]$ as a subring, and is therefore a $\bbZ[\hat{G}]$-algebra. 
We also see $b_G([G,\varphi])=Ind_{G}^{G}(\varphi)=\varphi$, then we can consider $b_G$ as a $\bbZ[\hat{G}]$-algebra 
homomorphism. We also have the following $\bbZ[\hat{G}]$-linear ring homomorphism:
\begin{equation}
 \pi_G:R_{+}(G)\to\bbZ[\hat{G}],\quad [H,\varphi]\to\begin{cases}
                                                     \varphi & \text{if $H=G$},\\
                                                     0, & \text{if $H<G$}.
                                                    \end{cases} 
\end{equation}

3. Using commutative multiplication (\ref{eqn 2.3}), the map $b_G$ turns into a {\bf surjective homomorphism of commutative 
rings}, and $Ker(b_G)\subset R_+(G)$ becomes an {\bf ideal}.

\end{rem}

We also need to mention some {\bf functionality properties of $b_G$:}

\item\textbf{Restriction map $\mathrm{Res}_{+H}^{G}$}. If $H\leq G$, we define a restriction map that is also 
a ring homomorphism 
\begin{center}
 $\mathrm{Res}_{+H}^{G}:R_{+}(G)\to R_{+}(H)$
\end{center}
by the double coset formula 
\begin{equation}
 \mathrm{Res}_{+H}^{G}([B,\mu]_{G})=\sum_{s\in H\backslash G/B}[H\cap B^s,\mu^s]_{H}
\end{equation}
where $\mu^s$ means $\mu^s|_{H\cap B^s}$.

\item \textbf{Induction $\mathrm{Ind}\,_{+H}^{G}$}.
For $H\leq G$, we define the map 
$\mathrm{Ind}_{+H}^{G}:R_{+}(H)\to R_{+}(G)$ is an inclusion, i.e.,  
$[\Gamma,\psi]_H\to [\Gamma,\psi]_G$.
This map is well-defined by the following commutative diagram 
 \[
\begin{CD}
R_{+}(H) @>{\mathrm{Ind}_{+H}^{G}}>> R_{+}(G) \\
@V{b_H}VV  @VV{b_G}V \\
R(H) @>>{\mathrm{Ind}_{H}^{G}}> R(G) \\
\end{CD}
\]
From this above diagram, we have $b_{G}\circ\mathrm{Ind}_{+H}^{G}=\mathrm{Ind}_{H}^{G}\circ b_H$.

\item {\bf Inflation $Inf_{+ G/N}^{G}$:} If $f:G\to G/N=:\overline{G}$ is the canonical surjection for a normal 
subgroup $N$ of $G$, we obtain
\begin{equation}
 Inf_{+ H/N}^{H}([H/N,\overline{\varphi}]_{\overline{G}})=[H,\varphi],
\end{equation}
where $N\leq H\leq G$, and $\varphi\in\widehat{H}$ vanishes on $N$. Thus $Inf_{+G/N}^{G}$ maps the basis 
$\mathcal{M}_{\overline{G}}/\overline{G}$ injectively into the basis $\mathcal{M}_G/G$.

\end{itemize}

\begin{prop}\label{Proposition 2.2}
 If $\rho\in R_+(G)$ is any element, then the product with $[H,\chi]\in R_+(G)$ can be rewritten as:
 $$\rho\cdot [H,\chi]=Ind_{+H}^{G}(Res_{+H}^G(\rho)\cdot [H,\chi]),$$
 where $Res_{+H}^{G}(\rho)\cdot [H,\chi]$ is the product in $R_{+}(H)$.
\end{prop}

\begin{proof}
 We may assume that $\rho=[H_1,\chi_1]$ is one of the generators of $R_+(G)$. Then, we obtain:
 $$Res_{+H}^{G}([H_1,\chi_1])=\sum_{g\in H_1\backslash G/H}[H_1^g\cap H,\chi_1^g|_{H_1^g\cap H}]_H,$$
 $$Res_{+H}^{G}([H_1,\chi_1]\cdot [H,\chi])=\sum_{g\in H_1\backslash G/H}[H_1^g\cap H,\chi_1^g\chi|_{H_1^g\cap H}]_H,$$
 where for the second equality, we have used the torsion formula in the case $G=H$. However, then applying $Ind_{+H}^{G}$
 to the second equality, we see from (\ref{eqn 2.3}) that the right side turns into $[H_1,\chi_1]\cdot[H,\chi]$.
\end{proof}

\begin{rem}
 If $\rho=\sum_{[H_i,\chi_i]\in\mathcal{M}_G/G}\alpha_{[H_i,\chi_i]}[H_i,\chi_i]\in R_{+}(G)$, 
 then $\rho\in Ker(b_G)$ means a relation
 $$\sum_{[H_i,\chi_i]\in\mathcal{M}_G/G}\alpha_{[H_i,\chi_i]}Ind_{H_i}^{G}(\chi_i)\cong 0\in R(G).$$
If $G'$ is a subquotient of $G$, which means $G' \overset{u}\twoheadleftarrow G''\le G$, then our functoriality properties fit 
together with a commutative diagram 
$$ \begin{CD}
R_+(G') @>Inf_+>> R_+(G'') @>{Ind_+}>> R_+(G) \\
@Vb_{G'}VV @Vb_{G''}VV @Vb_{G}VV \\
R(G') @>Inf>> R(G'') @>Ind>> R(G)
\end{CD} .$$
Moreover, because $Ker(b_G)$ is always an ideal, and taking into account the torsion operation, we see from the diagram that 
\begin{equation}
 \sigma\in R_+(G')\mapsto Ind_{+G''}^{G}(\chi\cdot Inf_{+G'}^{G''}(\sigma)),\quad \text{where $\chi\in\widehat{G''}$},
\end{equation}
will take $Ker(b_{G'})$ into $Ker(b_G)$. 
This gives an answer to the following question:\\
{\it Why is the kernel of $b_{G'}$ (hence the kernel of $b_G$) important for 
studying extendible functions?}
\end{rem}

\subsection{{\bf The Map $a_G$}}

Now we are in a position to define the section map $a_G$ of $b_G$
from $R(G)$ to $R_{+}(G)$, that is, here we define the function
$$a_G:R(G) \to R_{+}(G)$$
such that $b_G\circ a_G$ is the identity map on $R(G)$.


\underline{{\bf Axioms of $a_G$:}}\\
For a family of maps $a_G:R(G)\to R_{+}(G)$, we consider the following two conditions:
\begin{enumerate}
 \item 
 \[
\begin{CD}
R(G) @>{a_G}>> R_{+}(G) \\
@V{Res_{H}^{G}}VV  @VV{Res_{+H}^{G}}V \\
R(H) @>>{a_H}> R_{+}(H) \\
\end{CD}
\]
This diagram commutes for all subgroups relation $H\leq G$.
\item 
The following diagram is commutative for all groups $G$,
\[
\begin{tikzpicture}
\matrix (m) [matrix of math nodes, row sep=2em, column sep=2em, text height=1.5ex, text depth=0.5ex]
 {   {}    &  R_{+}(G) \\
   R(G)  &  \mathbb{Z}[\widehat{G}] \\
  };
\path[->,font=\scriptsize ]
(m-2-1) edge node[auto] {$\scriptstyle \rho_G$}(m-2-2)
(m-1-2) edge node[auto] {$\scriptstyle \pi_G$} (m-2-2); 
\path[->, font=\scriptsize]
(m-2-1) edge node[auto] {$\scriptstyle a_G$} (m-1-2);
\end{tikzpicture}
\]


Here $\rho_G:R(G) \to \bbZ[\hat{G}]$ is defined by 
$$\rho_G(\chi)=\sum_{\varphi\in\hat{G}}(\varphi,\chi)\cdot\varphi.$$

\end{enumerate}

\begin{thm}[\cite{RB2}, Theorem 2.1]
 There is a unique family of maps $a_G:R(G)\to R_{+}(G)$ satisfying conditions $(1)$, and $(2)$ (from above axioms),
 such that
\begin{equation}\label{eqn 2.4}
 a_{G}(\chi)=\sum_{[H,\varphi]\in \mathcal{M}_{G}/G}\alpha_{[H,\varphi]}
(\chi)\cdot[H,\varphi]\in R_{+}(G).
\end{equation}

\end{thm}
This family has the following properties (for details, cf. \cite{RB2}, pp. 38-40):
\begin{enumerate}
 \item[(a)]{\bf Description of $a_G$:} 
 \begin{enumerate}
  \item[(i)] The coefficients $\alpha_{[H,\varphi]}(\chi)$ are unique  for a given $\chi$.
  \item[(ii)] $a_G$ is the adjoint map of $b_G$ with respect to $[-,-]$, and $(-,-)$.
  \item[(iii)] We have the explicit formula 
  \begin{equation}
   a_G(\chi)=\frac{1}{|G|}\sum_{(H,\varphi)\le (H',\varphi')\in\mathcal{M}_G}|H|\mu_{(H,\varphi),(H',\varphi')}^{\mathcal{M}_G}
   (\varphi',\chi|_{H})\cdot [H,\varphi]
  \end{equation}
where $\mu_{(H,\varphi),(H',\varphi')}^{\mathcal{M}_G}$ is the M\"{o}bius function of $\mathcal{M}_G$.
 \end{enumerate}
\item[(b)] $b_G\circ a_G=id_{R(G)}$ ($a_G$ is a Brauer induction formula for $G$).
\item[(c)] $a_G$ is $\bbZ[\hat{G}]$-linear ($a_G$ commutes with twists by linear characters), and for each $\varphi\in\hat{G}$
we have $a_G(\varphi)=[G,\varphi]$.
\item[(d)] For each $\chi\in R(G)$, we have 
\begin{equation}
 \chi(1)\cdot 1_G=\sum_{[H,\varphi]\in\mathcal{M}_G/G}\alpha_{[H,\varphi]}(\chi) Ind_{H}^{G}(1_H)
\end{equation}
and for each $\chi\in R(G)$ with $\chi(1)=0$, we have 
\begin{equation}
 \chi=\sum_{[H,\varphi]\in\mathcal{M}_G/G}\alpha_{[H,\varphi]}(\chi) Ind_{H}^{G}(\varphi-1_H).
\end{equation}
\end{enumerate}

\begin{rem}
{\bf i).}
In 1947, Richard Brauer (cf. \cite{Br1}) proved that any virtual representation $\rho\in R(G)$ can be expressed as follows:
\begin{equation}\label{eqn Brauer Induction formula}
 \rho=\sum_{i}n_i Ind_{H}^{G}(\varphi_i),\qquad n_i\in\bbZ, H_i\le G, \varphi_i\in\widehat{H_i}.
\end{equation}
In \cite{RB}, Robert Boltje gives an explicit, and {\it canonical} formula for Brauer's induction theorem by algebraic, and 
combinatorial methods. {\it Under certain functorial properties of $a_G$, it can be proved that the expression 
(\ref{eqn 2.4}) is {\bf unique} among all the expressions for $\rho$ as above (\ref{eqn Brauer Induction formula}).
Because of its {\bf uniqueness}, the expression (\ref{eqn 2.4}) is called {\bf canonical}.}\\
{\bf ii).} Here, we also mention the modified canonical Brauer induction formula (cf. Property (h) on p. 40 of \cite{RB2}): for 
all $\rho\in R(G)$ with $\rho(1)=0$, we have 
\begin{equation}\label{eqn modified Canonical BIF}
 \rho=\sum_{[H,\varphi]\in\mathcal{M}_G/G}\alpha_{[H,\varphi]}(\rho)Ind_{H}^{G}(\varphi-1_H),
\end{equation}
where $1_H$ is the trivial character of $H$.\\
{\bf iii).} For all $\rho\in R(G)$, we have the expression (\ref{eqn 2.4}). Here, we should also mention two formulas
regarding the relationship between the coefficients $\alpha_{[H,\varphi]}(\rho)$
and the dimension of $\rho=\rho(1)$, and they are (cf. Property (j) of $a_G$ on p. 40 of \cite{RB2}):
\begin{equation}
 \sum_{[H,\varphi]\in\mathcal{M}_G/G}\alpha_{[H,\varphi]}(\rho)=\rho(1), \quad and 
\end{equation}
\begin{equation}\label{eqn 2.10}
 \sum_{[H,\varphi]\in\mathcal{M}_G/G}[G:H]\cdot \alpha_{[H,\varphi]}(\rho)=\rho(1). 
\end{equation}
\end{rem}

\subsection{Extendible  functions} \label{Subsection extendible function}

Let $G$ be any finite group. We denote $\mathcal{M}_G^v$ as the set of all pairs $(H,\rho)$,
where $H$ is a subgroup of $G$, and $\rho$ is a
virtual representation of $H$.
The group $G$ acts on $\mathcal{M}_G^v$ as follows:
\begin{center}
$(H,\rho)^g=(H^g,\rho^g)$,     $g\in G$,\\
$\rho^g(x)=\rho(gxg^{-1})$,   $x\in H^g:=g^{-1}Hg$.
\end{center}
Furthermore, we denote by $\widehat{H}$ the set of all one-dimensional representations of $H$, and 
by $\mathcal{M}_G$ the subset of $\mathcal{M}_G^v$ of pairs $(H,\chi)$ with 
$\chi\in \widehat{H}$. Here character $\chi$ of $H$ is always a \textbf{linear} 
character, i.e., $\chi:H\to \mathbb{C}^\times$. 

Now, we define a function $\mathcal{F}:\mathcal{M}_G \rightarrow \mathcal{A}$, where $\mathcal{A}$ is a multiplicative abelian 
group with
\begin{equation}\label{eqn 2.2.1}
\mathcal{F}((H,1_H))=1
\end{equation}
, and 
\begin{equation}\label{eqn 2.2.2}
\mathcal{F}((H^g,\chi^g))=\mathcal{F}((H,\chi))
\end{equation}
for all $(H,\chi)$, where $1_H$ denotes the trivial representation of $H$.\\
Here, a function $\mathcal{F}$ on $\mathcal{M}_G$, we mean a function which satisfies Equations (\ref{eqn 2.2.1}) 
and (\ref{eqn 2.2.2}).

A function $\mathcal{F}$ is said to be extendible  if $\mathcal{F}$ can be extended to 
an $\mathcal{A}$-valued 
function
on $\mathcal{M}_G^v$ satisfying: 
\begin{equation}\label{eqn 2.2.3}
 \mathcal{F}((H,\rho_1+\rho_2))=\mathcal{F}((H,\rho_1))\mathcal{F}((H,\rho_2))
\end{equation}
for all $(H,\rho_i)\in \mathcal{M}_G^v,i=1,2$, and if $(H,\rho)\in \mathcal{M}_G^v$ with
$\mathrm{dim}\,\rho=0$, and $\Delta$ is a subgroup of
$G$ containing 
$H$, then
\begin{equation}\label{eqn 2.2.4}
 \mathcal{F}((\Delta,\mathrm{Ind}_{H}^{\Delta}\rho))=\mathcal{F}((H,\rho)),
\end{equation}
where $\mathrm{Ind}_{H}^{\Delta}\rho$ is the virtual representation of $\Delta$ induced from $\rho$.

\begin{rem}[{\bf Weakly Extendible Functions, and Langlands $\lambda$-functions}]

We can always construct a zero-dimensional representation from a finite-dimensional representation. 
Let $\rho$ be a 
representation of $H$ with $\mathrm{dim}(\rho)\neq0$.
We can define a zero-dimensional representation of $H$ by $\rho$, and which is:
  $$\rho_0:=\rho-\mathrm{dim}(\rho)\cdot 1_H.$$ 
  So $\mathrm{dim}\,\rho_0$ is zero, then we use the Equation (\ref{eqn 2.2.4}) for
$\rho_0$, and we obtain
\begin{equation}\label{eqn 2.2.5}
 \mathcal{F}((\Delta,\mathrm{Ind}_{H}^{\Delta}\rho_0))=\mathcal{F}((H,\rho_0)).
 \end{equation}
 Now replace $\rho_0$ by $\rho-\mathrm{dim}\rho\cdot 1_H$ in the above Equation (\ref{eqn 2.2.5}), and we have
 \begin{align*}
   \mathcal{F}((\Delta,\mathrm{Ind}_{H}^{\Delta}(\rho-\mathrm{dim}\rho \cdot 1_H)))
   &=\mathcal{F}((H,\rho-\mathrm{dim}\rho\cdot1_H))\\\implies
   \frac{\mathcal{F}((\Delta,\mathrm{Ind}_{H}^{\Delta}\rho))}
   {\mathcal{F}((\Delta,\mathrm{Ind}_{H}^{\Delta}1_H))^{\mathrm{dim}\rho}}
   &=\frac{\mathcal{F}((H,\rho))}
   {\mathcal{F}((H,1_H))^{\mathrm{dim}\rho}}.
 \end{align*}
Therefore, 
\begin{align}
 \mathcal{F}((\Delta,\mathrm{Ind}_{H}^{\Delta}\rho))\nonumber
 &=\left\{\frac{\mathcal{F}((\Delta,\mathrm{Ind}_{H}^{\Delta}1_H))}{\mathcal{F}((H,1_H))}\right\}^{\mathrm{dim}\rho}
 \cdot\mathcal{F}((H,\rho))\\
 &=\lambda_{H}^{\Delta}(\mathcal{F})^{\mathrm{dim}\rho}\mathcal{F}((H,\rho)), \label{eqn 2.6}
\end{align}
where
\begin{equation}
 \lambda_{H}^{\Delta}(\mathcal{F}):=\frac{\mathcal{F}((\Delta,\mathrm{Ind}_{H}^{\Delta}1_H))}{\mathcal{F}((H,1_H))}.\label{eqn 2.7}
\end{equation}
However, by the definition of $\mathcal{F}$, we have 
$\mathcal{F}((H,1_H))=1$, so we can write 
\begin{equation}
 \lambda_{H}^{\Delta}(\mathcal{F})={\mathcal{F}((\Delta,\mathrm{Ind}_{H}^{\Delta}1_H})).\label{eqn 2.8}
\end{equation}
This $\lambda_{H}^{\Delta}(\mathcal{F})$
is called \textbf{Langlands $\lambda$-function} (or simply $\lambda$-function) which is independent of $\rho$.
A extendible  function $\mathcal{F}$ is called \textbf{strongly} extendible  if it satisfies
equation (\ref{eqn 2.2.3}), and fulfills Equation (\ref{eqn 2.2.4}) for all $(H,\rho)\in \mathcal{M}_G^v$,
and if Equation (\ref{eqn 2.2.4})
is fulfilled
only when $\mathrm{dim}\,\rho=0$, then
$\mathcal{F}$ is called a \textbf{weakly} extendible  function. 
\end{rem}

\begin{lem}[\cite{JT1}, p. 103]\label{Uniqueness of extendible function}
 If extendible functions exist, they are unique.
\end{lem}

\begin{proof}
Let $\mu$ be a function on $\mathcal{M}_G$ and satisfy Equations (\ref{eqn 2.2.1}), (\ref{eqn 2.2.2}).
Let $\mu_1, \mu_2$ be two extendible functions of $\mu$ on $\mathcal{M}_G^v$. Now we have to show $\mu_1=\mu_2$. 
By definition for one dimensional
representation $\mu_1=\mu_2$, and $\mu_i((H,1_H))=1, i=1,2$.
Again because $\mu_1,\mu_1$ are extendible functions of $\mu$, they satisfy equations (\ref{eqn 2.2.3}), and (\ref{eqn 2.2.4}).

Let $H\le G$ be a subgroup of $G$, and $\rho\in R(H)$.
\begin{align*}
 \mu_1((H,\rho))=\mu_1((H,\rho-\dim(\rho)\cdot 1_H))=
 \mu_1((H,\sum_{[U_i,\chi_i]\in\mathcal{M}_H/H}\alpha_{[U_i,\chi_i]}Ind_{U_i}^{H}(\chi_i-1_{U_i})))=\\
 \prod_{[U_i,\chi_i]\in\mathcal{M}_H/H}\mu_1((H,Ind_{U_i}^{H}(\chi_i-1_{U_i})))^{\alpha_{[U_i,\chi_i]}}=
 \prod_{[U_i,\chi_i]\in\mathcal{M}_H/H}\mu_1((U_i,\chi_i-1_{U_i}))^{\alpha_{[U_i,\chi_i]}}\\
 = \prod_{[U_i,\chi_i]\in\mathcal{M}_H/H}\mu_1((U_i,\chi_i))^{\alpha_{[U_i,\chi_i]}}\quad \text{because $\mu_1((U_i,1_{U_i}))=1$}\\
 = \prod_{[U_i,\chi_i]\in\mathcal{M}_H/H}\mu_2((U_i,\chi_i))^{\alpha_{[U_i,\chi_i]}}\quad\text{because $\mu_1=\mu_2$ on $\mathcal{M}_G$}\\
 =\mu_2((H,\rho)).
\end{align*}

This implies that if extendible functions exist, they are unique.
\end{proof}

In the following, we mention three {\bf functoriality properties of the $\lambda$-functions} which are attached to an 
extendible function (for arithmetic proof, cf. \cite{RL} when extendible function is local abelian root number/epsilon
factor). Group theoretically, it is not hard to 
 see its proof. They follow from the definition.
 
 \begin{lem}\label{Lemma 4.1.1}
  Let $H$ be a subgroup of a group $G$, and $\mathcal{F}$ be an extendible function on $\mathcal{M}_G$.
  Then we have the following
properties of $\lambda$-factor.
\begin{enumerate}
 \item If $g\in G$, then $\lambda_{g^{-1}Hg}^{G}(\mathcal{F})=\lambda_{H}^{G}(\mathcal{F})$, where $H\subseteq G$. 
 \item If $H'$ is a subgroup of $H$ then 
 $\lambda_{H'}^{G}(\mathcal{F})=\lambda_{H'}^{H}(\mathcal{F})\lambda_{H}^{G}(\mathcal{F})^{[H:H']}$, 
 where $[H:H']$ is the index
 of $H'$ in $H$.
 \item If $H'$ is a normal subgroup of $G$ contained in $H$, then 
 $\lambda_{H}^{G}(\mathcal{F})=\lambda_{H/{H'}}^{G/{H'}}(\mathcal{F})$. 
\end{enumerate}
 \end{lem}

\begin{rem}
For explicit computation for $\lambda$-functions, see Proposition 2 on p. 124 of \cite{GH84}, Saito's Theorem on p. 508 
of \cite{TS}, Theorems 1.1, and 1.2 on p. 182 of \cite{SAB1}, and Theorem 1.1 of \cite{SABLTRQ}.
If $F$ is a non-Archimedean local field, and $K$ is a wildly ramified quadratic extension of $F$, then computation of 
$\lambda_{\{1\}}^{Gal(K/F)}$ is still {\bf open}. When $F=\bbQ_2$, for explicit computation see pp. 60-64 of 
\cite{SABT}.

\end{rem}

\begin{exm}
In \cite{RL}, Langlands proves that local constants are weakly extendible  functions. This is the main theme of this paper.
In Section 3,  we prove this. The Artin root numbers
(also known as global constants) are strongly extendible  functions (cf. Appendix, Subsection 6.1 of this paper).
\end{exm}

In the following theorem, the {\bf criterion for the extendability} of a function using Langlands $\lambda$-function is shown. 
{\it The existence of a $\lambda$-function with appropriate behavior can be transformed into a criterion for extendibility.}
This is a crucial idea of Langlands' proof of the existence of nonabelian local constants.

\begin{thm}[Theorem 2.1 of \cite{RL}, Lemma 3.2 of \cite{HK}]\label{Lemma 3.2 of Koch}
 Let $G$ be a finite group, and $\mathcal{F}$ be a function on $\mathcal{M}_G$ with values in the abelian group 
 $\mathcal{A}$. Then, $\mathcal{F}$ is extendible  
 to $\mathcal{M}_G^v$ 
 if and only if 
 for all subgroups $H$ of $G$, there is a function 
 $$U\in\mathcal{U}(H)\mapsto \lambda_{U}^{H}(\mathcal{F})\in\mathcal{A}$$
 which is defined on the set $\mathcal{U}(H)$ of subgroups of $H$ such that:\\ 
 \begin{equation}\label{eqn 5.20}
\lambda_{H}^{H}(\mathcal{F})=1,
 \end{equation} 
and if there are $[U,\chi_U]\in \mathcal{M}_H/H$ such that
$$\sum_{[U,\chi_U]\in \mathcal{M}_H/H}\alpha_{[U,\chi_U]}(0)Ind_{U}^{H}(\chi_U)=0,$$
then 
\begin{equation}\label{eqn 5.21}
 \prod_{[U,\chi_U]\in \mathcal{M}_H/H}\mathcal{F}((U,\chi_U))^{\alpha_{[U,\chi_U]}(0)}
 \lambda_{U}^{H}(\mathcal{F})^{\alpha_{[U,\chi_U](0)}}=1.
\end{equation}

\end{thm}

\begin{proof}
 If $\mathcal{F}$ is extendible, then by the definition of Langlands' $\lambda$-function, and the 
 Equations (\ref{eqn 2.2.3}), and (\ref{eqn 2.2.5}), conditions (\ref{eqn 5.20}), and (\ref{eqn 5.21}) hold.
 
 Conversely, suppose that $\lambda_{U}^{H}(\mathcal{F})$ exists
 with (\ref{eqn 5.20}), (\ref{eqn 5.21}). Now for any $H\le G$, we must define a function 
 $\mathcal{F}:\mathcal{M}_H^v\to \mathcal{A}$ such that $\mathcal{F}$ is extendible to $\mathcal{M}_G^v$.
 Before this, we need these two 
 following properties of $\lambda$-functions:\\ 
1. For a fix $H\leq G$, the function $U\in\mathcal{U}(H)\mapsto\lambda_{U}^{H}(\mathcal{F})$ is {\bf uniquely} determined by 
(\ref{eqn 5.20}), (\ref{eqn 5.21}). 
Using Equation (\ref{eqn modified Canonical BIF}), we can write
\begin{equation}
 Ind_{U}^{H}(1_U)-[H:U]\cdot 1_H=\sum_{[U',\varphi_{U'}]\in\mathcal{M}_H/H}\alpha_{[U',\varphi_{U'}]}
 Ind_{U'}^{H}(\varphi_{U'}-1_{U'}).
\end{equation}
This implies that 
$$\lambda_{U}^{H}(\mathcal{F})=\prod_{[U',\varphi_{U'}]\in \mathcal{M}_H/H}
\mathcal{F}((U',\varphi_{U'}))^{\alpha_{[U',\varphi_{U'}]}}.$$
2. If $H\subseteq G'\subseteq G,$ then consider a function on $\mathcal{U}(H)$ as follows:
$$U\in\mathcal{U}(H)\mapsto\lambda_{U}^{G'}(\mathcal{F})\lambda_{H}^{G'}(\mathcal{F})^{-[H:U]}.$$
It turns out to be a function
satisfying Equations (\ref{eqn 5.20}), (\ref{eqn 5.21}). Again for 
$$\sum_{[U_i,\chi_i]\in\mathcal{M}_H/H}\alpha_{[U_i,\chi_i]}Ind_{U_i}^{H}(\chi_i)\cong 0,$$
from Equation (\ref{eqn 2.10})
$$\sum_{[U_i,\chi_i]\in\mathcal{M}_H/H}\alpha_{[U_i,\chi_i]}\cdot [H:U_i]=0,$$ and 
$$\sum_{[U_i,\chi_i]\in\mathcal{M}_H/H}\alpha_{[U_i,\chi_i]}Ind_{U_i}^{G}(\chi_i)\cong 0.$$
Therefore, 
\begin{align*}
 \prod_{[U_i,\chi_i]\in\mathcal{M}_H/H}\mathcal{F}((U_i,\chi_i))^{\alpha_{[U_i,\chi_i]}}
\big(\lambda_{U_i}^{G}(\mathcal{F}){\lambda_{H}^{G}}^{-[H:U_i]}\big)^{\alpha_{[U_i,\chi_i]}}\\
= \prod_{[U_i,\chi_i]\in\mathcal{M}_H/H}\mathcal{F}((U_i,\chi_i))^{\alpha_{[U_i,\chi_i]}}
 (\lambda_{U_i}^{G})^{\alpha_{[U_i,\chi_i]}}=1,
\end{align*}
where the second equality is Condition (\ref{eqn 5.21}) for our original $\lambda$.
Hence, {\bf uniqueness} implies
\begin{equation}\label{eqn 19}
 \lambda_{U}^{H}(\mathcal{F})=\lambda_{U}^{G'}(\mathcal{F})\lambda_{H}^{G'}(\mathcal{F})^{-[H:U]}.
\end{equation}
Now, we come to the proof.

Let $\rho$ be a virtual representation of $H$. Then, we can uniquely write:
\begin{equation}\label{equation 2.19}
 \rho=\sum_{[U',\varphi_{U'}]\in \mathcal{M}_H/H}\alpha_{[U',\varphi_{U'}](\rho)}Ind_{U'}^{H}(\varphi_{U'}).
\end{equation}
Then, from this above equation, we define
$$\mathcal{F}((H,\rho)):=\prod_{[U',\varphi_{U'}]\in\mathcal{M}_H/H}
\mathcal{F}((U',\varphi_{U'}))^{\alpha_[U',\varphi_{U'}](\rho)}{\lambda_{U'}^{H}(\mathcal{F})}^{\alpha_{[U',\varphi_{U'}](\rho)}}.$$
 
Furthermore, it is clear that Equation (\ref{eqn 2.2.3}) is satisfied. It remains to show property (\ref{eqn 2.2.5}).

Let $G'$ be a subgroup of $G$ containing $H$. Then 
$$Ind_{H}^{G'}(\rho)\cong \sum_{[U',\varphi_{U'}]\in\mathcal{M}_H/H}\alpha_{[U',\varphi_{U'}](\rho)}
Ind_{U_i}^{G'}(\varphi_{U'}).$$
Hence, by definition
$$\mathcal{F}((G', Ind_{H}^{G'}(\rho)))=\prod_{[U',\varphi_{U'}]\in\mathcal{M}_H/H}
\mathcal{F}((U',\varphi_{U'}))^{\alpha_{[U',\varphi_{U'}](\rho)}}
{\lambda_{U'}^{G'}(\mathcal{F})}^{\alpha_{[U',\varphi_{U'}](\rho)}}.$$
Therefore, it is sufficient to show 
$$\lambda_{U_i}^{G'}(\mathcal{F})=\lambda_{U_i}^{H}(\mathcal{F})\cdot\lambda_{H}^{G'}(\mathcal{F})^{[H:U_i]}$$
and it follows from equation (\ref{eqn 19}).

This completes the proof of the theorem.
\end{proof}

\begin{rem}
 Let $u:G\twoheadrightarrow G'$ be a {\bf surjective} group homomorphism. Then, we have 
 $$u^*: \mathcal{M}_{G'}\to \mathcal{M}_G,\qquad (H',\chi')\mapsto (u^{-1}(H'),\chi'\circ u),$$
and a function $\mathcal{F}$ on $\mathcal{M}_G$ induces a function $\mathcal{F}_{G'}$ on $\mathcal{M}_{G'}$ via
 $$\mathcal{F}\mapsto\mathcal{F}_{G'}=\mathcal{F}\circ u^*,\quad \mathcal{F}_{G'}((H',\chi')):=
 \mathcal{F}((u^{-1}(H'),\chi'\circ u)).$$
 Alternatively, let $G'\le G$ be a subgroup. Then, we have a natural injection
 $$\iota: \mathcal{M}_{G'}\subset \mathcal{M}_G,$$
and a function $\mathcal{F}$ on $\mathcal{M}_G$ has a well-defined restriction $\mathcal{F}_{G'}=\mathcal{F}\circ\iota$
 down to $\mathcal{M}_{G'}$:
 $$\mathcal{F}_{G'}((H',\chi')):=\mathcal{F}((H',\chi')),\quad\text{for $(H',\chi')\in\mathcal{M}_{G'}$}.$$
\end{rem}

\begin{prop}
 In both cases the maps $u^*$, and $\iota$ naturally extend to {\bf additive maps}:\\
 $u^*:\mathcal{M}_{G'}^v\to\mathcal{M}_G^v$, and $\iota: \mathcal{M}_{G'}^{v}\subset \mathcal{M}_G^v$ resp., and if 
 $\mathcal{F}$ is extendible onto $\mathcal{M}_G^v$ then $\mathcal{F}_{G'}$ extends from $\mathcal{M}_{G'}$ onto 
 $\mathcal{M}_{G'}^v$ using $\mathcal{F}\circ u^*$, and $\mathcal{F}\circ \iota$ resp.
\end{prop}
\begin{proof}
 As to check the property (\ref{eqn 2.2.4}) for $\mathcal{F}_{G'}=\mathcal{F}\circ u^*$, we note that for $u$ we have by 
 definition 
 $$\mathcal{F}_{G'}((H',\rho'))=\mathcal{F}((u^{-1}(H'),\rho'\circ u)),$$
 $$\mathcal{F}_{G'}((\Delta',Ind_{H'}^{\Delta'}(\rho')))=\mathcal{F}((u^{-1}(\Delta'), Ind_{H'}^{\Delta'}(\rho')\circ u)).$$
 We must see that the left sides are equal if $\dim(\rho')=0$. Thus, we turn to the right sides where we may use that 
 $$\dim(\rho'\circ u)=\dim(\rho')=0\quad, and\quad 
 Ind_{H'}^{\Delta'}(\rho')\circ u=Ind_{u^{-1}(H')}^{u^{-1}(\Delta')}(\rho'\circ u),$$
 hence, the right sides are equal because $\mathcal{F}$ is extendible.
\end{proof}


Before giving a list of generating relations for $Ker(b_G)$ in Section 3, we turn to one particular example which will be of 
crucial importance.

Let $G$ be a finite solvable group.
We consider $(H,\chi)\in \mathcal{M}_G$. Furthermore, let $C\leq G$ be any abelian 
normal subgroup. Then, we consider the set of characters
$$S:=S(C,\chi):=\{\mu\in\widehat{C}|\mu|_{H\cap C}=\chi|_{H\cap C}\}\subset \widehat{C}.$$
Here, $H$ acts by conjugation on $C$, and on $\widehat{C}$ which induces an action $\mu\mapsto \mu^h=h^{-1}\mu h$
on the subset $S$.

For $\mu\in S$, we let 
$$H_\mu:=\{h\in H| h\mu h^{-1}=\mu\}$$
be the isotropy group. Furthermore, we put
$$G_\mu'=H_\mu\cdot C,\quad \mu'\in \widehat{G_\mu'}\quad\text{the extension of
$\mu\in S\subset\widehat{C}$ by means of $\chi$},$$
more precisely: 
$$\mu'(hc)=\chi(h)\mu(c)\quad\text{for $h\in H_\mu, c\in C$}.$$
Further, let $T\subset S$ be a system of representatives for the orbit $S/H$ (action of $H$ on $S$ from the right).

\begin{rem}
In the following lemma, we see how the monomial representation $Ind_{H}^{G}(\chi)$ splits over the system of 
representative $T$ for $S/H$. In the remaining parts of Sections 2, and 3, this lemma will be used many times
because we have to deal with monomial representations.
For every element $(H,\chi)\in\mathcal{M}_G$, we will have $T, S, G_\mu'$, therefore the following lemma is important.
\end{rem}
\begin{lem}[Lemma 15.1 of \cite{RL}, Lemma 2.1 of \cite{HK}]\label{Lemma Koch 2.1}
 With notation, and assumptions as we have fixed above, we have 
 $$Ind_{H}^{G}(\chi)=\sum_{\mu\in T}Ind_{G_\mu'}^{G}(\mu').$$
\end{lem}
\begin{proof}
To prove the above assertion, we only need to show 
$$[G:H]=\sum_{\mu\in T}[G:G_\mu'],$$
because here $\chi$, and $\mu'$ are characters of $H$, and $G_\mu'$, respectively.\\

Here $G_{\mu}^{'}=H_\mu C$, then we can see that $G_\mu'\subseteq HC\subseteq G$.
For any subgroup $K$ of $G$ containg $H$, we can write
$$Ind_{H}^{G}\chi\equiv Ind_{K}^{G}(Ind_{H}^{K}\chi).$$
Since here we have $G_\mu'\subseteq HC\subseteq G$, 
 we can rewrite our assertion as 
 $$Ind_{HC}^{G}(Ind_{H}^{HC}(\chi))=Ind_{HC}^{G}(\sum_{\mu\in T}Ind_{G_\mu'}^{HC}(\mu')).$$
 This above relation implies that 
 it is sufficient to prove the lemma in the case $G=HC$.
 
 Let $V$ be a complex-vector space of all functions $f:G\to\bbC$ with
 $$f(hx)=\chi(h)f(x),\quad \text{for all}\quad h\in H, x\in G.$$
 G acts on $V$ as follows
 $$(gf)(x):=f(xg)\quad\text{for all $g\in G$}.$$
Because of $G=HC$ we have $|S|=[C:H\cap C]=[HC:H]=[G:H]$. For every $\mu\in S$, we define a function 
$f_\mu:G\to\bbC$ as follows:
$$\mu\in S\mapsto f_\mu,\quad f_\mu(hc):=\chi(h)\mu(c),\quad \text{for all}\quad h\in H, \, c\in C.$$
It can be proved (cf. Lemma 2.2 on p. 14 of \cite{HKEWZ}) that the set $\{f_\mu\}_{\mu\in S}$ is a basis of the space $V$. 

On the other hand, for every $\mu\in S$, we can construct the subspace denoted by $V_\mu$ of $V$ generated by the 
set 
$\{f_{\mu^g}:=f_{g\mu g^{-1}}\}\quad \text{for all $g\in G/{G_{\mu}'}$}.$
Furthermore, we can show (cf. \cite{HKEWZ}) that the subspace $V_\mu$ is G-isomorphic to $Ind_{G_\mu'}^{G}(\mu')$.

And the dimension of $Ind_{G_\mu'}^{G}(\mu')$ is $[G:G_\mu']=[H:H_\mu]$ if we assume $G=HC$, because 
$$G_\mu'\cap H=H_\mu C\cap H=H_\mu(C\cap H)=H_\mu,$$
hence 
$$G/G_\mu'=HC/H_\mu C\overset{\sim}\leftarrow H/H_\mu.$$
Hence, $[G:H]=\sum_{\mu\in T}[G:G_\mu']$. This completes the proof.
\end{proof}

\begin{rem}
 If $hc=1$, which means $h=c^{-1}\in H\cap C$, then $\mu\in S$ means $\mu(c^{-1})=\chi(h)$, and therefore 
 $$\mu'(hc)=\chi(h)\mu(c)=\mu(c^{-1})\mu(c)=1, \quad\text{as it should be.}$$
 If in particular $H\cap C=\{1\}$, and $H\cdot C=G$, then we have $S=\widehat{C}$, and $T=\widehat{C}/H$, hence 
 $$Ind_{H}^{G}(\chi)=\sum_{\mu\in \widehat{C}/H}Ind_{H_\mu\cdot C}^{G}(\chi|_{H_\mu}\cdot\mu),\quad 
 \widetilde{\chi}\otimes Ind_{H}^{G}(1_H)=
 \widetilde{\chi}\otimes\sum_{\mu\in \widehat{C}/H}Ind_{H_\mu\cdot C}^{G}(1_{H_\mu}\cdot\mu),$$
 where the second equality is only a reformation of the first one because using $G/C\cong H$, we can extend $\chi$ to the character
 $\widetilde{\chi}$ of $G$ which is trivial on $C$. According to Subsection 8.2 of \cite{JPSLR}, this is the decomposition of 
 $Ind_{H}^{G}(\chi)$ into irreducible components.
\end{rem}

\begin{lem}[Lemma 2.4 of \cite{HK}]\label{Lemma 2.4 of Koch}
 Let $C$ be a abelian normal subgroup of $G$, let $T$ be a set of representatives of $\widehat{C}/G$, and let $G_\mu$ be 
 the isotropy group of $\mu\in T$ in $G$. Furthermore, let $\{H_i:i\in I\}$ be a family of subgroup of $G$ containing $C$,
and let $\chi_i$ be a character of $H_i$ such that $\chi_i|_C\in T$. If we assume a relation
 \begin{equation}\label{eqn 2.27}
  \sum_{i\in I}n_i Ind_{H_i}^{G}(\chi_i)=0,\qquad \sum_{i\in I}n_i[H_i,\chi_i]_G\in Ker(b_G),
 \end{equation}
then for any fixed $\mu\in T$, we have 
$$\sum_{\chi_i;\chi_i|_C=\mu}n_i Ind_{H_i}^{G_\mu}(\chi_i)=0,
\qquad\sum_{\chi_i,\chi_i|_C=\mu}n_i[H_i,\chi_i]_{G_\mu}\in Ker(b_{G_\mu}).$$
\end{lem}

\begin{proof}
Because $C\subseteq H_i\subseteq G$, and $C$ is normal we have $C\backslash G/H_i=G/H_i$, and
$$Ind_{H_i}^{G}(\chi_i)|_C=\sum_{g\in G/H_i}\chi_i^g|_C.$$
Thus, we meet here only one orbit of $\widehat{C}/G$ hence, we meet only one $\mu$. Therefore, if $Ind_{H_i}^{G}(\chi_i)$
and $Ind_{H_j}^{G}(\chi_j)$ lead to different representatives in $T$ then their restrictions to $C$
are disjoint, and therefore the representations itself must be disjoint. Thus, if we fix one $\mu$, the
assumption (\ref{eqn 2.27}) will imply that the corresponding partial sum must vanish:
\begin{equation}\label{eqn 2.28}
 \sum_{\chi_i;\chi_i|_C=\mu}n_i Ind_{H_i}^{G}(\chi_i)=0,
\end{equation}
hence the assumption implies separate relations (\ref{eqn 2.28}) for each $\mu\in T$.
Furthermore, $\chi_i|_C=\mu$ implies $H_i\le G_\mu$, and:
$$Ind_{H_i}^{G}(\chi_i)|_{G_\mu}=\sum_{H_i\backslash G/G_\mu}Ind_{H_i^s\cap G_\mu}^{G_\mu}(\chi_i^s).$$
The direct sum on the right contains the term $Ind_{H_i}^{G_\mu}(\chi_)$ for $s\in H_i G_\mu$, and additional
components for $H_isG_\mu\ne H_i G_\mu$. And restricting further to $C\subseteq H_i^s\cap G_\mu$ we see that those
other components are disjoint from $Ind_{H_i}^{G_\mu}(\chi_i)$ because $s\not\in G_\mu$. Therefore, the separate
relations (\ref{eqn 2.28}) will imply the assertion of our Lemma for any fixed $\mu\in T$.
 
\end{proof}

\begin{rem}[{\bf Definition of Lambda-factors}]
To prove Theorem \ref{Theorem 1.2}, we will use Theorem \ref{Lemma 3.2 of Koch}.
 To use Theorem \ref{Lemma 3.2 of Koch}, we have to define 
 $$U\in \mathcal{U}(H)\mapsto \lambda_{U}^{H}(\mathcal{F})$$
 with properties
 (\ref{eqn 5.20}), (\ref{eqn 5.21}), where $U$, and $H$ are two arbitrary subgroups of $G$ with $U\subset H$.
 Note that any definition of $U\mapsto\lambda_{U}^H(\mathcal{F})$ such that 
 $\lambda_{U^h}^{H}(\mathcal{F})=\lambda_{U}^{H}(\mathcal{F})$ allows a linear extension of $\mathcal{F}$ onto the free 
 abelian group $R_+(H)$ using
 $$\mathcal{F}([U,\chi]_H):=\mathcal{F}((U,\chi))\cdot\lambda_{U}^{H}(\mathcal{F}),$$ because 
 $\mathcal{F}$ has the property (\ref{eqn 2.2.2}), but we need this extension to be trivial on $Ker(b_H)$.

 We fix a nontrivial minimal abelian normal subgroup $C\le G$ which exists because $G$ is solvable. We proceed by induction
 over $|G|$. If $|G|=1$, we put $\lambda_G^G=1$. In Lemma \ref{Lemma Koch 2.1}, put $\chi\equiv 1$ the trivial character of 
 $U$. Thus $S=\widehat{C/U\cap C}$, and $\mu':G_\mu'=U_\mu C\to\bbC^\times$ is the extension of $\mu\in S\subset\widehat{C}$
 by $1$, thus we can write
 \begin{equation}\label{eqn 2.29}
  Ind_U^G(1)\cong\sum_{\mu\in S/U}Ind_{G_\mu'}^{G}(\mu').
 \end{equation}
If $\mathcal{F}$ is extendible, then this implies 
\begin{equation}\label{eqn 2.30}
 \lambda_{U}^{G}(\mathcal{F})=\prod_{\mu\in T}\mathcal{F}((G_\mu',\mu'))\lambda_{G_\mu'}^{G}(\mathcal{F}).
\end{equation}
Now we turn (\ref{eqn 2.30}) into a definition of $\lambda_{U}^{G}(\mathcal{F})$. Indeed, we can define 
\begin{equation}\label{eqn 2.31}
 \lambda_{G_\mu'}^{G}(\mathcal{F}):=\lambda_{G_\mu'/C}^{G/C}(\mathcal{F}),
\end{equation}
because $G_\mu'=U_\mu C\supseteq C$, which brings us down to groups of smaller order where we may use the induction hypothesis.
And $\mathcal{F}((G_\mu',\mu'))$ is defined anyway because $\mu'$ is one-dimensional.
Therefore, we may consider Equation (\ref{eqn 2.30})
as a definition of $\lambda_U^G(\mathcal{F})$.

Using a fixed nontrivial minimal abelian normal subgroup $C\le G$, and the isomorphism (\ref{eqn 2.29}) we define:
$$\lambda_U^G(\mathcal{F})=\mathcal{F}(Ind_U^G(1)):
=\prod_{\mu\in T}\mathcal{F}((G_\mu',\mu'))\lambda_{G_\mu'/C}^{G/C}(\mathcal{F}),$$
 where $\lambda$ on the right side has already been defined by the induction hypothesis.

Note here that the choice of $C$ is part of the definition because a priori we do not know that $\lambda_U^G(\mathcal{F})$
will be unique. In the particular case where $U=\{1\}$ we have $G_\mu'=C$, and our
Definitions (\ref{eqn 2.30}), (\ref{eqn 2.31})
turn into 
$$\lambda_{\{1\}}^{G}(\mathcal{F}):=\prod_{\mu\in\widehat{C}}\big(\mathcal{F}((C,\mu))\lambda_{\{1\}}^{G/C}(\mathcal{F})\big)=
\big(\prod_{\mu\in\widehat{C}}\mathcal{F}((C,\mu))\big)\cdot \lambda_{\{1\}}^{G/C}(\mathcal{F})^{|C|}.$$
If $G$ is of prime order, then the only choice is $C=G$, and therefore the definition comes down to 
$$\lambda_{\{1\}}^G(\mathcal{F}):=\prod_{\mu\in\widehat{G}}\mathcal{F}((G,\mu)).$$
In the particular case, where $U\supseteq C$, the definition turns into 
$$\lambda_U^G(\mathcal{F}):=\lambda_{U/C}^{G/C}(\mathcal{F}).$$
This includes our original Definition 
(\ref{eqn 2.31}).



\end{rem}

As to the definition of $\lambda_U^H(\mathcal{F})$ for proper subgroups $H<G$, we may assume that is defined by the induction 
hypothesis. Now, we have the following lemma.
\begin{lem}\label{Lemma 2.18}
 (i). The above definition of $\lambda_{U}^{G}(\mathcal{F})$ is independent of the choice of the set $T$ of representatives 
 for $S/U$.\\
 (ii). $$\lambda_{U^g}^{G}(\mathcal{F})=\lambda_{U}^{G}(\mathcal{F})\quad\text{for all $g\in G$}.$$
 (iii). If $N\le U$ is a non-trivial normal subgroup of $G$ which sits in $U$, then 
 $\lambda_U^G(\mathcal{F})=\lambda_{U/N}^{G/N}(\mathcal{F})$ where the right side is already given by the induction hypothesis.
 In particular,
 $$\lambda_U^G(\mathcal{F})=\lambda_{e}^{G/U}(\mathcal{F}), \quad\text{if $U$ is a normal subgroup of $G$}.$$
 
\end{lem}
\begin{proof}
(i). If $u\in U,$ then 
 $$G_{\mu^u}'=u^{-1}G_\mu' u,\quad (\mu^u)'=(\mu')^u.$$
 Hence, by the condition (\ref{eqn 2.2.2}) of the extendible function,
 $$\mathcal{F}((G_{\mu^u}',(\mu^u)'))=\mathcal{F}((u^{-1}G_\mu' u, (\mu')^u))=\mathcal{F}((G_\mu',\mu')),$$
and 
 $$\lambda_{G_{\mu^u}'}^{G}(\mathcal{F})=\lambda_{u^{-1}G_\mu' u}^{G}(\mathcal{F})=\lambda_{G_\mu'}^{G}(\mathcal{F})$$
 by Lemma \ref{Lemma 4.1.1}(1).\\
(ii). We must compare  Equation (\ref{eqn 2.29}) with the corresponding formula for $U^g=g^{-1}Ug$, instead of $U$. Then, we obtain 
$$S^g:=\widehat{C/U^g\cap C},\, \mu\in S/U\mapsto \mu^g\in S^g/U^g$$
and $Stab_{U^g}(\mu^g)=g^{-1}Stab_{U}(\mu)g$, hence $G_{\mu^g}'=Stab_{U^g}(\mu^g)\cdot C=g^{-1}G_\mu' g$, and 
$(\mu^g)'=(\mu')^g$ which implies 
$$\lambda_{U^g}^{G}(\mathcal{F})=\prod_{\mu}\mathcal{F}((g^{-1}G_\mu' g, (\mu')^g))\lambda_{(G_\mu')^g/C}^{G/C}(\mathcal{F})
=\lambda_U^G(\mathcal{F})$$
because $\mathcal{F}$ has the property (\ref{eqn 2.2.2}), and $(G_\mu')^g/C=(G_\mu'/C)^g$ such that 
$\lambda_{(G_\mu')^g/C}^{G/C}(\mathcal{F})=\lambda_{G_\mu'/C}^{G/C}(\mathcal{F})$ by the induction hypothesis.\\
(iii). If $U$ contains a normal subgroup $N$ of $G$, then using the notation of (\ref{eqn 2.30}), we have $N\subseteq U_\mu:$
indeed $\mu$ is a character of $C/U\cap C$, and $[N,C]\subseteq N\cap C\subseteq U\cap C$, hence 
$$nxn^{-1}x^{-1}\in N\cap C\subseteq U\cap C\quad \text{for $x\in C, n\in N$},$$
and therefore 
$$\frac{\mu^n}{\mu}(x)=\mu^n\mu^{-1}(x)=\mu^n(x)\mu^{-1}(x)=\mu(nxn^{-1}x^{-1})=1.$$
Thus, we obtain $NC\le U_\mu C=G_\mu'$, and $\mu'$ (which is the extension of $\mu$ by $1$) is trivial on $N$, hence 
(\ref{eqn 2.29}) is rewritten as
$$Ind_U^G(1)=Ind_{U/N}^{G/N}(1)=\sum_{\mu\in T}Ind_{G_\mu'/NC}^{G/NC}(\mu'),\,\lambda_{U/N}^{G/N}(\mathcal{F})=
\prod_{\mu\in T}\mathcal{F}((G_\mu'/N,\mu'))\cdot \lambda_{G_\mu'/N}^{G/N}(\mathcal{F}),$$
because for the group $G/N$ of a smaller order, we have all properties available by the induction hypothesis. Finally:
$$\lambda_{G_\mu'/N}^{G/N}(\mathcal{F})=\lambda_{G_\mu'/NC}^{G/NC}(\mathcal{F})=\lambda_{G_\mu'/C}^{G/C}(\mathcal{F}),$$
because all groups occurring here are of a smaller order. Therefore, we obtain 
$$\lambda_{U/N}^{G/N}(\mathcal{F})=\lambda_{U}^{G}(\mathcal{F})$$
as we have defined it using (\ref{eqn 2.31}).

 This completes the proof.
 
 \end{proof}

\begin{rem}
 In general, definition (\ref{eqn 2.30}) is true if $U$ contains a normal subgroup $N$ of $G$.
 In this case $N\subseteq G_\mu'$ because 
 $N\cap C\subseteq U\cap C$. Hence, $nxn^{-1}x^{-1}=:[n,x]\in U\cap C$, and $\mu([n,x])=1$ for $x\in C, n\in N$. Therefore,
 (\ref{eqn 2.30}) is the inflation of the corresponding equation for $G/N, U/N, CN/N$, which is valid by the 
 induction assumption.
\end{rem}

The definition of $\lambda_U^G(\mathcal{F})$ has been completed, and we have proved that the definition (\ref{eqn 2.30})
implies
Lemma (\ref{Lemma 4.1.1}(1)), and (\ref{Lemma 4.1.1}(3)). Because of Lemma \ref{Lemma 2.18}(ii), we have now a well 
defined $\bbZ$-linear map 
$$R_+(G)\to\mathcal{A},\quad [U,\chi]_G\mapsto\mathcal{F}([U,\chi]_G):=\mathcal{F}((U,\chi))\lambda_U^G(\mathcal{F}).$$
Now, we need to check that the definition of $\lambda_{U}^{G}(\mathcal{F})$ satisfies property Lemma \ref{Lemma 4.1.1} (2).
And we check it in the following lemma.

\begin{lem}[Lemma 3.4 of \cite{HK}]\label{Lemma 3.4 of Koch}
 If $G'$ is a subgroup of $G$ containing $U$, then 
 $$\lambda_{U}^{G}(\mathcal{F})=\lambda_{U}^{G'}(\mathcal{F})\lambda_{G'}^{G}(\mathcal{F})^{[G':U]}.$$
\end{lem}

\begin{proof}
By the definition $\lambda_{G}^{G}(\mathcal{F})=1$. Therefore, if $G'=G$, then the assertion is trivial.

Let $G'\ne G$, then 
$\lambda_{U}^{G'}(\mathcal{F}_{G'})$ is defined by the induction assumption.
Using $G_\mu'=U_\mu C\subseteq G'C\subseteq G$, and the induction hypothesis, we obtain:
\begin{align*}
 \lambda_{G_\mu'}^{G}(\mathcal{F}):=\lambda_{G_\mu'/C}^{G/C}(\mathcal{F})\\
 =\lambda_{G_\mu'/C}^{G'C/C}(\mathcal{F}_{G'C})\cdot\lambda_{G'C/C}^{G/C}(\mathcal{F})^{[G'C:G_\mu']}\\
 =\lambda_{G_\mu'}^{G'C}(\mathcal{F}_{G'C})\cdot\lambda_{G'C}^{G}(\mathcal{F})^{[G'C:G_\mu']}.
\end{align*}
This together with the relation 
$$[G'C:G']=\sum_{\mu\in T}[G'C:G_\mu'] \quad(\text{because $Ind_{G'}^{G'C}(1)\cong\sum_{\mu\in T}Ind_{G_\mu'}^{G'C}(\mu')$})$$
implies
\begin{align}\label{eqn 5.25}
 \lambda_{G'}^{G}(\mathcal{F}):=\prod_{\mu\in T}\mathcal{F}((G_\mu',\mu'))\lambda_{G_\mu'}^{G}(\mathcal{F})\\\nonumber
 =\lambda_{G'C}^{G}(\mathcal{F})^{[G'C:G']}
 \prod_{\mu\in T}\mathcal{F}((G_\mu',\mu'))\lambda_{G_\mu'}^{G'C}(\mathcal{F}_{G'C})\\\nonumber
 =\lambda_{G'}^{G'C}(\mathcal{F}_{G'C})\lambda_{G'C}^{G}(\mathcal{F})^{[G'C:G']}.
\end{align}
Similarly, we can see
\begin{equation}\label{eqn 5.26}
 \lambda_{U}^{G}(\mathcal{F})=\lambda_{U}^{G'C}(\mathcal{F}_{G'C})\lambda_{G'C}^{G}(\mathcal{F})^{[G'C:U]}.
\end{equation}
Furthermore, if $G'C\ne G$, we have 
$$\lambda_{U}^{G'C}(\mathcal{F}_{G'C})=\lambda_{U}^{G'}(\mathcal{F}_{G'})\lambda_{G'}^{G'C}(\mathcal{F}_{G'C})^{[G':U]}$$
by the induction assumption. We multiply the last equation by $\lambda_{G'C}^{G}(\mathcal{F})^{[G'C:U]}$:
$$\lambda_{U}^{G'C}(\mathcal{F}_{G'C})\lambda_{G'C}^{G}(\mathcal{F})^{[G'C:U]}
=\lambda_{U}^{G'}(\mathcal{F}_{G'})\lambda_{G'}^{G'C}(\mathcal{F}_{G'C})^{[G':U]}\lambda_{G'C}^{G}(\mathcal{F})^{[G'C:U]}.$$
This together with Equation (\ref{eqn 5.25}), and Equation (\ref{eqn 5.26}), proves the lemma in the {\bf case $G\ne G'C.$}

To complete the proof, we are left with the case: {\bf when $G=G'C$.} To prove this, we need the following lemma. 

Now, suppose that $G=G'C$. Then $G'\cap C$ is a normal subgroup of $G$, and it is contained in $C$.
Because $C$ has been chosen as minimal, abelian, normal in $G$ the only possibilities are $G'\supseteq C$, hence $G=G'C=G$
(which is trivial case; see above) or $G'\cap C=\{1\}$ hence $G=G'\cdot C$ is semidirect product. This is the case 
we are going to proceed with, and before that we need the following lemma.

\begin{lem}[Lemma 3.5 of \cite{HK}]\label{Lemma 3.5 of Koch}
For $U\le G'$, let $\Omega= UC\le G=G'C$ semidirect. Let $\chi:\Omega\to \bbC^\times$ be a character of $\Omega$.
Then, there are subgroups $\Omega_i$ of $\Omega$ 
 with $C\subseteq\Omega_i$, and characters $\mu_i\in\widehat{\Omega_i}$ such that 
 \begin{equation}\label{eqn 5.27}
  Ind_{U}^{G}(\chi_U)\cong\sum_{i}Ind_{\Omega_i}^{G}(\mu_i\chi)
 \end{equation}
\begin{equation}\label{eqn 5.28}
 \mathcal{F}((U,\chi_U))\lambda_{U}^{G}(\mathcal{F})=\prod_{i}\mathcal{F}((\Omega_i,\mu_i\chi))\lambda_{\Omega_i}^{G}(\mathcal{F}).
\end{equation}
\end{lem}
\begin{proof}
 The existence of (\ref{eqn 5.27}) is proved in Lemma \ref{Lemma Koch 2.1}. It remains to show that Equation (\ref{eqn 5.27})
 implies Equation (\ref{eqn 5.28}). \\
 If $U\ne G'$, equivalently $G\ne \Omega$, the induction assumption applies to
 $$Ind_{U}^{\Omega}(\chi_U)\cong\sum_{i}Ind_{\Omega_i}^{\Omega}(\mu_i\chi)$$
 hence 
\begin{equation}\label{eqn 5.29}
 \mathcal{F}((U,\chi_U))\lambda_{U}^{\Omega}(\mathcal{F})=
 \prod_{i}\mathcal{F}((\Omega_i,\mu_i\chi))\lambda_{\Omega_i}^{\Omega}(\mathcal{F}).
\end{equation}
On the other hand, we have already proved Lemma \ref{Lemma 3.4 of Koch} if $G'C\ne G$, and therefore we may use it for 
$G'=\Omega$ because $G'C=G'\ne G$, hence we have:
\begin{equation}\label{eqn 2.37}
 \lambda_{U}^{\Omega}(\mathcal{F})\lambda_{\Omega}^{G}(\mathcal{F})^{[\Omega:U]}=\lambda_{U}^{G}(\mathcal{F}),\qquad
\lambda_{\Omega_i}^{\Omega}(\mathcal{F})\lambda_{\Omega}^{G}(\mathcal{F})^{[\Omega:\Omega_i]}=
\lambda_{\Omega_i}^{G}(\mathcal{F}).
\end{equation}
And the identity $[\Omega:U]=|C|=\sum_{i}[\Omega:\Omega_i]$ where the index stands for the different orbits of 
$U\backslash\widehat{C}$ yields
\begin{equation}\label{eqn 2.38}
\lambda_{\Omega}^{G}(\mathcal{F})^{[\Omega:U]}=\prod_{i}\lambda_{\Omega}^{G}(\mathcal{F})^{[\Omega:\Omega_i]}.
\end{equation}
Now, multiplying Equation (\ref{eqn 5.29}) by Equation (\ref{eqn 2.38}), then using Equalities (\ref{eqn 2.37}), we obtain
(\ref{eqn 5.28}).

Now, assume $U=G'$, equivalently $\Omega=G$, hence $\chi$ is a character of $G=UC$, and in (\ref{eqn 5.27}), we have now
$$Ind_{U}^{G}(\chi_U)=\chi\otimes Ind_{U}^{G}(1_U),$$
where $1_U$ is the trivial character of $U$. Therefore, Equation (\ref{eqn 5.27}) is now the same as 
Equation (\ref{eqn 2.29}) tensored by the 
character $\chi$ of $G$. Because $U\cap C=\{1\}$ we have now $S=\widehat{C}$, $T=\widehat{C}/U$.

If $U$ contains a nontrivial subgroup $N$ which is normal in $G$, then as we have seen in the proof of Lemma \ref{Lemma 2.18},
 Equation (\ref{eqn 2.29}) rewrites as 
$$Ind_{U/N}^{G/N}(1_{U/N})=\sum_{\mu\in T}Ind_{G_\mu'/N}^{G/N}(\mu')$$
and therefore, Equation (\ref{eqn 5.27}) turns into 
$$\chi\otimes Ind_{U/N}^{G/N}(1_{U/N})=\sum_{\mu\in T}\chi\otimes Ind_{G_\mu'/N}^{G/N}(\mu'),
\quad\text{for $\mu\in\widehat{G}$}.$$
Because $\chi$ is trivial on $[G,G]$, and the other tensor factors are trivial on $N$, we are dealing here with 
a representation
of $G/(N\cap[G,G])$. Thus, for $N\cap [G,G]\ne \{ 1\}$, $G/(N\cap[G,G])$ is a group of smaller order, and we can deduce 
Equation (\ref{eqn 5.28}) by the induction hypothesis.

On the other hand, if $N\cap [G,G]=\{1\}$, then $[N,G]\subseteq N\cap[G,G]$ implies that $N$ must be contained in the center
$Z=Z(G)$, hence $N\le Z(G)\cap U\ne\{1\}$.
Then, $N=Z(G)\cap U$ is a nontrivial subgroup of $U$ which is normal in $G$; and then by the induction hypothesis for 
$G/N$,
we have Equation (\ref{eqn 5.28}).\\
The other way around: If $Z(G)\cap U=\{1\}$, and $N\le U$ is normal in $G$, then 
$N\cap[G,G]\ne \{1\}$, and we come down to the induction hypothesis.

To complete the proof of Lemma \ref{Lemma 3.5 of Koch}, we have to check the case: $G=U\ltimes C$, and $U$ does not contain
subgroups $N$ which are normal in $G$.
Again, 
\begin{equation}\label{eqn 5.30}
  Ind_{U}^{G}(\chi_U)\cong\sum_{\mu\in T} Ind_{G_\mu'}^{G}(\mu'\chi)
 \end{equation}
 implies
\begin{equation}\label{eqn 5.31}
 \mathcal{F}((U,\chi_U))\prod_{\mu\in T}\mathcal{F}((G_\mu',\mu'))=\prod_{\mu\in T}\mathcal{F}((G_\mu',\mu'\chi))
\end{equation}
because we are assuming $\mathcal{F}$ has the property (\ref{eqn 1.3}). On the other hand, by definition in 
Equation (\ref{eqn 2.30})
we have
\begin{equation}\label{eqn 5.32}
 \lambda_{U}^{G}(\mathcal{F})=\prod_{\mu\in T}\mathcal{F}((G_\mu',\mu'))\lambda_{G_\mu'}^{G}(\mathcal{F}),
\end{equation}
where the last factor is defined using $\lambda_{G_\mu'/C}^{G/C}(\mathcal{F})$. Now, multiply Equation (\ref{eqn 5.32})
by $\mathcal{F}((U,\chi_U))$, and then using Equation (\ref{eqn 5.31}), we obtain
\begin{equation}\label{eqn 5.33}
 \mathcal{F}((U,\chi_U))\lambda_{U}^{G}(\mathcal{F})=
 \prod_{\mu\in T}\mathcal{F}((G_\mu',\mu'\chi))\lambda_{G_\mu'}^{G}(\mathcal{F}).
\end{equation}
This ends the proof of Lemma \ref{Lemma 3.5 of Koch} because Equations (\ref{eqn 5.30}), and (\ref{eqn 5.33}) are the desired
Equations (\ref{eqn 5.27}), and (\ref{eqn 5.28}), respectively.
\end{proof}

Now we come back to the proof of Lemma \ref{Lemma 3.4 of Koch}, where we are left to show that when 
$U\subseteq G'\subseteq G=G'C$, we must have 
$$\lambda_{U}^{G}(\mathcal{F})=\lambda_{U}^{G'}(\mathcal{F})\lambda_{G'}^{G}(\mathcal{F})^{[G':U]}.$$
From the proof of Theorem \ref{Lemma 3.2 of Koch},
it is sufficient to show that 
$$U\in\mathcal{U}(G')\mapsto \lambda_{U}^{G'}:=\lambda_{U}^{G}(\mathcal{F})\lambda_{G'}^{G}(\mathcal{F})^{-[G':U]}$$ 
satisfies the conditions (\ref{eqn 5.20}), 
(\ref{eqn 5.21}) for $H=G'$. The uniqueness of such a function shows then 
$$\lambda_{U}^{G}(\mathcal{F})\lambda_{G'}^{G}(\mathcal{F})^{-[G':U]}=\lambda_{U}^{G'}(\mathcal{F}).$$

Using Boltje's canonical Brauer induction formula, we can write $0\in R(G')$ as follows
\begin{equation}\label{eqn 2.43}
 \sum_{[U_i,\chi_i]\in \mathcal{M}_{G'}/G'}\alpha_{[U_i,\chi_i]}Ind_{U_i}^{G'}(\chi_i)=0.
\end{equation}
Now, we must show that Equation (\ref{eqn 2.43}) implies
\begin{equation}\label{eqn 2.44}
 \prod_{[U_i,\chi_i]\in \mathcal{M}_{G'}/G'}\{\mathcal{F}((U_i,\chi_i))\lambda_{U_i}^{G'}\}^{\alpha_{[U_i,\chi_i]}}=1,\qquad
 \text{where $\lambda_{U_i}^{G'}=\lambda_{U_i}^{G}(\mathcal{F})\lambda_{G'}^{G}(\mathcal{F})^{-[G':U]}$}.
\end{equation}
First, from Equation (\ref{eqn 2.43}) we see:
$$\sum_{[U_i,\chi_i]\in \mathcal{M}_{G'}/G'}\alpha_{[U_i,\chi_i]}[G':U_i]=0,\quad
\text{hence}\quad \prod_{[U_i,\chi_i]\in \mathcal{M}_{G'}/G'}\lambda_{G'}^{G}(\mathcal{F})^{-\alpha_{[U_i,\chi_i]}[G':U_i]}=1.$$
Thus, to verify Equation (\ref{eqn 2.44}), it is enough to see that Equation (\ref{eqn 2.43}) will imply
\begin{equation}\label{eqn 2.45}
 \prod_{[U_i,\chi_i]\in \mathcal{M}_{G'}/G'}\{\mathcal{F}((U_i,\chi_i))
 \lambda_{U_i}^{G}(\mathcal{F})\}^{\alpha_{[U_i,\chi_i]}}=1.
\end{equation}
For each pair $(U_i,\chi_i)$, we may consider the subgroup $\Omega_i=U_i\cdot C\le G=G'\cdot C$, and the character $\chi_i$
of $U_iC$ which is trivial on $C$, and to this situation apply Lemma \ref{Lemma 3.5 of Koch}. Then we obtain relations 
of types (\ref{eqn 5.27}), and (\ref{eqn 5.28}), respectively, which we will write as:
\begin{equation}\label{eqn 2.46}
Ind_{U_i}^{G}(\chi_i)\cong \sum_{j=1}^{r_i}Ind_{U_{ij}}^{G}(\mu_{ij}'\chi_i), 
\end{equation}
where $\{\mu_{ij}\}_{j}\subset \widehat{C}$ are representatives of the cosets $U_i\backslash\widehat{C}$, and where 
$U_{ij}:=Stab_{U_i}(\mu_{ij})C$, and $\chi_i\in\widehat{U_i}$ is understood to be trivial on $C$, and $\mu_{ij}'$
extends $\mu_{ij}$ hence is non-trivial on $C$. And Equation (\ref{eqn 5.28}) applied to $U=U_i$ reads as follows:
\begin{equation}\label{eqn 2.47}
 \mathcal{F}((U_i,\chi_i))\lambda_{U_i}^{G}(\mathcal{F})=\prod_{j=1}^{r_i}\mathcal{F}((U_{ij},\mu_{ij}'\chi_i))
 \lambda_{U_{ij}}^{G}(\mathcal{F}).
\end{equation}
Implementing this into equation (\ref{eqn 2.45}), the assertion becomes 
\begin{equation}\label{eqn 2.48}
   \prod_{[U_i,\chi_i]\in \mathcal{M}_{G'}/G'}\{ \prod_{j=1}^{r_i}\mathcal{F}((U_{ij},\mu_{ij}'\chi_i))
 \lambda_{U_{ij}}^{G}(\mathcal{F})\}^{\alpha_{[U_i,\chi_i]}}=1.
\end{equation}
On the other hand, applying $Ind_{G'}^{G}$ to equation (\ref{eqn 2.43}), and using equation (\ref{eqn 2.46}), we obtain:
\begin{equation}\label{eqn 2.49}
 \sum_{[U_i,\chi_i]\in \mathcal{M}_{G'}/G'}\alpha_{[U_i,\chi_i]}(\sum_{j=1}^{r_i}Ind_{U_{ij}}^{G}(\mu_{ij}'\chi_i))\cong 0.
\end{equation}
Again here all $U_{ij}\ge C$, and therefore $\lambda_{U_{ij}}^{G}(\mathcal{F})=\lambda_{U_{ij}/C}^{G/C}(\mathcal{F})$;
and $\mathcal{F}((U_{ij},\mu_{ij}'\chi_i))$ is always well defined because $\mu_{ij}'\chi_i$ is $1$-dimensional.
Therefore,  we can write
$$\mathcal{F}((U_{ij},\mu_{ij}'\chi_i))\cdot\lambda_{U_{ij}/C}^{G/C}(\mathcal{F})=
\mathcal{F}((U_{ij},\mu_{ij}'\chi_i))\cdot\lambda_{U_{ij}}^{G}(\mathcal{F})=\mathcal{F}(Ind_{U_{ij}}^{G}(\mu_{ij}'\chi_i)).$$
Furthermore, here $\mathcal{F}$ is an extendible function. Applying $\mathcal{F}$ to Equation (\ref{eqn 2.49}), we obtain
the above assertion (\ref{eqn 2.48}).

This completes the proof.

\end{proof}

\section{{\bf Existence proof of non-abelian local constants}}

\begin{rem}[{\bf Why Kernel of $b_G$?}]\label{Remark on Kernel of bG}

In this subsection, we first discuss why we should study the kernel $Ker(b_G)$, and the reasons are as follows:
\begin{enumerate}
\item It can be completely described the generating 
set of $Ker(b_A)$ for $A:=B/K$ with $G\supset B\supset K$, and $K$ normal in $G$, 
then inductively we will have all information about the generating set of $Ker(b_G)$. Similar to any arithmetic
invariants (e.g., $L$-functions, $\gamma$-factors etc.), to study the properties of the local constants, we have to know 
how they behave under induction, and inflation. 
\item Again, $A\twoheadleftarrow B\leq G$, and $\sigma\in Ker(b_A)\subset R_+(A)$
implies 
$$\rho=Ind_{+B}^{G}(\chi\cdot Inf_{+A}^{B}(\sigma))\in Ker(b_G)\subset R_+(G)\quad\text{for all $\chi\in\widehat{B}$}.$$

For this, we have to repeat how {\it inflation}, and {\it induction} are defined for the rings $R_+$. 
If $[H,\psi]\in R_+(A)$, and 
$A\twoheadleftarrow B$ then we consider $H_B$ the full preimage of $H$ in $B$, and $\psi_H$ the lift of $\psi$ from $H$ 
to $H_B$. Then, we have 
$$Inf_{+A}^{B}([H,\psi])=[H_B,\psi_B]\in R_+(B),$$
and this inflation map is compatible with the usual inflation $Inf_A^B: R(A)\to R(B)$. Therefore, it takes $Ker(b_A)$ to 
$Ker(b_B)$. 

Next, if $(H_B,\psi_B)\in\mathcal{M}_B$, and $\chi\in\widehat{B}$, and put 
$$\chi\cdot (H_B,\psi_B)=(H_B,Res_{H_B}^{B}(\chi)\cdot\psi_B),$$
where $Res_{H_B}^{B}$ is the restriction from $B$ to the subgroup $H_B$. 
Under $b_B: R_+(B)\to R(B)$ this is compatible with the 
usual $\chi$-twist for virtual representations of $B$; therefore again the $\chi$-twist on $R_+(B)$ takes 
$Ker(b_B)$ into itself.

As to the induction, the map $Ind_{+B}^{G}: R_+(B)\to R_+(G)$, it is induced by the identity 
$$[H_B,\psi_B]\mapsto [H_B,\psi_B]$$
due to the fact that any subgroup $H_B$ of $B$ may also be considered as a subgroup of $G$. This map is compatible with
the usual induction $Ind_B^G: R(B)\to R(G)$, and therefore again the map $Ind_{+B}^{G}$ takes $Ker(b_B)$ to $Ker(b_G)$.

 \item 
For any 
$$\rho\in Ker(b_G)\subset R_{+}(G),$$
we always have $\deg(\rho)=\deg(b_G(\rho))=\deg(0)=0.$
This implies that for any $\rho\in Ker(b_G)$, $b_G(\rho)$ is a virtual representation of dimension {\bf zero}. Defining local 
constants for dimension zero representations is simple because it remains invariant under induction. Because the
{\it strong extendability} implies {\it weak extendability}, proving {\it the existence of non-abelian local constants, 
is equivalent to prove that the extendible is strong on the virtual representations of dimension zero.}
\end{enumerate}

\end{rem}

\begin{dfn}[{\bf Type I, Type II, and Type III}]\label{Definition of type I, II, III}
An element $\rho\in Ker(b_G)$ is called of {\bf type I, type II, and type III}, if there is a subquotient 
$$A\twoheadleftarrow B\leq G,\quad \text{such that $\rho=Ind_{+B}^G(\chi\cdot Inf_{+A}^B(\sigma))$},$$
where $\chi\in \widehat{B}, \sigma\in Ker(b_A)$, with the following requirements on $A$, and $\sigma$:\\
{\bf Type I.} Let $\ell$ be a prime, and $A$ be the cyclic group of order $\ell$, then denoting $e:=\{1\}<A$ we have 
$$\sigma=[e, 1_e]-\sum_{\mu\in \widehat{A}}[A,\mu]\in Ker(b_A)\subset R_+(A),$$
because 
$$ Ind_{e}^{A}(1_e)=\sum_{\mu\in\widehat{A}}\mu.$$
{\bf Type II.} $A$ is a central extension of an abelian group of type $\bbZ/\ell\bbZ\times\bbZ/\ell\bbZ$ where again $\ell$
is a prime. Thus, we have $Z$ central in $A$, and $A/Z=X_1\times X_2$ is abelian, where $X_i\cong \bbZ/\ell\bbZ$, hence $A$ is 
two-step nilpotent, and the commutator in $A$ induces a bilinear alternating map:
$$[.,.]: A/Z\wedge A/Z\to Z,\quad \bar{a_1}\wedge\bar{a_2}\mapsto [a_1,a_2],$$
which must be trivial if we restrict to the cyclic subgroups $X_i\subset A/Z$. We assume that $A$ is non-abelian which means 
the commutator map is non-trivial. Then obviously  
\begin{itemize}
 \item the subgroups $H_i<A$ which are the full preimages of $X_i\subset A/Z$ must be maximal abelian subgroups such that 
 $H_1\cap H_2=Z$, and
 \item the commutator subgroup is $[A,A]=[H_1,H_2]\cong \bbZ/\ell\bbZ$.
\end{itemize}
We fix characters $\chi_i$ of $H_i$ which agree on $Z:\quad \chi_1|_Z=\chi_2|_Z$, and which are nontrivial on the commutator
subgroup $[A,A]$. Then from Lemma \ref{Lemma Koch 2.1} with $G, H, C$ replaced by $A, H_1, H_2$ we obtain 
$$Ind_{H_1}^{A}(\chi_1)\cong Ind_{H_2}^{A}(\chi_2)\in R(A)$$
and we take 
$$\sigma:=[H_1,\chi_1]-[H_2,\chi_2]\in Ker(b_A)$$
as our relation of {\bf type II}.

{\bf Type III}.  Based on Lemma \ref{Lemma Koch 2.1}, we now come
to the {\bf third type} of requirements on $A$, and $\sigma$:\\  $A=H\ltimes C$ is a semidirect product,
where $C\ne \{1\}$ is an abelian normal subgroup contained in 
all nontrivial abelian normal subgroups of $A$.( In particular this implies that $C\ne \{1\}$ is minimal.) Making use of 
Lemma \ref{Lemma Koch 2.1} with $G=A$, and $(H,1_H)\in\mathcal{M}_A$, we have 
$$Ind_{H}^{A}(1_H)=\sum_{\mu\in T}Ind_{A_\mu'}^{A}(\mu'),$$
and therefore 
$$\sigma:=[H,1_H]-\sum_{\mu\in T}[A_\mu',\mu']\in Ker(b_A)\subset  R_{+}(A).$$

\end{dfn}

\begin{rem}
For the three relations, we have indeed $b_A(\sigma)=0$, as follows from Lemma \ref{Lemma Koch 2.1}: 
For {\bf type I} 
 set $H=\{e\}, C=A$, for {\bf type II} set $H=H_1, C=H_2$, for {\bf type III} set $H=H, C=C$.
 
 We denote the $\bbZ$-module generated by relations of {\bf type I} by $\mathcal{R}(I,G)$, for {\bf type I}, and 
 {\bf type II} by 
 $\mathcal{R}(II, G)$, and by {\bf type I, type II}, and {\bf type III} by $\mathcal{R}(III,G)$. This means 
 $\mathcal{R}(-,G)\subseteq Ker(b_G)\subset R_+(G)$, and $\mathcal{R}(-,G)$ consists of all elements $\rho\in R_+(G)$ which 
 are given as
 $$\rho=\sum_{i=1}^{r}n_i\cdot\rho_i=\sum_{i=1}^{r}n_i\cdot Ind_{+B_i}^{G}(\chi_i\cdot Inf_{+A_i}^{B_i}(\sigma_i)),
 \quad n_i\in\bbZ,$$
 where $A_i=B_i/K_i$ is a subquotient of $G$, $\chi_i\in\widehat{B_i}$, and $\sigma_i\in Ker(b_{A_i})$ as specified above 
 in the {\bf type I, type II, and type III}.
 
 We also can see that $\mathcal{R}(-,G)$ is stable under torsion because 
 $$(G,\chi_G)\cdot Ind_{+B}^{G}(\chi_B\cdot Inf_{+A}^{B}(\sigma))= Ind_{+B}^{G}(Res_B^G(\chi_G)\chi_B\cdot Inf_{+A}^{B}(\sigma)).$$



\end{rem}

{\Large {\bf Proof of Theorem \ref{Theorem 1.1}:}}

Now we are ready to prove Theorem \ref{Theorem 1.1}. And we prove into three cases, and they are:\\
{\bf Case-1: When $G$ is an abelian group.}
\begin{lem}[Lemma 2.3 of \cite{HK}]\label{Lemma 2.3 of Koch}
 If $G$ is abelian, then $Ker(b_G)=\mathcal{R}(I,G).$
\end{lem}
\begin{proof}
If $\rho\in \mathcal{R}(I,G)$, then from Remark \ref{Remark on Kernel of bG}(2) it is simple to see that 
$\rho\in Ker(b_G)$. Now to complete the proof, we have to show that $Ker(b_G)\subset\mathcal{R}(I,G)$.






If $G$ is commutative then the action by conjugation is always trivial, hence $[H,\chi]=(H,\chi)$, and therefore 
$$R_+(G)=\text{free abelian group over $\mathcal{M}_G$}=\bigoplus_{H\le G}R(H)$$
where the second equality holds because all irreducible representations of $H$ are one-dimensional.
Now assume that
$$\rho=\sum_{[H_i,\chi_i]\in\mathcal{M}_G/G}\alpha_{[H_i,\chi_i]}\cdot (H_i,\chi_i)\in Ker(b_G)\subset R_+(G).$$
Thus we have:
$$\sum_{[H_i,\chi_i]\in\mathcal{M}_G/G}\alpha_{[H_i,\chi_i]} Ind_{H_i}^{G}(\chi_i)\cong 0\in R(G).$$
Moreover, because we are in the abelian case: $Ind_{H_i}^{G}(\chi_i)=\sum_j\widetilde{\chi_{ij}}$, 
where $\widetilde{\chi_{ij}}$ are the characters of $G$ which extend $\chi_i$. Therefore, our relation rewrites as
$$\sum_{[H_i,\chi_i]\in\mathcal{M}_G/G}\alpha_{[H_i,\chi_i]}(\sum_j\widetilde{\chi_{ij}})\cong 0\in R(G).$$
However, $R(G)$ naturally embeds into $R_{+}(G)$, and therefore this is the same as
$$\sum_{[H_i,\chi_i]\in\mathcal{M}_G/G}\alpha_{[H_i,\chi_i]}(\sum_j(G,\widetilde{\chi_{ij}}))\cong 0\in R_+(G),$$
such that our original $\rho\in Ker(b_G)$ rewrites as
$$\rho=\sum_{[H_i,\chi_i]\in\mathcal{M}_G/G}\alpha_{[H_i,\chi_i]}\cdot (H_i,\chi_i)=
\sum_{[H_i,\chi_i]\in\mathcal{M}_G/G}\alpha_{[H_i,\chi_i]}\cdot \big((H_i,\chi_i)-\sum_j(G,\widetilde{\chi_{ij}})\big),$$
hence $Ker(b_G)$ is generated by the elements 
$$ (H,\chi)-\sum_{\tilde{\chi}}(G,\tilde{\chi})$$ 
where $\tilde{\chi}\in\widehat{G}$ runs over the 
extensions of $\chi\in\widehat{H}$. Finally, to prove our Lemma we are left to show that such elements
can always be rewritten as
\begin{equation}\label{eqn ?}
 (H,\chi)-\sum_{\tilde{\chi}}(G,\tilde{\chi})=\sum_i Ind_{+B_i}^{G}(\chi_i\cdot Inf_{+A_i}^{B_i}(\sigma_i)),
\end{equation}
where $A_i$
is a subquotient of $G$, and $\sigma_i\in Ker(A_i)$ is of {\bf type I}.

If $H$ is of index $[G:H]=p$ a prime then we may take $B=G, A=G/H$, and (\ref{eqn ?}) rewrites
$$(H,\chi)-\sum_{\tilde{\chi}}(G,\tilde{\chi})=
\widetilde{\chi}\cdot Inf_{+A}^{G}\big((e,1_e)-\sum_{\mu\in\widehat{A}}(A,\mu)\big),$$
where $\widetilde{\chi}$ is one fixed extension of $\chi$, hence all other extensions are given as 
$\widetilde{\chi}\cdot\mu$ for
characters $\mu$ of $A=G/H$.

If $[G:H]$ is not a prime then we argue by induction on $[G:H]$ where we choose a
subgroup $H'$ such that $H<H'<G$, and $[H':H]=p$ a prime. Then, we have:
\begin{equation}\label{eqn ??}
 (H,\chi)-\sum_{\widetilde{\chi}}(G,\widetilde{\chi})=\big((H,\chi)-\sum_{\chi'}(H',\chi')\big)+\sum_{\chi'}\big((H',\chi')-
 \sum_{\widetilde{\chi}|\chi'}(G,\widetilde{\chi})\big),
\end{equation}
where $\chi'$ runs over the extensions of $\chi$ onto $H'$.
Now the first term on right rewrites as
$$(H,\chi)-\sum_{\chi'}(H',\chi')=Ind_{+H'}^{G}\big(\chi_0'\cdot Inf_{+A}^{H'}((e,1_e)-\sum_{\mu\in\widetilde{A}}(A,\mu))\big),$$
where $\chi_0'\in\widehat{H'}$ is a fixed $\chi'$, and where $A=H'/H, e=H/H$, and $H'$ plays the role of $B<G$.
Then under $u:H'\twoheadrightarrow A=H'/H$ we have $u^{-1}(e)=H, u^{-1}(A)=H'$
such that 
$$\chi_0'\cdot Inf_{+A}^{H'}((e,1_e)-\sum_{\mu\in\widetilde{A}}(A,\mu))=(H,\chi)-\sum_{\chi'}(H',\chi')\in R_+(H'),$$
, and $Ind_{+H'}^{G}: R_+(H')\to R_+(G)$ is the identity.
We are left with the other terms on the right of (\ref{eqn ??}), and here we can use the induction
hypothesis because $[G:H']<[G:H]$.


\end{proof}


{\bf Case-2: When $G$ is a nilpotent group.}
\begin{thm}[Langlands, Section 18 of \cite{RL}]\label{Theorem 2.2 of Koch}
 If $G$ is nilpotent, then $Ker(b_G)=\mathcal{R}(II,G)$.
\end{thm}
\begin{proof}
 From Remark \ref{Remark on Kernel of bG}(2), it can be seen that $\mathcal{R}(II,G)\subset Ker(b_G)$. Now we have to prove
 $Ker(b_G)\subset \mathcal{R}(II,G)$. When $G$ is abelian, in Lemma \ref{Lemma 2.3 of Koch} it is proved. So to complete 
 the proof, we have to prove this when $G$ is nonabelian nilpotent group.

 Because abelian groups are nilpotent, and in Lemma \ref{Lemma 2.3 of Koch} it is studied, here we consider $G$ as a nonabelian
 nilpotent group.
 Let $Z$ be the center of $G$. We prove Theorem \ref{Theorem 2.2 of Koch} by induction on $|G/Z|$. When $|G/Z|=1$, that is,
 $G$ is abelian, and it is proved above. Let $|G/Z|\ne 1$.
 Take an arbitrary element
\begin{equation}\label{eqn 5.11}
\rho=\sum_{[H_i,\chi_i]\in\mathcal{M}_G/G}\alpha_{[H_i,\chi_i]}[H_i,\chi_i]\in Ker(b_G).
\end{equation}
Now we have to show that $\rho\in \mathcal{R}(II,G)$.

We can express $b_G([H_i,\chi_i])=Ind_{H_i}^{G}(\chi_i)$ as follows:
$$Ind_{H_i}^{G}(\chi_i)=\sum_{\chi_i'|_{H_i}=\chi_i}Ind_{H_iZ}^{G}(\chi_i')+Ind_{H_iZ}^{G}(Ind_{H_i}^{H_iZ}(\chi_i)-
\sum_{\chi_i'|_{H_i}=\chi_i}\chi_i'),$$
and 
$$Ind_{H_i}^{H_iZ}(\chi_i)-\sum_{\chi_i'|_{H_i}=\chi_i}\chi_i'=b_{H_iZ}(\rho_1),$$
where 
$$\rho_1:=[H_i,\chi_i]-\sum_{\chi_i'|_{H_i}=\chi_i}[H_iZ,\chi_i'].$$
Again $\rho_1$ can be expressed as follows:
$$\rho_1=Inf_{+H_iZ/Ker(\chi_i)}^{H_iZ}([H_i/Ker(\chi_i),\chi_i]-\sum_{\chi_i'|_{H_i}=
\chi_i}[H_iZ/Ker(\chi_i),\chi_i'])$$
and by Lemma \ref{Lemma 2.3 of Koch}
$$[H_i/Ker(\chi_i),\chi_i]-\sum_{\chi_i'|_{H_i}=
\chi_i}[H_iZ/Ker(\chi_i),\chi_i']\in \mathcal{R}(I,H_iZ/Ker(\chi_i)), $$
hence $\rho_1\in \mathcal{R}(I, H_iZ)$. This shows that when $Z\not\subset H_i$, $\rho\in\mathcal{R}(II,G)$.\\
{\bf To complete the proof, we have to show this when $Z\subseteq H_i$.}

Let $C$ be a abelian normal subgroup of $G$ containing $Z$ such that $[C:Z]=\ell$ is a prime, and $C/Z$ is in the center of 
$G/Z$. With the notation of Lemma \ref{Lemma Koch 2.1} for $H=H_i, \chi=\chi_i$, and $G=G$ we have 
$$[H_i,\chi_i]=\sum_{\mu\in T}[G_\mu',\mu']+ Ind_{+H_iC}^{G}(\rho_2)$$
with 
$$\rho_2=[H_i,\chi_i]_{H_iC}-\sum_{\mu\in T}[G_\mu',\mu']_{H_iC}.$$
Therefore, a relation (\ref{eqn 5.11}) with $C\subseteq H_i$ are in $\mathcal{R}(II,G)$. 

Lemma \ref{Lemma 2.4 of Koch} shows that a relation of type (\ref{eqn 5.11}) with $C\subseteq H_i$ is a linear combination 
of relations induced from relations of the form 
$$\sum_{[H_i,\chi_i]\in \mathcal{M}_{G_\mu}/G_\mu}\alpha_{[H_i,\chi_i]}[H_i,\chi_i]\in Ker(b_{G_\mu})$$
with $H_i\supseteq C$, and $\chi_i|_C=\mu$, where $G_\mu$ denotes the isotropy group of $\mu\in\widehat{C}$ in $G$.

Furthermore
$$\sum_{[H_i,\chi_i]\in \mathcal{M}_{G_\mu}/G_\mu }\alpha_{[H_i,\chi_i]}[H_i,\chi_i]=Inf_{+G_\mu/Ker(\mu)}^{G_\mu}(\rho_3)$$
with 
$$\rho_3=\sum u_i[H_i/Ker(\mu),\chi_i]\in Ker(\varphi_{G_\mu/Ker(\mu)}).$$
Because $C/Ker(\mu)$ is central in $G_\mu/Ker(\mu)$ by definition of $G_\mu$, we can apply the induction assumption to $\rho_3$.

It remains to show $\rho_2\in\mathcal{R}(II,H_iC)$. We omit the index $i$. If $HC\ne G$, this follows from the induction 
assumption applied to $HC$. Hence we assume $HC=G$. Then $H$ is a normal subgroup of $G$ because $Z\subseteq H$, and $C/Z$ is 
central in $G/Z$.
We set 
$$X:=\cap_{g\in G}Ker(\chi)^g.$$
Then $\rho_2$ is the inflation of the corresponding relation in $Ker(b_{G/X})$. Hence it is sufficient to prove: \\
{\it if $X=\{1\}$, the relation $\rho_2$ is of {\bf type II}.}

$H$ is abelian because $[H,H]\subseteq X=\{1\}$, $Z$ is cyclic because $Ker(\chi|_Z)\subseteq X=\{1\}$. Because $C/Z$ is central 
in $G/Z$, the commutator $[c,h]:=c^{-1}h^{-1}ch$, $c\in C, h\in H$, lies in $Z$, and depends bi-additively on $c$, and $h$. We 
define a bilinear form $u:C/Z\otimes H/Z\to Z$ by 
$$u(\bar{c},\bar{h})=[c,h].$$
Let $\bar{c}$ be a generator of $C/Z$. Because $Z$ is the center of $G$, and $G\ne Z$, the map 
$u(\bar{c},):H/Z\to Z$ is injective, and image is not trivial. Since
$$[c,h^\ell]=[c,h]^\ell=[c^\ell,h]=1\quad\text{for $h\in H$},$$
we have $[H:Z]=[C:Z]=\ell$. Because $H$ is not central, and $X=\{1\}$, $\chi$ is distinct from some of its conjugates, hence 
$\chi([h,h])\ne \{1\}$, and we are in the situation of {\bf type II}.

\end{proof}

{\bf Case-3: When $G$ is a solvable group.}
\begin{thm}[Theorem 2.6 of \cite{HK}]\label{Theorem 2.6 of Koch}
 If $G$ is solvable, $Ker(b_G)=\mathcal{R}(III,G).$
\end{thm}
\begin{proof}
As before from Remark \ref{Remark on Kernel of bG}(2), it is simple to see that $\mathcal{R}(III,G)\subset Ker(b_G)$. 
Now to complete this proof, we have to show the inclusion: $Ker(b_G)\subset \mathcal{R}(III,G)$. When $G$
is nilpotent, we have proved this above. So we only need to show that for non-nilpotent solvable group $G$, we have 
$Ker(b_G)\subset\mathcal{R}(III,G)$.

We will prove this by induction on the order of $G$. If $|G|=1$, then the assertion is trivial. Now let $H$ be a proper subgroup
of $G$, and for $H$ the assertion is true, that is,
$$Ker(b_H)=\mathcal{R}(III,H)\quad\text{for all proper subgroups $H$ of $G$}.$$

{\bf Step 1:} 
By using this above assumption, we can see that {\it $\mathcal{R}(III,G)$ is an ideal of the ring $R_{+}(G)$}.

Let $\rho\in\mathcal{R}(III,G)$, and $[H,\chi]\in R_{+}(G)$. By definition (see torsion in $R_{+}(G)$ on p. 5),
if $H=G$, then $\rho\cdot [G,\chi]$ belongs 
 to $\mathcal{R}(III,G)$. 

If $H\ne G$, we can see that for $\rho\in\mathcal{R}(III,G)$, hence $\rho\in Ker(b_G)$ we have 
$$b_H(Res_{+H}^{G}(\rho))=Res_{H}^{G}(b_G(\rho))=Res_{H}^{G}(0)=0.$$
This implies $Res_{+H}^{G}(\rho)\in Ker(b_H)=\mathcal{R}(III,H)$, hence 
$Res_{+H}^{G}(\rho)\cdot[H,\chi]\in\mathcal{R}(III,H)$.
 
Now from Proposition \ref{Proposition 2.2} we can write
 $$\rho\cdot[H,\chi]=Ind_{+H}^{G}(Res_{+H}^{G}(\rho)\cdot[H,\chi])\in \mathcal{R}(III,G).$$
%
 This implies that $\mathcal{R}(III,G)$ is an ideal of $R_{+}(G)$.

Let $Z$ be the center of $G$.\\
{\bf Step 2:} Suppose that $Z$ be nontrivial. 
We use the Brauer Induction theorem for $G/Z$ (see \cite{JPSLR} Section 10.5, Theorems 19, and 20). There are {\bf nilpotent} 
subgroups $H_i\supseteq Z$, and characters $\chi_i$ of $H_i/Z$ such that 
\begin{equation}
 1_{G/Z}=\sum_{i}n_i\cdot Ind_{H_i/Z}^{G/Z}(\chi_i)\quad\text{with $n_i\in\bbZ$}.
\end{equation}
Here $G$, and $G/Z$ are not nilpotent, but the $H_i/Z$ are nilpotent, hence they are proper subgroups of $G/Z$
which ensures that we obtain a nontrivial relation.

Because we do induction over the order of $G$, we may assume 
$$Ker(b_{G/Z})=\mathcal{R}(III,G/Z),$$
and by inflation we obtain
\begin{equation}
 \sigma:=[G,1_G]-\sum_{i}n_i\cdot [H_i,\chi_i]\in\mathcal{R}(III,G).
\end{equation}
Given $\rho\in Ker(b_G)$, by using Proposition \ref{Proposition 2.2} we have 
\begin{equation}
 \rho\sigma=\rho-\sum_{i}n_i\cdot (\rho\cdot[H_i,\chi_i])=
 \rho-\sum_{i}n_i Ind_{+H_i}^{G}(Res_{+H_i}^{G}(\rho)\cdot[H_i,\chi_i]).
\end{equation}
However, by assumption $\mathcal{R}(III,G)$ is an ideal in $R_{+}(G)$, hence $\rho\sigma\in \mathcal{R}(III,G)$, and we have seen 
$Res_{+H_i}^{G}(\rho)\in Ker(b_{H_i})=\mathcal{R}(III,H_i)$ because $H_i$ are proper subgroups of $G$. Therefore, it implies
\begin{equation}
 \rho=\rho\sigma+\sum_{i} n_i Ind_{+H_i}^{G}(Res_{+H_i}^{G}(\rho)\cdot[H_i,\chi_i])\in\mathcal{R}(III,G)
\end{equation}
which proves our assertion in the case $Z\ne \{1\}$.

 {\bf Step 3:} When $Z=\{1\}$.\\
Because $G$ is a solvable group,  let $C$ be a minimal abelian normal subgroup of $G$, which exists because $G$ is solvable.
We proceed 
as in the proof of Theorem \ref{Theorem 2.2 of Koch}. Every relation $\rho\in Ker(b_G)$ is a linear combination of the 
form $\rho_2$ (of the proof of Theorem \ref{Theorem 2.2 of Koch}) for $HC$, where $H$ is a subgroup of $G$
which does not contain $C$, and $\chi$ a character of $H$, and of relations induced from relations of the form 
\begin{equation}\label{eqn 5.12}
 \sum_{[H_i,\chi_i]\in \mathcal{M}_{G_\mu}/G_\mu}\alpha_{[H_i,\chi_i]}[H_i,\chi_i]\in Ker(b_{G_\mu})
\end{equation}
with $ G_\mu\supseteq H_i \supseteq C$, and $\chi_i|_C=\mu\in\widehat{C}$.
Since
$[G_\mu:Ker(\mu)]<|G|$, we can apply the induction assumption to (\ref{eqn 5.12}). It remains to consider a relation
of the form $\rho_2$.

If $HC\ne G$ we can apply the induction assumption. If $HC=G$, then $H\cap C$ is normal in $H$
and $C$, hence normal in $G$. The minimality of $C$ implies $H\cap C=\{1\}$. Thus we have $G=H\ltimes C$ is semidirect 
and using $G/C\cong H$ we can extend $\chi\in\widehat{H}$ by setting 
$\chi(c)=1$ for $c\in C$ to a character $\widetilde{\chi}$ of $G$. Then our relation is the torsion of 
\begin{equation}\label{eqn 5.13}
 [H,1_H]-\sum_{\mu\in T}[G_\mu',\mu']
\end{equation}
with $[G,\widetilde{\chi}]$, where $\mu'$ denotes the extension of $\mu$ to $G_\mu'$ which is trivial on $G_\mu'\cap H$.

Suppose there is a nontrivial normal subgroup $H_1$ of $H$ lying in the centralizer of $C$. Then $H_1$ is normal in $G$
because $G=HC$. Then $[H_1,C]\subseteq H_1\cap C\subset H\cap C=\{1\}$. Therefore,  the group $H_1$ is contained in $G_\mu'$
for all $\mu$ hence (\ref{eqn 5.13}) is the inflation of 
$$[H/H_1,1_{H/H_1}]-\sum_{\mu\in T}[G_\mu'/H_1,\mu],$$
which is contained in $\mathcal{R}(III,G/H_1)$ by induction assumption, and then the inflation, and its torsion by
$[G,\widetilde{\chi}]$ must be in $\mathcal{R}(III,G)$.

Now we assume that $G$ is the semidirect product of $H$, and $C$, that $C$ is a minimal nontrivial abelian normal subgroup
of $G$, and that $H$ contains no normal subgroup lying in the centralizer $Z$ of $C$. Then $Z=(Z\cap H)C$, and $Z\cap H$ is a 
normal subgroup of $H$. Hence $Z=C$. If $D$ is a normal subgroup of $G$, and $D$ does not contain $C$, then, by assumption, 
$D\cap C=\{1\}$. The commutator of $C$, and $D$ is contained in $D\cap C$, therefore $D$ is contained in $Z=C$, and $D=\{1\}$.
This means that 
$$[G,\widetilde{\chi}]\cdot([H,1_H]_{G}-\sum_{\mu\in T}[G_\mu',\mu']_{G})$$
 is of {\bf type III}.

\end{proof}

Now we are prepared to prove Theorem \ref{Theorem 1.2}.
\begin{proof}[\large{{\bf Proof of Theorem \ref{Theorem 1.2}}:}]
 Let $\mathcal{F}$ be extendible. Then in the situation of {\bf type I} we have 
 $$Ind_{K}^{B}(\chi)=\sum_{\mu\in\widehat{A}}\mu\chi.$$
 Hence by Equation (\ref{eqn 2.2.3}), and Equation (\ref{eqn 2.2.5})
 $$\mathcal{F}((B,Ind_{K}^{B}(\chi)))=\lambda_{K}^{B}(\mathcal{F})\mathcal{F}((K,\chi))=
 \prod_{\mu\in\widehat{A}}\mathcal{F}((B,\mu\chi)).$$
 Furthermore,
 $$\lambda_{K}^{B}(\mathcal{F})=\mathcal{F}((B,Ind_{K}^{B}(1_K)))=\prod_{\mu\in\widehat{A}}\mathcal{F}((B,\mu))=
 \prod_{\mu\in\widehat{A}}\mathcal{F}((A,\mu)).$$
 This implies Equation (\ref{eqn 5.15}).
 
 In the situation of {\bf type II}, we have 
 $$Ind_{G_1}^{B}(\chi_1\chi)=Ind_{G_2}^{B}(\chi_2\chi).$$
 Hence by Equation (\ref{eqn 2.2.5})
 $$\lambda_{G_1}^{B}(\mathcal{F})\mathcal{F}((G_1,\chi_1\chi))=\lambda_{G_2}^{B}(\mathcal{F})\mathcal{F}((G_2,\chi_2\chi)),$$
 $$\lambda_{G_i}^{B}(\mathcal{F})=\mathcal{F}((B,Ind_{G_i}^{B}(1_{G_i})))=\prod_{\mu\in\widehat{X_i}}\mathcal{F}((X_i,\mu)),\quad
 i=1,2,$$
 which implies Equation (\ref{eqn 5.16}).

In the situation of {\bf type III}, by Lemma \ref{Lemma Koch 2.1} we have
$$Ind_{G'}^{B}(\chi)=\sum_{\mu\in T}Ind_{B_\mu'}^{B}(\mu\chi).$$
Hence 
\begin{equation}\label{eqn 5.18}
 \lambda_{G'}^{B}(\mathcal{F})\mathcal{F}((G',\chi))=
 \prod_{\mu\in T}\lambda_{B_\mu'}^{B}(\mathcal{F})\mathcal{F}((B_\mu',\mu\chi)).
\end{equation}
In particular for $\chi=1_B$
\begin{equation}\label{eqn 5.19}
 \lambda_{G'}^{B}(\mathcal{F})=\prod_{\mu\in T}\lambda_{B_\mu'}^{B}(\mathcal{F})\mathcal{F}((A_\mu',\mu)).
\end{equation}
 Equations (\ref{eqn 5.18}), and (\ref{eqn 5.19}) imply  Equation \ref{eqn 1.3}.

{\bf Conversely,} we assume that $\mathcal{F}$ is a function which satisfies (\ref{eqn 5.15})-(\ref{eqn 1.3}).
We want to show that 
$\mathcal{F}$ is (weakly) extendible .

Now we are able to show that $U\to\lambda_{U}^{G}(\mathcal{F})$ satisfies equation (\ref{eqn 5.20}) for $H=G$. 
Theorem \ref{Lemma 3.2 of Koch} shows that it is sufficient to prove that 
\begin{equation}\label{eqn 5.38}
 \sum_{[U_i,\chi_i]\in\mathcal{M}_G/G}\alpha_{[U_i,\chi_i]}(0) Ind_{U_i}^{G}(\chi_i)\cong 0
\end{equation}
implies
\begin{equation}\label{eqn 5.39}
 \prod_{[U_i,\chi_i]\in\mathcal{M}_G/G}\mathcal{F}((U_i,\chi_i))^{\alpha_{[U_i,\chi_i](0)}}
 \lambda_{U_i}^{G}(\mathcal{F})^{\alpha_{[U_i,\chi_i]}(0)}=1
\end{equation}
if (\ref{eqn 5.38}) is of the form
$$Ind_{B}^{G}(\chi\cdot Inf_{A}^{B}(\sigma'))\cong 0,$$
where $\sigma'$ is of the form {\bf type I, type II}, or {\bf type III}:\\
{\bf type I.}
$$\sigma'=Ind_{e}^{A}(1_e)-\sum_{\mu\in\widehat{A}}\mu.$$
Then 
$$Ind_{B}^{G}(\chi\cdot Inf_{A}^{B}(\sigma'))=Ind_{K}^{G}(\chi)-\sum_{\mu\in\widehat{A}}Ind_{B}^{G}(\mu\chi).$$
Hence we have to show
$$\mathcal{F}((K,\chi))\lambda_{K}^{G}(\mathcal{F})=\prod_{\mu\in\widehat{A}}\mathcal{F}((B,\mu\chi))
\lambda_{B}^{G}(\mathcal{F}).$$
By Lemma \ref{Lemma 3.4 of Koch} this follows from 
\begin{equation}\label{eqn 5.40}
 \mathcal{F}((K,\chi))\lambda_{K}^{G}(\mathcal{F})=\prod_{\mu\in\widehat{A}}\mathcal{F}((B,\mu\chi)).
\end{equation}
By definition of $\lambda_{K}^{B}(\mathcal{F})$,
$$\lambda_{K}^{B}(\mathcal{F})=\lambda_{e}^{A}(\mathcal{F})=\prod_{\mu\in\widehat{A}}\mathcal{F}((A,\mu)).$$
Therefore, Equation (\ref{eqn 5.15}) implies Equation(\ref{eqn 5.40}).

In the case {\bf type II},
$$\sigma=Ind_{H_1}^{A}(\chi_1)-Ind_{H_2}^{A}(\chi_2),$$
and {\bf type III},
$$\sigma=Ind_{H}^{A}(1_H)-\sum_{\mu\in T}Ind_{A_\mu'}^{A}(\mu'),$$
one proceeds in a similar way as in the case {\bf type I}. This concludes the proof of the existence of 
$U\to\lambda_{U}^{H}(\mathcal{F})$ with properties (\ref{eqn 5.20}), (\ref{eqn 5.21}), and by Theorem \ref{Lemma 3.2 of Koch}
the proof of Theorem \ref{Theorem 1.2}.

\end{proof}

\section{{\bf Local root numbers and their properties}}

\subsection{\textbf{Local Fields and their finite extensions}}

Let $F$  be a non-archimedean local field, i.e., a finite extension of the field $\mathbb{Q}_p$ (field of $p$-adic numbers),
where $p$ is a prime.
Let $K/F$ be a finite extension of the field $F$. Let $e_{K/F}$ be the ramification index for the extension $K/F$, and $f_{K/F}$ be 
the residue degree of the extension $K/F$. The extension $K/F$ is called \textbf{unramified} 
if $e_{K/F}=1$; equivalently $f_{K/F}=[K:F]$. The extension $K/F$ is \textbf{totally ramified} if 
$e_{K/F}=[K:F]$; equivalently $f_{K/F}=1$. 
Let
$q_F$ be the cardinality of the residue field $k_F$ of $F$. If $\mathrm{gcd}(p,[K:F])=1$, then the extension 
$K/F$ is called \textbf{tamely ramified}, otherwise \textbf{wildly ramified}. The extension $K/F$ is \textbf{totally tamely ramified}
if it is both totally ramified, and tamely ramified.

For a tower of {\bf local} fields $K/L/F$, we have (cf. \cite{FV}, p. 39, Lemma 2.1) 
\begin{equation}
 e_{K/F}(\nu_K)=e_{K/L}(\nu_K)\cdot e_{L/F}(\nu_L),
\end{equation}
where $\nu_K$ is a valuation on $K$, and $\nu_L$ is the induced 
valuation on $L$, i.e., $\nu_L=\nu_K|_{L}$. For 
the tower of fields $K/L/F$ we simply write $e_{K/F}=e_{K/L}\cdot e_{L/F}$.
Let $O_F$ be the 
ring of integers in the local field $F$, and $P_F=\pi_F O_F$ is the unique prime ideal in $O_F$ 
and $\pi_F$ is a uniformizer, i.e., an element in $P_F$ whose valuation is one, i.e.,
 $\nu_F(\pi_F)=1$.
Let $U_F=O_F-P_F$ be the group of units in $O_F$.
Let $P_{F}^{i}=\{x\in F:\nu_F(x)\geq i\}$, and for $i\geq 0$ define $U_{F}^i=1+P_{F}^{i}$
(with proviso $U_{F}^{0}=U_F=O_{F}^{\times}$).
We also consider that $a(\chi)$ is the conductor of 
 nontrivial character $\chi: F^\times\to \mathbb{C}^\times$, i.e., $a(\chi)$ is the smallest integer $\geq 0$ such 
 that $\chi$ is trivial
 on $U_{F}^{a(\chi)}$. We say $\chi$ is unramified if the conductor of $\chi$ is zero, and otherwise ramified.
Now onwards when $K/F$
is unramified we choose uniformizers $\pi_K=\pi_F$. When $K/F$ is totally ramified (both tame, and wild) we choose
uniformizers $\pi_F=N_{K/F}(\pi_K)$, where $N_{K/F}$ is the norm map from $K^\times$ to $F^\times$.

\begin{dfn}[\textbf{Different, and Discriminant}] 
 Let $K/F$ be a finite separable extension of non-archimedean local field $F$. We define the \textbf{inverse different (or codifferent)}
 $\mathcal{D}_{K/F}^{-1}$ of $K$ over $F$ to be $\pi_{K}^{-d_{K/F}}O_K$, where $d_{K/F}$ is the largest integer (this is the 
 exponent of the different $\mathcal{D}_{K/F}$) such that 
 \begin{center}
  $\mathrm{Tr}_{K/F}(\pi_{K}^{-d_{K/F}}O_K)\subseteq O_F$,
 \end{center}
 where $\rm{Tr}_{K/F}$ is the trace map from $K$ to $F$.
Then the \textbf{different} is defined by:
\begin{center}
 $\mathcal{D}_{K/F}=\pi_{K}^{d_{K/F}}O_K$
\end{center}
and the \textbf{discriminant} $D_{K/F}$ is 
\begin{center}
 $D_{K/F}=N_{K/F}(\pi_{K}^{d_{K/F}})O_F$.
\end{center}
 Thus it is simple to see that $D_{K/F}$ is an \textbf{ideal of $O_F$}.

We know that if $K/F$ is unramified, then $D_{K/F}$ is \textbf{a unit in $O_F$}. If $K/F$ is 
tamely ramified, then 
\begin{equation}\label{eqn 2.2}
 \nu_K(\mathcal{D}_{K/F})=d_{K/F}=e_{K/F} - 1.
\end{equation}
\end{dfn}
(see \cite{JPS}, Chapter III, for details about different, and discriminant of the extension $K/F$.)
We need to mention a very important result of J-P. Serre for our purposes.

\begin{lem}[\cite{JPS}, p. 50, Proposition 7]\label{Lemma 2.21}
Let $K/F$ be a finite separable extension of the field $F$. Let $I_F$ (resp. $I_K$) be a fractional ideal of $F$ (resp. $K$)
relative to $O_F$ (resp. $O_K$). Then the following two properties are equivalent:
\begin{enumerate}
 \item $\mathrm{Tr}_{K/F}(I_K)\subset I_F$.
 \item $I_K\subset I_F.\mathcal{D}_{K/F}^{-1}$.
\end{enumerate}

 \end{lem}

\begin{dfn}[\textbf{Canonical additive character}] 

We define the non trivial additive character of $F$, $\psi_F:F\to\mathbb{C}^\times$ as the composition of the following 
four maps:
\begin{center}
 $F\xrightarrow{\mathrm{Tr}_{F/\mathbb{Q}_p}}\mathbb{Q}_p\xrightarrow{\alpha}\mathbb{Q}_p/\mathbb{Z}_p
 \xrightarrow{\beta}\mathbb{Q}/\mathbb{Z}\xrightarrow{\gamma}\mathbb{C}^\times$,
\end{center}
where
\begin{enumerate}
 \item $\mathrm{Tr}_{F/\mathbb{Q}_p}$ is the trace from $F$ to $\mathbb{Q}_p$,
 \item $\alpha$ is the canonical surjection map,
 \item $\beta$ is the canonical injection which maps $\mathbb{Q}_p/\mathbb{Z}_p$ onto the $p$-component of the 
 divisible group $\mathbb{Q}/\mathbb{Z}$, and 
 \item $\gamma$ is the exponential map $x\mapsto e^{2\pi i x}$, where $i=\sqrt{-1}$.
\end{enumerate}
For every $x\in\mathbb{Q}_p$, there is a rational $r$, uniquely determined modulo $1$, such that $x-r\in\mathbb{Z}_p$.
Then $\psi_{\bbQ_p}(x)=\psi_{\bbQ_p}(r)=e^{2\pi i r}$.
The nontrivial additive character  $\psi_F=\psi_{\bbQ_p}\circ \rm{Tr}_{F/\bbQ_p}$ of $F$ 
is called the \textbf{canonical additive character} (cf. \cite{JT1}, p. 92).
\end{dfn}
The {\bf conductor} of any nontrivial additive character $\psi$ of the field $F$ is an integer $n(\psi)$ if $\psi$ is trivial
on $P_{F}^{-n(\psi)}$, but nontrivial on $P_{F}^{-n(\psi)-1}$. So, from Lemma \ref{Lemma 2.21} we can 
observe that 
\begin{center}
 $n(\psi_F)=n(\psi_{\bbQ_p}\circ\rm{Tr}_{F/\bbQ_p})= \nu_F(\mathcal{D}_{F/\bbQ_p})$,
\end{center}
because $d_{\bbQ_p/\bbQ_p}=0$, and hence $n(\psi_{\bbQ_p})=0$.

\subsection{\textbf{Local Constants (or local root numbers/local epsilon factors)}}

Let $F$ be a non-archimedean local field, and $\chi$ be a character of $F^\times$.
 The $L(\chi)$-functions are defined as follows:
\begin{align*}
 L(\chi)
 &=\begin{cases}
            (1-\chi(\pi_F))^{-1} & \text{if $\chi$ is unramified},\\
            1 & \text{if $\chi$ is ramified}.
           \end{cases}
\end{align*}

 We denote by $dx$ a Haar measure on $F$, by $d^\times x$ a Haar measure on $F^\times$, and the relation between these two 
 Haar measure is: 
 \begin{center}
  $d^\times x=\frac{dx}{|x|}$,
 \end{center}
for arbitrary Haar measure $dx$ on $F$. For a given additive character $\psi$ of $F$, and Haar measure $dx$ on $F$, we have a 
\textbf{Fourier transform} as:
\begin{equation}
 \hat{f}(y)=\int f(x)\psi(xy)dx.
\end{equation}
where $f\in L^{1}(F^{+})$ (that is, $|f|$ is integrable), and the Haar measure is normalized 
such that $\hat{\hat{f}}(y)=f(-y)$, i.e., $dx$ is self-dual with respect to $\psi$.
By Tate (cf. \cite{JT2}, p. 13), for any character $\chi$ of $F^\times$, there exists 
a number $W(\chi,\psi,dx)\in\mathbb{C}^\times$ such that it satisfies the following local functional equation:
\begin{equation}\label{eqn 4.4}
 \frac{\int\hat{f}(x)w_1\chi^{-1}(x)d^\times x}{L(w_1\chi^{-1})}=W(\chi,\psi,dx)
 \frac{\int f(x)\chi(x)d^\times x}{L(\chi)}.
\end{equation}
for any such function $f$ for which the both sides make sense. Here $w_s(x)=|x|_F^s=q_F^{-s\nu_F(x)}$ is
unramified character of $F^\times$. The number $W(\chi,\psi, dx)$ is called the
{\bf local epsilon factor or local constant} of $\chi$.


For a nontrivial multiplicative character $\chi$ of $F^\times$, and nontrivial additive character $\psi$ of $F$,
we have (cf. \cite{RL}, p. 4)
\begin{equation}
 W(\chi,\psi,c)=\chi(c)\frac{\int_{U_F}\chi^{-1}(x)\psi(x/c) dx}{|\int_{U_F}\chi^{-1}(x)\psi(x/c) dx|}\label{label1}
\end{equation}
where the Haar measure $dx$ is normalized such that the measure of $O_F$ is $1$, and 
where $c\in F^\times$ with valuation $n(\psi)+a(\chi)$.
The \textbf{modified} formula of local constant (cf. \cite{JT1}, p. 94) is:
\begin{equation}\label{eqn 4.6}
 W(\chi,\psi,c)=\chi(c)q^{-a(\chi)/2}\sum_{x\in\frac{U_F}{U_{F}^{a(\chi)}}}\chi^{-1}(x)\psi(x/c).
\end{equation}
where $c=\pi_{F}^{a(\chi)+n(\psi)}$. Now if $u\in U_F$ is unit, and replace $c=cu$, then we have 
\begin{equation}
 W(\chi,\psi,cu)=\chi(c)q^{-\frac{a(\chi)}{2}}\sum_{x\in\frac{U_F}{U_{F}^{a(\chi)}}}\chi^{-1}(x/u)\psi(x/cu)=W(\chi,\psi,c).
\end{equation}
Therefore, $W(\chi,\psi,c)$ \textbf{depends} only on the exponent $\nu_{F}(c)=a(\chi)+n(\psi)$. Therefore, we can 
simply write $W(\chi,\psi, c)=W(\chi,\psi)$, because $c$ is determined by 
$\nu_F(c)=a(\chi)+n(\psi)$ up to a unit $u$ which has \textbf{no influence on} $W(\chi,\psi,c)$.
If $\chi$ is unramified, i.e., $a(\chi)=0$, therefore $\nu_F(c)=n(\psi)$. Then from the formula of $W(\chi,\psi,c)$, we can write
\begin{equation}\label{eqn 2.3.5}
 W(\chi,\psi,c)=\chi(c),
\end{equation}
and therefore $W(1,\psi,c)=1$ if $\chi=1$ is the trivial character.

\begin{thm}[{\bf Lamprecht-Tate formula, Proposition 1 of \cite{JT1}}]\label{Theorem 6.1.1}
\footnote{ Before Tate, in 1953 in \cite{EL} Erich Lamprecht  gave 
a formula for local epsilon factors of linear characters. Then Tate generalizes the formula for epsilon factors.}
Let $F$ be a non-Archimedean local field.
Let $\chi$ be a character of $F^\times$ of conductor $a(\chi)$, and let $m$ be a natural number 
such that $2m\le a(\chi)$. Let $\psi$ be a nontrivial additive character of $F$. 
Then there exists $c\in F^\times$, $\nu_F(c)=a(\chi)+n(\psi)$ such that 
\begin{equation}\label{eqn 5.4.5}
 \chi(1+y)=\psi(c^{-1}y)\qquad\text{for all $y\in P_{F}^{a(\chi)-m}$},
\end{equation}
and for such a $c$ we have:
\begin{equation}\label{eqn 6.0.9}
 W(\chi,\psi)=\chi(c)\cdot q_{F}^{-\frac{(a(\chi)-2m)}{2}}
 \sum_{x\in U_F^m/U_F^{a(\chi)-m}}\chi^{-1}(x)\psi(c^{-1}x).
\end{equation}
{\bf Remark:} Note that the assumption (\ref{eqn 5.4.5}) is obviously fulfilled for $m=0$ because then both sides are $=1$,
and the resulting formula for $m=0$ is the Tate formula (\ref{eqn 4.6}).
\end{thm}
For the proof we refer \cite{SABLT}.

\subsection{\textbf{Some properties of $W(\chi,\psi)$}}

 \begin{enumerate}
 
  \item Let $b\in F^\times$ be the uniquely determined element such that $\psi'=b\psi$. Then 
  \begin{equation}
   W(\chi,\psi',c')=\chi(b)\cdot W(\chi,\psi,c).
  \end{equation}
\begin{proof}
 Here $\psi'(x)=(b\psi)(x):=\psi(bx)$ for all $x\in F$. It is an additive character of $F$. The existence, and uniqueness of $b$
 is clear. From the definition of conductor of an 
 additive character we have 
 \begin{center}
  $n(\psi')=n(b\psi)=n(\psi)+\nu_F(b)$.
 \end{center}
Here $c'\in F^\times$ is of valuation 
$$\nu_F(c')=a(\chi)+n(\psi')=a(\chi)+\nu_F(b)+n(\psi)=\nu_F(b)+\nu_F(c)=\nu_F(bc).$$
Therefore, $c'=bcu$ where $u\in U_F$ is some unit. Now 
\begin{align*}
 W(\chi,\psi',c')
 &=W(\chi,b\psi,bcu)\\
 &=W(\chi,b\psi,bc)\\
 &=\chi(bc)q_{F}^{-\frac{a(\chi)}{2}}\sum_{x\in\frac{U_F}{U_{F}^{a(\chi)}}}\chi^{-1}(x)((bc)^{-1}(b\psi))(x)\\
 &=\chi(b)\cdot\chi(c)q_{F}^{-\frac{a(\chi)}{2}}\sum_{x\in\frac{U_F}{U_{F}^{a(\chi)}}}\chi^{-1}(x)\psi(xc^{-1})\\
 &=\chi(b)\cdot W(\chi,\psi,c).
\end{align*}

\end{proof}

\item Let $F/\bbQ_p$ be a local field inside $\overline{\bbQ_p}$.
Let $\chi$, and $\psi$ be a character of $F^\times$, and $F^{+}$ respectively, and $c\in F^\times$ with 
valuation $\nu_F(c)=a(\chi)+n(\psi)$. If $\sigma\in\mathrm{Gal}(\overline{\mathbb{Q}_p}/\mathbb{Q}_p)$ is an automorphism, then:
\begin{center}
 $W_{F}(\chi,\psi,c)=W_{\sigma^{-1}(F)}(\chi^{\sigma},\psi^{\sigma},\sigma^{-1}(c))$,
\end{center}
where $\chi^{\sigma}(y):=\chi(\sigma(y))$, $\psi^{\sigma}(y):=\psi(\sigma(y))$, for all $y\in\sigma^{-1}(F)$.

\begin{proof}
Let $L:=\sigma^{-1}(F)$. Because $\sigma$ is an automorphism of $\overline{\bbQ_p}$, then we have $O_F/P_F\cong O_L/P_L$, hence 
$q_F=q_L$. We also can see that $a(\chi^{\sigma})=a(\chi)$, and $n(\psi^{\sigma})=n(\psi)$. Then from the formula of local constant
we have 
\begin{align*}
 W_{\sigma^{-1}(F)}(\chi^{\sigma},\psi^{\sigma},\sigma^{-1}(c))
 &=W_{L}(\chi^{\sigma},\psi^{\sigma},\sigma^{-1}(c))\\
 &=\chi^{\sigma}(\sigma^{-1}(c))q_{L}^{-\frac{a(\chi^{\sigma})}{2}}\sum_{y\in\frac{U_{L}}{U_{L}^{a(\chi^{\sigma})}}}
 (\chi^{\sigma})^{-1}(y)\cdot((\sigma^{-1}(c))^{-1}\psi^{\sigma}(y)\\
 &=\chi(c)q_{F}^{-\frac{a(\chi)}{2}}\sum_{x\in \frac{U_F}{U_{F}^{a(\chi)}}}\chi^{-1}(x)\psi(\frac{x}{c})\\
 &=W_{F}(\chi,\psi,c).
\end{align*}
Here we put 
$y=\sigma^{-1}(x)$, and use $(\sigma^{-1}(c))^{-1}\psi^{\sigma}=(c^{-1}\psi)^{\sigma}$.

\end{proof}
\begin{rem}
 We can simply write as before $W_{F}(\chi,\psi)=W_{\sigma^{-1}(F)}(\chi^{\sigma},\psi^{\sigma})$.
 Tate in his paper \cite{JT1} on local constants defines the local root number as:
\begin{center}
 $W_{F}(\chi):=W_F(\chi,\psi_F)=W_F(\chi,\psi_F,d)$,
\end{center}
where $\psi_F$ is the canonical character of $F^\times$, and $d\in F^\times$ with
$\nu_F(d)=a(\chi)+n(\psi_F)$.
Therefore, after fixing canonical additive character $\psi=\psi_F$, we can rewrite 
\begin{center}
 $W_F(\chi)=\chi(d(\psi_F))$, if $\chi$ is unramified,\\
 $W_F(\chi)=W_{\sigma^{-1}(F)}(\chi^{\sigma})$.
\end{center}
The last equality follows because the canonical character $\psi_{\sigma^{-1}(F)}$ is related to the canonical character 
$\psi_F$ as: $\psi_{\sigma^{-1}(F)}=\psi_{F}^{\sigma}$.\\

So we see that 
\begin{equation*}
 (F,\chi)\to W_F(\chi)\in\mathbb{C}^\times
\end{equation*}
is a function with properties (\ref{eqn 2.2.1}), (\ref{eqn 2.2.2}) of extendible  functions.
\end{rem}

\item If $\chi\in\widehat{F^\times}$, and $\psi\in\widehat{F}$, then 
\begin{equation*}
 W(\chi,\psi)\cdot W(\chi^{-1},\psi)=\chi(-1).
\end{equation*}
Furthermore if the character $\chi:F^\times\to\mathbb{C}^{\times}$ is unitary (in particular, if $\chi$ is of finite order), then 
 \begin{center}
  $|W(\chi,\psi)|^{2}=1$.
  \end{center}

\begin{proof}
We prove this properties by using equation (\ref{eqn 4.6}). 
We know that the additive characters are always unitary, hence
\begin{center}
 $\psi(-x)=\psi(x)^{-1}=\overline{\psi}(x)$.
\end{center}
On the other hand we write $\psi(-x)=((-1)\psi)(x)$, where $-1\in F^\times$. Therefore, $\overline{\psi}=(-1)\psi$.
We also have $a(\chi)=a(\chi^{-1})$. Therefore, by using equation (\ref{eqn 4.6}) we have 
\begin{align*}
 W(\chi,\psi)\cdot W(\chi^{-1},\psi)
 &=\chi(-1)\cdot q_{F}^{-a(\chi)}\sum_{x,y\in \frac{U_F}{U_{F}^{a(\chi)}}}\chi^{-1}(x)\chi(y)\psi(\frac{x-y}{c})\\
 &=\chi(-1)\cdot q_{F}^{-a(\chi)}\sum_{x,y\in \frac{U_F}{U_{F}^{a(\chi)}}}\chi^{-1}(x)\psi(\frac{xy-y}{c}),
 \quad\text{replacing $x$ by $xy$}\\
 &=\chi(-1)\cdot q_{F}^{-a(\chi)}\sum_{x\in \frac{U_F}{U_{F}^{a(\chi)}}}\chi^{-1}(x)\varphi(x),
\end{align*}
where
\begin{equation}
\varphi(x)=\sum_{y\in\frac{U_F}{U_{F}^{a(\chi)}}}\psi(y\frac{x-1}{c}).
\end{equation}
Because $\frac{U_F}{U_{F}^{a(\chi)}}=(\frac{O_F}{P_{F}^{a(\chi)}})^\times=
\frac{O_F}{P_{F}^{a(\chi)}}\setminus\frac{P_F}{P_{F}^{a(\chi)}}$, therefore $\varphi(x)$
can be written as the difference 
\begin{align*}
\varphi(x)
&=\sum_{y\in\frac{U_F}{U_{F}^{a(\chi)}}}\psi(y\frac{x-1}{c})\\
&=\sum_{y\in\frac{O_F}{P_{F}^{a(\chi)}}}\psi(y\frac{x-1}{c})-
\sum_{y\in\frac{P_F}{P_{F}^{a(\chi)}}}\psi(y\frac{x-1}{c})\\
&=\sum_{y\in\frac{O_F}{P_{F}^{a(\chi)}}}\psi(y\frac{x-1}{c})-
\sum_{y\in\frac{O_F}{P_{F}^{a(\chi)-1}}}\psi(y\frac{(x-1)\pi_F}{c})\\
&=A-B,
\end{align*}
where $A=\sum_{y\in\frac{O_F}{P_{F}^{a(\chi)}}}\psi(y\frac{x-1}{c})$
and $B=\sum_{y\in\frac{O_F}{P_{F}^{a(\chi)-1}}}\psi(y\frac{(x-1)\pi_F}{c})$. It is simple to see that (cf. \cite{M}, p. 28, Lemma 2.1)
\begin{align*}
 \sum_{y\in\frac{O_F}{P_{F}^{a(\chi)}}}\psi(y\alpha)=\begin{cases}
                                                      q_{F}^{a(\chi)} & \text{when $\alpha\in P_{F}^{-n(\psi)}$}\\
                                                      0 & \text{otherwise}
                                                     \end{cases}
\end{align*}
Therefore, $A=q_{F}^{a(\chi)}$ when $x\in U_{F}^{a(\chi)}$, and $A=0$ otherwise. Similarly
$B=q_{F}^{a(\chi)-1}$ when $x\in U_{F}^{a(\chi)-1}$
, and $B=0$ otherwise. Therefore, we have
 \begin{align*}
  W(\chi,\psi)\cdot W(\chi^{-1},\psi)
  &=\chi(-1)\cdot q_{F}^{-a(\chi)}\cdot\{q_{F}^{a(\chi)}-q_{F}^{a(\chi)-1}\sum_{x\in\frac{U_{F}^{a(\chi)-1}}{U_{F}^{a(\chi)}}}\chi^{-1}(x)\}\\
  &=\chi(-1)-\chi(-1)\cdot q_{F}^{-1}\sum_{x\in\frac{U_{F}^{a(\chi)-1}}{U_{F}^{a(\chi)}}}\chi^{-1}(x).
 \end{align*}
Because the conductor of $\chi$ is $a(\chi)$, then it can be proved that
$\sum_{x\in\frac{U_{F}^{a(\chi)-1}}{U_{F}^{a(\chi)}}}\chi^{-1}(x)=0$.
Thus we obtain
\begin{equation}\label{eqn 2.3.9}
 W(\chi,\psi)\cdot W(\chi^{-1},\psi)=\chi(-1).
\end{equation}
\vspace{.3cm}

The right side of equation (\ref{eqn 2.3.9}) is a sign, hence we may rewrite (\ref{eqn 2.3.9}) as 
\begin{center}
 $W(\chi,\psi)\cdot\chi(-1)W(\chi^{-1},\psi)=1$.
\end{center}
However, we also know from our earlier property that 
$$\chi(-1)W(\chi^{-1},\psi)=W(\chi^{-1},(-1)\psi)=W(\chi^{-1},\overline{\psi}).$$
So the identity (\ref{eqn 2.3.9}) rewrites as 
\begin{center}
 $W(\chi,\psi)\cdot W(\chi^{-1},\overline{\psi})=1$.
\end{center}
Now we assume that $\chi$ is unitary, hence 
\begin{center}
 $W(\chi^{-1},\overline{\psi})=W(\overline{\chi},\overline{\psi})=\overline{W(\chi,\psi)}$
\end{center}
where the last equality is obvious. Now we see that for unitary $\chi$ the identity (\ref{eqn 2.3.9}) rewrites as 
\begin{center}
 $|W(\chi,\psi)|^{2}=1$.
\end{center}
\end{proof}

\begin{rem}
From the functional equation (\ref{eqn 4.4}), we can directly see the first part of the above property of local constant. 
Denote 
\begin{equation}\label{eqn 66}
 \zeta(f,\chi)=\int f(x)\chi(x)d^\times x.
\end{equation}
Now replacing  $f$ by $\hat{\hat{f}}$ in equation (\ref{eqn 66}), and we get
\begin{equation}\label{eqn 77}
 \zeta(\hat{\hat{f}},\chi)=\int \hat{\hat{f}}(x)\chi(x)d^\times x=\chi(-1)\cdot\zeta(f,\chi),
\end{equation}
because $dx$ is self-dual with respect to $\psi$, hence $\hat{\hat{f}}(x)=f(-x)$ for all $x\in F^{+}$.

Again the functional equation (\ref{eqn 4.4}) can be written as follows:
\begin{equation}\label{eqn 88}
 \zeta(\hat{f},w_1\chi^{-1})=W(\chi,\psi,dx)\cdot\frac{L(w_1\chi^{-1})}{L(\chi)}\cdot\zeta(f,\chi).
\end{equation}

 Now we replace $f$ by $\hat{f}$, and $\chi$ by $w_1\chi^{-1}$ in equation (\ref{eqn 88}), and we obtain
 \begin{equation}\label{eqn 99}
  \zeta(\hat{\hat{f}},\chi)=W(w_1\chi^{-1},\psi,dx)\cdot\frac{L(\chi)}{L(w_1\chi^{-1})}\cdot\zeta(\hat{f},w_1\chi^{-1}).
 \end{equation}
Then by using equations (\ref{eqn 77}), (\ref{eqn 88}), from the above equation (\ref{eqn 99}) we obtain:
\begin{equation}\label{eqn 100}
 W(\chi,\psi,dx)\cdot W(w_1\chi^{-1},\psi,dx)=\chi(-1).
\end{equation}
Moreover, the convention $W(\chi,\psi)$ is actually as follows (cf. \cite{JT2}, p. 17, equation (3.6.4)):
$$W(\chi w_{s-\frac{1}{2}},\psi)=W(\chi w_s,\psi,dx).$$
By using this relation from equation (\ref{eqn 100}) we can conclude
$$W(\chi,\psi)\cdot W(\chi^{-1},\psi)=\chi(-1).$$

\end{rem}

\item \textbf{Twisting formula of abelian local constants:}\label{Twisting formula for characters}

\begin{enumerate}
 \item If $\chi_1$, and $\chi_2$ are two unramified characters of $F^\times$, and $\psi$ is
 a nontrivial additive character of $F$, then from equation (\ref{eqn 2.3.5}) we have
 \begin{equation}
  W(\chi_1\chi_2,\psi)=W(\chi_1,\psi)W(\chi_2,\psi).
 \end{equation}
 \item  Let $\chi_1$ be ramified, and $\chi_2$ unramified then (cf. \cite{JT2}, (3.2.6.3))
\begin{equation}
 W(\chi_1\chi_2,\psi)=\chi_2(\pi_F)^{a(\chi_1)+n(\psi)}\cdot W(\chi_1,\psi).
\end{equation}
\begin{proof}
  By the given condition
$a(\chi_1)>a(\chi_2)=0$. Therefore, $a(\chi_1\chi_2)=a(\chi_1)$. Then we have
\begin{align*}
 W(\chi_1\chi_2,\psi)
 &=\chi_1\chi_2(c)q_{F}^{-a(\chi_1)/2}\sum_{x\in\frac{U_F}{U_{F}^{a(\chi)}}}(\chi_1\chi_2)^{-1}(x)\psi(x/c)\\
 &=\chi_1(c)\chi_2(c)q_{F}^{-a(\chi_1)/2}\sum_{x\in\frac{U_F}{U_{F}^{a(\chi)}}}\chi_{1}^{-1}(x)\chi_{2}^{-1}(x)\psi(x/c)\\
 &=\chi_2(c)\chi_1(c)q_{F}^{-a(\chi_1)/2}\sum_{x\in\frac{U_F}{U_{F}^{a(\chi)}}}\chi_{1}^{-1}(x)\psi(x/c),
 \quad\text{because $\chi_2$ unramified}\\
 &=\chi_2(c)W(\chi_1,\psi)\\
 &=\chi_2(\pi_F)^{a(\chi_1)+n(\psi)}\cdot W(\chi_1,\psi).
\end{align*}
\end{proof}
\item 
 We also have a twisting formula of local constants by Deligne (cf. \cite{D1}, Lemma 4.16)
 under some special condition, and which is as follows (for proof, cf. Corollary 3.2 (2) of \cite{SABLT}):\\
Let $\alpha$, and $\beta$ be two multiplicative characters of a local field $F$ such that $a(\alpha)\geq 2\cdot a(\beta)$.
Let $\psi$ be an additive character of $F$.
Let $y_{\alpha,\psi}$ be an element of $F^\times$ such that 
$$\alpha(1+x)=\psi(y_{\alpha,\psi}x)$$
for all $x\in F$ with valuation $\nu_F(x)\geq\frac{a(\alpha)}{2}$ (if $a(\alpha)=0$, $y_{\alpha,\psi}=\pi_{F}^{-n(\psi)}$). Then 
\begin{equation}\label{eqn 2.3.17}
 W(\alpha\beta,\psi)=\beta^{-1}(y_{\alpha,\psi})\cdot W(\alpha,\psi).
\end{equation}

\item {\bf General twisting formula for characters:} 
In the following theorem, one can see a generalized twisting formula of local constants using local Jacobi sums.

\begin{thm}[Theorem 3.5 on p. 592 of \cite{SABTF}]\label{Theorem 5.1} 
Let $F$ be a non-Archimedean local field with $q$ as the cardinality of the residue field of $F$.
Let $\psi$ be a nontrivial additive character of $F$.
Let $\chi_1$, and $\chi_2$ be two ramified 
characters of $F^\times$ with conductors $n$, and $m$ respectively. Let $r$ be the conductor of character $\chi_1\chi_2$. Then
\begin{align}\label{eqn 4.8}
 W(\chi_1\chi_2,\psi)
 =\begin{cases}
   \frac{q^{\frac{n}{2}}W(\chi_1,\psi)W(\chi_2,\psi)}{J_1(\chi_1,\chi_2,n)} & \text{when $n=m=r$},\\ 
    \frac{q^{\frac{r}{2}}\chi_1\chi_2(\pi_{F}^{r-n})W(\chi_1,\psi)W(\chi_2,\psi)}{J_1(\chi_1,\chi_2,n)} & \text{when $n=m>r$},\\
     \frac{q^{n-\frac{m}{2}}W(\chi_1,\psi)W(\chi_2,\psi)}{J_1(\chi_1,\chi_2,n)} & \text{when $n=r>m$},\\
   \end{cases}
\end{align}
Here the local Jacobi sum is:
\begin{equation}
  J_t(\chi_1,\chi_2,n)=\sum_{\substack{x\in \frac{U_F}{U_{F}^{n}}\\t-x\in U_F}}\chi_{1}^{-1}(x)\chi_{2}^{-1}(t-x).
 \end{equation}
\end{thm}

\end{enumerate}

\end{enumerate}

\subsection{\textbf{Connection of different conventions for local constants} }

Mainly there are two conventions for local constants. They are due to Langlands (\cite{RL}), and Deligne (\cite{D1}). 
Recently Bushnell, and Henniart (\cite{BH}) also give a convention of local constants. In this subsection we shall show the connection
among all three conventions for local constants\footnote{The convention $W(\chi,\psi)$ is actually due to Langlands \cite{RL},
and it is:
\begin{center}
 $\epsilon_{L}(\chi,\frac{1}{2},\psi)=W(\chi,\psi).$
\end{center}
See equation (3.6.4) on p. 17 of \cite{JT2} for $V=\chi$.} 
We denote $\epsilon_{BH}$ as local constant of Bushnell-Henniart (introduced in
Bushnell-Henniart, \cite{BH}, Chapter 6).

On page 142 of \cite{BH}, the authors define a rational function 
$\epsilon_{BH}(\chi,\psi,s)\in\mathbb{C}(q_{F}^{-s})$. From Theorem 23.5 on p. 144 of \cite{BH} for ramified character 
$\chi\in\widehat{F^\times}$, and conductor\footnote{The definition of level of an additive character $\psi\in\widehat{F}$ 
in \cite{BH} on p. 11 is the negative sign with our conductor $n(\psi)$, i.e., level of $\psi=-n(\psi)$.} $n(\psi)=-1$ we have 
\begin{equation}\label{eqn 2.3.12}
 \epsilon_{BH}(\chi,s,\psi)=
 q_{F}^{n(\frac{1}{2}-s)}\sum_{x\in\frac{U_F}{U_{F}^{n+1}}}\chi(\alpha x)^{-1}\psi(\alpha x)/q_{F}^{\frac{n+1}{2}},
\end{equation}
where  $n=a(\chi)-1$, and $\alpha\in F^\times$ with $\nu_{F}(\alpha)=-n$.

Also from the Proposition 23.5 of \cite{BH} on p. 143 for unramified character $\chi\in\widehat{F^\times}$, and $n(\psi)=-1$ we have 
\begin{equation}\label{eqn 2.3.13}
\epsilon_{BH}(\chi,s,\psi)=q_{F}^{s-\frac{1}{2}}\chi(\pi_F)^{-1}. 
\end{equation}

\begin{enumerate}
 \item \textbf{Connection between $\epsilon_{BH}$, and $W(\chi,\psi)$.}
 \begin{center}
 $W(\chi,\psi)=\epsilon_{BH}(\chi,\frac{1}{2},\psi)$.
\end{center}
 \begin{proof}
 From \cite{BH}, p. 143, Lemma 1 we see:
 \begin{center}
  $\epsilon_{BH}(\chi,\frac{1}{2},b\psi)=\chi(b)\epsilon_{BH}(\chi,\frac{1}{2},\psi)$
 \end{center}
for any $b\in F^\times$. However, we have seen already that $W(\chi,b\psi)=\chi(b)W(\chi,\psi)$ has the same transformation rule. If we fix
one nontrivial $\psi$ then all other nontrivial $\psi'$ are uniquely given as $\psi'=b\psi$ for some $b\in F^\times$. Because of the
parallel transformation rules it is now enough to verify our assertion for a single $\psi$. Now we take $\psi\in\widehat{F^{+}}$
with $n(\psi)=-1$, hence $\nu_F(c)=a(\chi)-1$. Then we obtain
\begin{equation*}
 W(\chi,\psi)=W(\chi,\psi,c)=\chi(c)q_{F}^{-\frac{a(\chi)}{2}}\sum_{x\in\frac{U_F}{U_{F}^{a(\chi)}}}\chi^{-1}(x)\psi(c^{-1}x).
\end{equation*}
We compare this to the equation (\ref{eqn 2.3.12}). There the notation is $n=a(\chi)-1$, and the assumption is $n\geq 0$. This means
we have $\nu_F(c)=n$, hence we may take $\alpha=c^{-1}$, and then comparing our formula with equation (\ref{eqn 2.3.12}), we see that
\begin{center}
 $W(\chi,\psi)=\epsilon_{BH}(\chi,\frac{1}{2},\psi)$
\end{center}
in the case when $n(\psi)=-1$.\\
We are still left to prove our assertion if $\chi$ is unramified, i.e., $a(\chi)=0$. Again we can reduce to the case where 
$n(\psi)=-1$. Then our assertion follows from equation \ref{eqn 2.3.13}.

\end{proof}
\begin{rem}
From Corollary $23.4.2$ of \cite{BH}, on p. 142, for $s\in\mathbb{C}$, we 
 can write 
 \begin{align}
  \epsilon_{BH}(\chi,s,\psi)=q_{F}^{(\frac{1}{2}-s)n(\chi,\psi)}\cdot\epsilon_{BH}(\chi,\frac{1}{2},\psi),
 \end{align}
 for some $n(\chi,\psi)\in\mathbb{Z}$. In fact here $n(\chi,\psi)=a(\chi)+n(\psi)$.
From above connection, we only see $W(\chi,\psi)=\epsilon_{BH}(\chi,\frac{1}{2},\psi)$. Thus for arbitrary $s\in\mathbb{C}$, we 
obtain
\begin{equation}\label{eqn 2.3.10}
\epsilon_{BH}(\chi,s,\psi)=q_{F}^{(\frac{1}{2}-s)(a(\chi)+n(\psi))}\cdot W(\chi,\psi).
\end{equation}
This equation (\ref{eqn 2.3.10}) is very important for us. We shall use this to connect with Deligne's convention.

In \cite{JT2} there is defined a number $\epsilon_{D}(\chi,\psi,dx)$ depending on $\chi$, $\psi$, and a Haar measure
$dx$ on $F$. This notion is due to Deligne \cite{D1}. We write 
 $\epsilon_{D}$ for Deligne's convention in order to distinguish it from the $\epsilon_{BH}(\chi,\frac{1}{2},\psi)$ introduced
in Bushnell-Henniart \cite{BH}.

In the next Lemma we give the connection between Bushnell-Henniart, and Deligne conventions for local constants.
\end{rem}

\item {\bf The connection between $\epsilon_D$, and $\epsilon_{BH}$}:
\begin{lem}
We have the relation
 \begin{center}
 $\epsilon_{BH}(\chi,s,\psi)=\epsilon_{D}(\chi\cdot\omega_{s},\psi,dx_{\psi})$,
\end{center}
where $\omega_{s}(x)=|x|_{F}^{s}=q^{-s\nu_F(x)}$ is unramified character of $F^\times$ corresponding to complex number $s$, and where 
$dx_{\psi}$ is the self-dual Haar measure corresponding to the additive character $\psi$.
\end{lem}

\begin{proof}

From equation equation (3.6.4) of \cite{JT2}, we know that
\begin{equation}\label{eqn 2.3.11}
 \epsilon_{L}(\chi,s,\psi):=\epsilon_{L}(\chi\omega_{s-\frac{1}{2}},\psi)=\epsilon_{D}(\chi\omega_{s},\psi,dx_{\psi}).
\end{equation}
 We prove this connection by using the relations (\ref{eqn 2.3.10}), and (\ref{eqn 2.3.11}). From equation (\ref{eqn 2.3.11}) we can write 
 our $W(\chi,\psi)=\epsilon_{D}(\chi\omega_{\frac{1}{2}},\psi,dx_{\psi})$. Therefore, when $s=\frac{1}{2}$, we have the 
 relation:
 \begin{equation}
  \epsilon_{BH}(\chi,\frac{1}{2},\psi)=\epsilon_{D}(\chi\omega_{\frac{1}{2}},\psi,dx_{\psi}),
 \end{equation}
because $W(\chi,\psi)=\epsilon_{BH}(\chi,\frac{1}{2},\psi)$.

We know that $\omega_{s}(x)=q_{F}^{-s\nu_F(x)}$ is an unramified character of $F^\times$.
So when $\chi$ is also unramified, we can write 
\begin{equation}
 W(\chi\omega_{s-\frac{1}{2}},\psi)=
 \omega_{s-\frac{1}{2}}(c)\cdot\chi(c)=q_{F}^{(\frac{1}{2}-s)n(\psi)}\epsilon_{BH}(\chi,\frac{1}{2},\psi)
 =\epsilon_{BH}(\chi,s,\psi).
\end{equation}
And when $\chi$ is ramified character, i.e., conductor $a(\chi)>0$,
from Tate's theorem for unramified twist (see property (4.20)) , we can write 
\begin{align*}
 W(\chi\omega_{s-\frac{1}{2}},\psi)
 &=\omega_{s-\frac{1}{2}}(\pi_{F}^{a(\chi)+n(\psi)})\cdot W(\chi,\psi)\\
 &=q_{F}^{-(s-\frac{1}{2})(a(\chi)+n(\psi))}\cdot W(\chi,\psi)\\
 &=q_{F}^{(\frac{1}{2}-s)(a(\chi)+n(\psi))}\cdot \epsilon_{BH}(\chi,\frac{1}{2},\psi)\\
 &=\epsilon_{BH}(\chi,s,\psi).
\end{align*}
Furthermore from equation (\ref{eqn 2.3.11}), we have 
\begin{equation}
 W(\chi\omega_{s-\frac{1}{2}},\psi)=\epsilon_{D}(\chi\omega_{s},\psi,dx_{\psi}).
\end{equation}
Therefore, finally we can write 
\begin{equation}
 \epsilon_{BH}(\chi,s,\psi)=\epsilon_{D}(\chi\omega_{s},\psi,dx_{\psi}).
\end{equation}
\end{proof}

 \begin{cor}\label{Corollary 4.1}
  For our $W$ we have :
  \begin{center}
   $W(\chi,\psi)=\epsilon_{BH}(\chi,\frac{1}{2},\psi)=\epsilon_{D}(\chi\omega_{\frac{1}{2}},\psi,dx_{\psi})$\\
   $W(\chi\omega_{s-\frac{1}{2}},\psi)=\epsilon_{BH}(\chi,s,\psi)=\epsilon_{D}(\chi\omega_{s},\psi,dx_{\psi})$.
  \end{center}
\end{cor}
\begin{proof}
 From the equations (3.6.1), and (3.6.4) of \cite{JT2} for $\chi$, and above two connections the assertions follow.
\end{proof}
\end{enumerate}

\subsection{{\bf Local constants for virtual representations}}\label{Subsection 2.5}

\begin{enumerate}
 \item To extend the concept of local constant, we need to go from one-dimensional to other virtual representations $\rho$ of 
 the Weil groups $W_F$ of non-archimedean local field $F$.
 According to Tate \cite{JT1}, the root 
number $W(\chi):=W(\chi,\psi_F)$ extensions to $W(\rho)$, where $\psi_F$ is the canonical additive character of $F$.
More generally, $W(\chi,\psi)$ extends to $W(\rho,\psi)$, and if
$E/F$ is a finite separable extension then one has to take 
$\psi_{E}=\psi_{F}\circ \mathrm{Tr}_{E/F}$ for the extension field $E$. 

According to Bushnell-Henniart \cite{BH}, Theorem on p. 189, the functions
$\epsilon_{BH}(\chi,s,\psi)$ extend to $\epsilon_{BH}(\rho,s,\psi_E)$, where $\psi_E=\psi\circ\rm{Tr}_{E/F}$
\footnote{ Note that they fix a base field $F$, and a nontrivial $\psi=\psi_F$
(which not to be the canonical character used in Tate \cite{JT1}) but then if $E/F$ is an extension they always use 
$\psi_E=\psi\circ\mathrm{Tr}_{E/F}$.}. According to Tate \cite{JT2}, Theorem (3.4.1) the functions $\epsilon_{D}(\chi,\psi,dx)$
extends to $\epsilon_{D}(\rho,\psi,dx)$. In order to get {\bf weak inductivity}  we have again to use $\psi_E=\psi\circ\mathrm{Tr}_{E/F}$
if we consider extensions. Then according to Tate \cite{JT2} (3.6.4) the Corollary \ref{Corollary 4.1} turns into 
\begin{cor}\label{Corollary 4.2}
 For the virtual representations of the Weil groups we have 
 \begin{center}
  $W(\rho\omega_{E,s-\frac{1}{2}},\psi_E)=\epsilon_{BH}(\rho,s,\psi_E)=\epsilon_{D}(\rho\omega_{E,s},\psi_E,dx_{\psi_E})$.\\
  $W(\rho,\psi_E)=\epsilon_{BH}(\rho,\frac{1}{2},\psi_E)=\epsilon_{D}(\rho\omega_{E,\frac{1}{2}},\psi_E,dx_{\psi_E})$.
 \end{center}
\end{cor}
Note that on the level of field extension $E/F$ we have to use $\omega_{E,s}$ which is defined as 
\begin{center}
 $\omega_{E,s}(x)=|x|_{E}^{s}=q_{E}^{-s\nu_E(x)}.$
\end{center}
We also know that $q_{E}=q_{F}^{f_{E/F}}$, and $\nu_E=\frac{1}{f_{E/F}}\cdot\nu_F(N_{E/F})$ (cf. \cite{FV}, p. 41, Theorem 2.5), therefore
we can easily see that 
$$\omega_{E,s}=\omega_{F,s}\circ N_{E/F}.$$

Because the norm map $N_{E/F}:E^\times\to F^\times$ corresponds using class field theory to the injection map $G_E\hookrightarrow G_F$,
Tate \cite{JT2} beginning from (1.4.6), simply writes $\omega_{s}=||^s$, and consider $\omega_{s}$ as an unramified character of the Galois
group (or of the Weil group) instead as a character on the field.
Then Corollary \ref{Corollary 4.2} turns into 
\begin{equation}\label{eqn 2.3.15}
 W(\rho\omega_{s-\frac{1}{2}},\psi_E)=\epsilon_{BH}(\rho,s,\psi_E)=\epsilon_{D}(\rho\omega_{s},\psi_E,dx_{\psi_E}),
\end{equation}
for all field extensions, where $\omega_{s}$ is to be considered as one-dimensional representation of the Weil group $W_E\subset G_E$
if we are on the $E$-level. The left side equation (\ref{eqn 2.3.15}) is the $\epsilon$-factor of Langlands
(see \cite{JT2}, (3.6.4)). 

\item The functional equation (\ref{eqn 2.3.9}) extends to 
\begin{equation}\label{eqn 2.3.23}
 W(\rho,\psi)\cdot W(\rho^{\vee},\psi)=\mathrm{det}_{\rho}(-1),
\end{equation}
where $\rho$ is any virtual representation of the Weil group $W_F$, $\rho^{\vee}$ is the contragredient, and $\psi$ is any nontrivial additive
character of $F$. This is formula (3) on p. 190 of \cite{BH} for $s=\frac{1}{2}$.

\item Moreover, the transformation law \cite{JT2} (3.4.5) can (on the $F$-level) be written as\\
\textbf{unramified character twist}
\begin{equation}
 \epsilon_{D}(\rho\omega_{s},\psi,dx)=\epsilon_{D}(\rho,\psi,dx)\cdot \omega_{F,s}(c_{\rho,\psi})
\end{equation}
for any $c=c_{\rho,\psi}$ such that $\nu_F(c)=a(\rho)+n(\psi)\mathrm{dim}(\rho)$. 
It implies that also for the root number on the $F$-level
we have 
\begin{equation}
 W(\rho\omega_{s},\psi)=W(\rho,\psi)\cdot \omega_{F,s}(c_{\rho,\psi}).
\end{equation}
\item {\bf Deligne-Henniart's twisting formula:}\label{Deligne-Henniart's twisting formula}
Let $\rho_1$, and $\rho_2$ be two finite-dimensional representations of $W_F$. Now the question is\\
{\it Is there any explicit 
formula for $W(\rho_1\otimes\rho_2,\psi)$?} \\
The answer to this question is not yet known. However,  under some special conditions, when any of $\rho_1$
and $\rho_2$ is one dimensional, then Deligne gives an explicit formula for $W(\rho_1\otimes\rho_2,\psi)$
(cf. \cite{D2}), and which is as follows:\\
Let $\rho_1=\rho$ be a finite-dimensional representation of $W_F$, and let $\rho_2=\chi$ be any nontrivial character
of $F^\times$. For each $\chi$ there exists an element $c\in F^\times$ such that
$$\chi(1+y)=\psi(cy)\quad\text{for sufficiently small $y$}.$$
For all $\chi$ with sufficiently large {\bf conductor}, we have the following formula:
\begin{equation}\label{Deligne's general twisting formula}
 W(\rho\otimes\chi,\psi)=W(\chi,\psi)^{\dim(\rho)}\cdot \det(\rho)(c^{-1}).
\end{equation}
Now if we define a virtual representation $\rho_0:=\rho-\dim(\rho)\cdot 1_{W_F}$, where $1_{W_F}$ is the trivial representation
of $W_F$, then from above equation (\ref{Deligne's general twisting formula}) we have 
$$W(\rho_0\otimes\chi,\psi)=\det(\rho_0)(c^{-1}).$$
In \cite{DH}, Deligne, and Henniart generalize the above result (see Section 4 of \cite{DH}), for 
for virtual representations $\rho$ of dimension $0$, in which $\chi$ is 
replaced by a representation $\rho'$.

\end{enumerate}

\section{{\bf Applications, and Open Problems}}

In this section, we study some applications of local constants, and some open problems regarding local constants.
Because local constants can be attached to every finite-dimensional complex representation
of a local Galois group, in the Langlands program, the local constants play a very important
role. These local constants appear in various places in modern number theory (in general, in algebraic curves,
local/global Galois representations etc.). In fact,  it is believed that if a mathematical object which
has $L$-function, then we can attach this constant with it.

However, the construction of these root
numbers for every mathematical object which has $L$-function is a difficult problem in
number theory. For instance, so far, we are not able to bring root numbers in the mod-$p$ Langlands
correspondence (cf. \cite{CB}), and the geometric Langlands correspondence (cf. \cite{GL}, \cite{VL}).

The explicit computation of these root numbers has many applications in modern number theory.
For instance, in \cite{MJT}, Taylor showed that in the theory of the structure of ring of algebraic
integers as a Galois module, the local root numbers determine whether or not this projective
module is free, in the case of a tame extension. Furthermore, because we know that the local
Langlands correspondence preserves the local root number, in the Langlands program, the local
root number plays an important role for checking the local Langlands conjecture.

Besides this, on the automorphic side, by extensive study of root numbers, we can classify
the automorphic representations (cf. \cite{DJ}, \cite{JWC}, \cite{CPS}, \cite{CPS98}). 
And it is also expected that on the Galois side, we also can
do the same. Although, so far on the Galois side, we do not have any such complete result
except Heiermann's result \cite{VH}. 

\subsection{{\bf Applications}}

As to application, in this article, we will only review {\bf converse theorems},
and Taylor's result \cite{MJT}. \\

{\bf Classical Converse Theorem in number theory}\\
{\it How to construct a modular form from a given Dirichlet series with 'nice' properties (e.g., analytic continuation,
moderate growth, functional equation), i.e., starting with the series 
$$L(s)=\sum_{n=1}^{\infty}\frac{a_n}{n^s},$$
under what conditions is the function 
$$f(z)=\sum_{n=1}^{\infty}a_n e^{2\pi i nz}$$
a modular form for some Fuchsian group?}\\
The answer to this question is known as the {\it classical converse theorem} in number
theory (cf. \cite{HH}, \cite{EH}, \cite{AW}).
The classical converse theorems establish 
a one-to-one correspondence between ``nice'' Dirichlet series, and automorphic functions. 
Traditionally, the converse theorems 
have provided a way to characterize Dirichlet series associated to modular forms in terms of their analytic properties.

The modern version of classical converse theorems are stated in terms of automorphic representations 
instead of modular
forms. Again, we know that using Langlands local correspondence that automorphic representations are associated with 
Galois representations. Therefore, one can ask the following questions:\\
{\it 1. Are there any converse theorems for automorphic representations (automrphic side of the converse theorem)?\\
2. Similarly, are there any converse theorems for Galois representations (Galois side of the converse theorem)?}\\

{\bf Note:} Because local root numbers for complex Galois representations are the main theme of the paper,
here we only discuss 
converse theorems on the Galois side. For local converse theorems on the automorphic side, one can
follow Dihua Jiang's paper {\it on local $\gamma$-factors} \cite{DJ}.

{\bf 1.  Converse theorem on the Galois side:}\\

So far on the Galois side, we do not have any converse theorem similar to $GL_n$ side except Volker Heiermann's \cite{VH}
work. Here, we summarize his work.

Let $k$ be a field, and $G$ be a group, on which a decreasing filtration by normal subgroups $G^\theta$ with 
$\theta\in [-1,+\infty)$ is defined. Furthermore, assume $G=G^{-1}$, and $G/G^0$ is {\bf cyclic}.
Volker Heiermann considers representations of $G$ in $k$-vector spaces of finite dimension which are trivial on one 
$G^\theta$. For such an indecomposable representation $\sigma$, define 
$$S(\sigma):=\mathrm{sup}\{\{\epsilon:\, \sigma|_{W^\epsilon}\not\supset{\bf 1}\}, 0\}\cdot\deg(\sigma),$$
(for arbitrary $\sigma$, define $S(\sigma)$ by additivity). Then he proves the following inequality: 
\begin{equation}\label{Heiermann's inequality}
 \frac{S(\sigma\otimes\tau^{\vee})}{\deg(\sigma)\deg(\tau)}\le\mathrm{Max}
 \Big\{\frac{S(\sigma\otimes\rho^{\vee})}{\deg(\sigma)\deg(\rho)}, \frac{S(\rho\otimes\tau^{\vee})}{\deg(\rho)\deg(\tau)}\Big\},
\end{equation}
where $\sigma,\tau,\rho$ are three indecomposable $k$-representations of $G$ of finite dimension which are trivial on one 
$G^\theta$, and where it is assumed that these restrictions to $G^\theta$ are semisimple for every $\theta>0$.

Finally, the author also proves, under some of conditions on the residue field of $F$, that an irreducible representation
$\sigma$ on $W_F$ is determined by the set of $\epsilon$-factors of the form $W(\sigma\otimes\tau)$ 
with $\deg(\sigma)\le\deg(\tau)$, up to equivalence.

\begin{thm}[Lemma 1, and Proposition on p. 4 of \cite{VH}]\label{Heiermann's Lemma regarding epsilon}
 Suppose $\sigma$, and $\sigma'$ are two irreducible representations of $W_F$. Suppose that we have 
 $$W(\sigma\otimes\tau,s)=W(\sigma'\otimes\tau,s)$$
 for all irreducible representations $\tau$ with 
 $$\deg(\tau)\le \mathrm{sup}\{\deg(\sigma),\deg(\sigma')\}$$
 of $W_F$, and for all $s\in\bbC$. Then there exists an unramified character $\mu$ of $F^\times$ for which 
 $$\sigma'=\sigma\otimes\mu.$$
 Moreover, if the residue field of $F$ contains more than $2$ elements or if $\sigma$ is induced from 
 a character of a field of residue field of $>2$ elements, then $\sigma$, and $\sigma'$ are isomorphic.
\end{thm}




\begin{rem}
 In Theorem 1.5 in \cite{SABHTF}, one can see a converse theorem on the Galois side regarding local Heisenberg Galois
 representations.
\end{rem}

{\bf 2. Taylor's result regarding the Galois module:}\\ 
Let $F$ be a number field, and $K$ be a tame Galois extension with Galois group $\Gamma:=Gal(K/F)$.
By Noether, it was known that
the ring of integers $O_K$ of $K$ is locally free $\bbZ\Gamma$-module. For a character $\chi\in \widehat{\Gamma}$, we can 
associate the Artin's root number (see Appendix 6.1 below) $W(\chi)$. Given a locally free $\bbZ\Gamma$-module $X$, $(X)$ denotes the 
class of $X$. Let $Cl(\bbZ\Gamma)$ denote the class group of locally free $\bbZ\Gamma$-modules, and $t(W)$ denotes a class in 
$Cl(\bbZ\Gamma)$ which is defined in terms of the values of Artin's root numbers of symplectic type characters. 

In \cite{AF}, A. Fr\"{o}hlich showed that the class of $O_K$ (when $F=\bbQ$, and $\Gamma$ is a quaternion group of order $8$),
is determined by the sign of the Artin root number of the irreducible symplectic character of $\Gamma$. In \cite{MJT}, Taylor
proved the following Theorem:
\begin{thm}[Taylor, Theorem 1 of \cite{MJT}]
 $(O_K)=t(W)$, and so in particular,\\
 (a) as the Artin root numbers of symplectic type characters are $\pm 1$, $(O_K)^2=1.$\\
 (b) the only obstructions to the vanishing of the class of $O_K$, are the signs of the Artin root numbers of symplectic 
 characters.
\end{thm}

\subsection{{\bf Open Problems}}

\vspace{.1cm}

{\bf {(1) Local constants under restriction:}}\\

Let $F$ be a non-Archimedean local field of characteristic zero. Let $K/L/F$ be a tower of field extensions, 
where $K/F$ is Galois, and $L/F$ is finite (need not be Galois).
Let $G=Gal(K/F)$, and $H=Gal(K/L)$ which is a subgroup of 
$G$ with finite index in $G$. Let $\rho$ be finite-dimensional complex representation of $G$, then \\
{\it What is the relationship between the local constants $W(\rho, \psi)$, and 
$W(Res_{H}^{G}(\rho),\psi_K)$,
where $\psi$ is a nontrivial additive character of F, $\psi_K=\psi\circ Tr_{K/F}$, and $Res_{H}^{G}(\rho)$
is the restriction of $\rho$ to $H$?}

\begin{rem}
\begin{enumerate}
 \item

One can think this restriction problem using Robert Boltje's canonical Brauer extension (cf. \cite{RB}),
because canonical Brauer 
induction commutes with restriction map. If we have a formula for 
$W(\text{Res}_{H}^{G}(\rho),\psi\circ \text{Tr}_{L/F})$, 
then it can be used in the {\bf Gan-Gross-Prasad
conjecture} (cf. \cite{GGP}, \cite{GP1}, \cite{GP2}).
\item Let $K/F$ be a finite extension of $F$. Let $\chi_K$ be a nontrivial multiplicative character of $K$.
Denote 
$$\chi_F:=\chi_K|_{F^\times},\quad \text{i.e., restriction of $\chi_K$ to $F^\times$}.$$
Now the question is: {\it Is there any relation between $W(\chi_F,\psi_F)$ and $W(\chi_K,\psi_K)$, 
where $\psi_K=\psi_F\circ Tr_{K/F}$?}\\
If we can answer this question, it will be useful in the Langlands program. But, so far we do 
not have any explicit answer to this question. In \cite{SABFM}, one can see some results in this direction.

\end{enumerate}

\end{rem}

{\bf (2). General Twisting formula:}\\
Let $\rho_1, \rho_2$ be two complex irreducible representations of $G_F$, and let $\psi$ be a nontrivial additive character 
of $F$. Now the question is: 
\begin{center}
 {\it What is the explicit formula for $W(\rho_1\otimes\rho_2,\psi)$?}
\end{center}

When both $\rho_1,\rho_2$ are one-dimensional, we can see the property \ref{Twisting formula for characters} of 
root numbers in Section 4. When 
one of the $\rho_1,\rho_2$ is one dimensional, we have the 
Deligne-Henniart's twisting Formula (\ref{Deligne-Henniart's twisting formula}).
For a particular 
Heisenberg representation, there is an extension of Deligne-Henniart's twisting formula 
(cf. \cite{SABHTF}).
Besides these, so far, we do not have any known
explicit formula for $W(\rho_1\otimes\rho_2,\psi)$ in literature. When $\rho_1$ and $\rho_2$ are one-dimensional, then 
via local Jacobi sum, we can give explicit formula for $W(\chi_1\chi_2,\psi)$ (cf. \cite{SABTF}).\\
{\bf{(3). Geometric analog, and arithmetic connection of local constants (or epsilon factors):}}\\

{\bf Motivation:} 
The grand unification between various mathematical objects is the main theme of the modern mathematical research. From 
Andr\'{e} Weil's philosophy, we (conjecturally) know that there are bijective correspondences between number
theory, and geometry over finite fields (using zeta function, and local information of number theory), number theory, and 
complex geometry  (using geometric Langlands), and geometry over finite fields, and complex geometry (using cohomology).
And these correspondences preserve some analytic objects (e.g., $L$-functions, $\varepsilon$-factors)
which contain many arithmetic information. Further, we also know that the local root numbers remain invariant
under local Langlands
correspondence. Then it is a natural question to ask:\\

{\it {\bf Question (a):} What is the geometric interpretation/connection
of the local constants? Do they have any role in the geometric Langlands program?} 

Because $L$-functions are very important in modern number theory, and arithmetic geometry,
hence, local root numbers are.
For instance, 
let $E$ be an elliptic curve over $\bbQ$, and $N_{E/\bbQ}$ be the conductor of $E$. We know that the $L$-function $L(E,s)$ has
an analytic continuation to the entire complex plane, and it satisfies the following functional equation:
\begin{equation}\label{eqn 5.8}
 \Lambda(E,s)=\omega(E/\bbQ)\cdot \Lambda(E,2-s),\qquad\text{with $\omega(E/\bbQ)=\pm 1$},
\end{equation}
where 
\begin{equation*}
 \Lambda(E,s):=N_{E/\bbQ}^{\frac{s}{2}}(2\pi)^{-s}\Gamma(s)L(E,s),\qquad \Gamma(s)=\int_{0}^{\infty}t^{s-1}e^{-t}dt.
\end{equation*}
The number $\omega(E/\bbQ)$ in the functional equation (\ref{eqn 5.8}), is called the {\bf root number} 
of $E$, and has a very important conjectural meaning 
in the Birch-Swinnerton-Dyer (BSD) conjecture. The {\bf Parity conjecture} (a weaker version of the BSD conjecture) claims
that 
$$\omega(E/\bbQ)=(-1)^{\text{Rank}(E/\bbQ)}.$$
If we notice, we can see that the exact sign of $\omega(E/\bbQ)$ will tell us whether the rank of $E/\bbQ$ is even or odd, 
hence (partially) the structure of Weil-Mordell group of $E/\bbQ$, and for elliptic curves over $\bbQ$ 
this is the ultimate goal of number theorists.\\

{\bf Question (b): (regarding relationship between $\omega(E/\bbQ)$, and $W(\chi,\psi)$):}

Let $F$ be a local field, and $\chi$ be a
multiplicative character of $F$. Let $\psi$ be a nontrivial additive character of $F$.\\
{\bf (i).} Is there any relation between $W(\chi,\psi)$, and $\omega(E/\bbQ)$?\\
{\bf (ii).} Because $\omega(E/\bbQ)$ is a {\bf sign}, so for a particular choice of elliptic curve with some 
admissible conductor, can we 
give any relation between $\omega(E/\bbQ)$, and the classical Gauss sum?\\
{\bf (iii).} Suppose that for a $1$-dimensional representation $\chi$, we are able to find a connection between 
$W(\chi,\psi)$, and $\omega(E/\bbQ)$. Then can we extend that result for any finite-dimensional 
Galois representations?

\begin{rem}
In \cite{KT}, the authors provide some results regarding elliptic curves, and local root numbers. One also can see
\cite{DER} for more information about the geometric analog of root numbers.
\end{rem}

\section{{\bf Appendix}}

\subsection{Global constants (or global epsilon factors/Artin root numbers)}

Let $F$ be a global field, and 
$$\psi: \mathbb{A}_F/F\to\bbC^\times$$
be a nontrivial additive character, and $dx$ the Haar measure on $\bbA_F$ such that 
$$\int_{\bbA_F/F}dx=1.$$
Here $\mathbb{A}_F$ is the adele of $F$.
This is called the {\bf Tamagawa measure}. Call $\psi_v$ the local component  of $\psi$ at a place $v$, and let 
$$dx=\prod_v dx_v$$
be any factorization of $dx$ into a product of local measures such that the ring $O_v$ of integers in $F_v$ gets measure $1$
for almost all $v$.

Let $\rho$ be a representation of the global Weil group $W_F$, and put 
\begin{equation}\label{eqn 5.1}
 L(\rho,s):=\prod_v L(\rho_v\otimes\omega_s)
 \end{equation}
where $s\in\bbC$, and $\omega_s:F_v^\times\to\bbC^\times$ is the unramified character such that 
$$\omega_s(\pi_v)=q_v^{-s}$$
for $\pi_v$ a local prime element, and $q_v$ the order of the corresponding residue field.

Moreover, put 
\begin{equation}\label{eqn 5.2}
 W(\rho,s):=\prod_v W(\rho_v\otimes\omega_s,\psi_v,dx_v),
\end{equation}
where the existence of the non-abelian local constant on the right side is assumed. Then
\begin{thm}
 The product (\ref{eqn 5.1}) converges for $s$ in some right half-plane, and defines a function $L(\rho,s)$ which is meromorphic
 in the whole $s$-plane, and satisfies the functional equation
 \begin{equation}\label{Global functional equation}
  L(\rho,s)=W(\rho,s)\cdot L(\rho^\vee,1-s)
 \end{equation}
where $\rho^\vee$ is the dual of $\rho$.
\end{thm}
\begin{rem}
 The global functional equation $(\ref{Global functional equation})$ does {\bf not} come from local functional equation
 (local (\ref{Global functional equation}) is wrong):
 $$L(\rho_v\otimes\omega_s)\ne W(\rho_v\otimes\omega_s,\psi_v,dx_v)\cdot L(\rho_v^*\otimes\omega_{1-s}).$$
 A local functional equation is known only for $\dim(\rho)=1$ (cf. \cite{JTPT}, \cite{TT})
 but not for higher dimensions, and the local functional
 equation for dimension $1$ looks different (cf. Equation (\ref{eqn 4.4}))
 \begin{equation}\label{eqn 2.3.2}
 \frac{\int\hat{f}(x)w_1\chi^{-1}(x)d^\times x}{L(w_1\chi^{-1})}=W(\chi,\psi,dx)
 \frac{\int f(x)\chi(x)d^\times x}{L(\chi)},
\end{equation}
where $\chi$ stands for the $1$-dimensional representation $\rho_v\otimes\omega_s$, and $\chi^{-1}$ stands for 
$\rho_v^*\otimes\omega_{-s}$.

One obtains (\ref{Global functional equation}) by taking the product of the local functional equation (\ref{eqn 2.3.2}) 
over all places $v$, and then verifying that the product of the numerators on the left (over all places $v$) is the {\bf same}
as the product of the numerators on the right. Therefore, in the global version of the numerators can be canceled, and the 
(\ref{Global functional equation}) follows. For higher dimensional $\rho$ we have no local functional equation at all.
Deligne's proof for existence, and for product formula (\ref{eqn 5.2}) is different.
\end{rem}


\begin{rem}

Let $F$ be an algebraic number field, and $K$ be a finite normal extension of $F$ with Galois group $Gal(K/F)=:G$.

This  global constant satisfies the following properties:\\ 
{\bf 1. Additivity:}
$$W(\rho_1\oplus\rho_2)=W(\rho_1)W(\rho_2)$$
for representations $\rho_1,\rho_2\in R(G)$.\\
{\bf 2. Invariant under inflation:} Let $E/F$ be a finite normal extension with $K\subseteq E$. Then 
$$W(Inf_{G/N}^{G}(\rho))=W(\rho)$$
where $N=Gal(K/E)$, and $Inf_{G/N}^{G}(\rho)$ is the representation of $N$ defined by inflation from $\rho$.\\
{\bf 3. Invariant under induction:} Let $L/F$ be a subextension in $K/F$, and $\rho_L$ be a representation of $Gal(K/L)$. Then 
$$W(Ind_{Gal(K/L)}^{Gal(K/F)}(\rho_L))=W(\rho_L).$$

\end{rem}

\subsection{The canonical Brauer induction formula on the Galois side}.

Let $F$ be an arbitrary field, and $\bar{F}$ its 
separable algebraic closure. In this paper, we denote $K/F$ for the tower $F\subset K\subset\bar{F}$, and $K$ 
is a finite Galois 
extension of $F$ whose Galois group is denoted by $G_{K/F}=\mathrm{Gal}(K/F)$. Let $G_F$ be the absolute
Galois group, which is a profinite group with $G_F=\varprojlim_{K/F}G_{K/F}$.

We denote $R(G_F)$ the Grothendieck ring of finite-dimensional complex representations of $G_F$ with kernel of finite 
index, and let $\widehat{G_F}$ be the group of linear complex representations of $G_F$ with kernel of finite index, 
then we have 
direct limits
\begin{itemize}
  \item $R(G_F)=\varinjlim_{K/F}R(G_{K/F})$.
  \item $\widehat{G_F}=\varinjlim_{K/F}\widehat{G_{K/F}}$ the group of linear characters.
  \item $R_{+}(G_F)=\varinjlim_{K/F}R_{+}(G_{K/F})$, free abelian group with basis $\mathcal{M}_{G_F}/G_F$
  where $\mathcal{M}_{G_F}=\varinjlim_{K/F}\mathcal{M}_{G_{K/F}}$.
  \item $b_{G_F}=\varinjlim_{K/F}b_{G_{K/F}}:R_{+}(G_F)\rightarrow R(G_F)$, which is well defined because $b_G$ 
  commutes with \emph{inflations}.
\item The section $a_{G_F}$ is $R(\widehat{G_F})=\mathbb{Z}[\widehat{G_F}]$-linear, and commutes with restrictions to 
subgroups.
\end{itemize}
  \begin{lem}[\cite{RB}, p. 27, Corollary 2.21]\label{Lemma 2.2}
    Let $N$ be a normal subgroup of $G$, and $\bar{G}=G/N$. Let $\mathrm{Inf}_{\overline{G}}^{G}(\chi)\in R(G)$ come
    by inflation along $G\to \bar{G}$ 
    for some $\chi\in R(\bar{G})$. Then we have 
   \begin{equation}\label{eqn 2.17}
 a_{G}(\mathrm{Inf}_{\overline{G}}^{G}(\chi))=\sum_{[\bar{H},\bar{\varphi}]}\alpha_{[\bar{H},\bar{\varphi}]}(\chi)[H,\varphi].
\end{equation}
This means $\alpha_{[H,\varphi]}(\mathrm{Inf}_{\overline{G}}^{G}(\chi))=0$ unless $N\leq H\leq G$
and $\varphi|_N=1_N$, and in this case 
$\alpha_{[H,\varphi]}(\mathrm{Inf}_{\overline{G}}^{G}(\chi))=\alpha_{[\bar{H},\bar{\varphi}]}(\chi)$.
\end{lem}
According to Lemma \ref{Lemma 2.2}, the section $a_G$ commutes with inflation, and 
if $\chi\in R(\overline{G})$ lives on a factor group $G\rightarrow \overline{G}$, and where $H\subset G$ denotes the full
preimage of $\overline{H}\subset\overline{G}$, and $[H,\varphi]$ is always inflated from 
$[\overline{H},\overline{\varphi}]$. Therefore, we can go to the profinite Galois group $G_F$, and obtain a Brauer 
induction formula for $G_F$ by
\begin{equation*}
 a_{G_F}=\varinjlim_{K/F}a_{G_{K/F}}:R(G_F)\rightarrow R_{+}(G_F).
\end{equation*}

\subsection{{\bf Extendible  function on the Galois side}}

Let $K/F$ be a finite Galois extension of the field $F$.
Let $\omega=(\omega_K)_{K/F}$ be an abelian invariant for $F$. Then,
\begin{itemize}
 \item An \textbf{extension} $W=(W_K)_{K/F}$ of $\omega$ is a family of maps
 \begin{center}
  $W_{K}:R(G_K)\rightarrow \mathcal{A}$,
 \end{center}
where $W_K$ extends $\omega_K:\widehat{G_K}\rightarrow \mathcal{A}$, and such that $W_{sK}(^s\chi)=W_{K}(\chi)$.
 \item A \textbf{strong (resp. weak) extension} of $\omega$ is an extension $W=(W_K)_{K/F}$ of $\omega$
 such that $W$ is invariant under induction (resp. induction in dimension zero), i.e., for all finite field 
 extensions $F\subset K\subset L$, and $\chi\in R(G_L)$(plus $\chi(1)=0$ in the weak case), we have:
 \begin{equation}\label{eqn 2.17}
  W_{K}(\mathrm{Ind}_{G_L}^{G_K}(\chi))=W_{L}(\chi).
 \end{equation}
 Suppose the $\mathrm{dim}\,\chi\neq0$, let $\chi_0=\chi-\mathrm{dim}\,\chi\cdot 1_{G_L}$ which a virtual representation
 of $G_L$ whose dimension is zero. Then now we use equation (\ref{eqn 2.17}) for $\chi_0$, and have:
 \begin{equation}
   W_{K}(\mathrm{Ind}_{G_L}^{G_K}(\chi))=\big(\lambda_{G_L}^{G_K}\big)^{\mathrm{dim}\,\chi}W_{L}(\chi).
 \end{equation}
where 
\begin{equation}
 \lambda_{G_L}^{G_K}:=\frac{W_{K}(\mathrm{Ind}_{G_L}^{G_K}1_{G_L})}{W_{L}(1_{G_L})}.
\end{equation}

\item $\omega$ is called strongly (resp. weakly) extendible , if there exists a strong (resp. weak)
extension of $\omega$.
\end{itemize}
\begin{rem}
 By definition, strong extendibility implies weak extendibility. For finite groups $G$, and $\chi\in R(G)$, we have
 \begin{equation}\label{eqn 2.19}
  \chi-\chi(1)\cdot 1_G=\sum_{[H,\varphi]\in \mathcal{M}_G/G}\alpha_{[H,\varphi]}(\chi)\mathrm{Ind}_{H}^{G}(\varphi-1_H).
 \end{equation}
Now then if $\chi(1)=0$, we have from equation (\ref{eqn 2.19})
\begin{equation}
 \chi=\sum_{[H,\varphi]\in \mathcal{M}_G/G}\alpha_{[H,\varphi]}(\chi)\mathrm{Ind}_{H}^{G}(\varphi-1_H).
\end{equation}
Moreover,
\begin{center}
 $W(\mathrm{Ind}_{H}^{G}(\varphi-1_H))=W(\varphi-1_H)=w(\varphi-1_H)$,
\end{center}
because $\varphi-1_H$ is of dimension $0$, and then $W$ can be replaced by $w$ because 
$\varphi$, and $1_H$ are linear characters.
Therefore, by Equation (\ref{eqn 2.19}), we see that $W(\chi)$ is uniquely determined for characters $\chi$ of finite
factor group $G\leftarrow \bar{G}$, and as we know each $\chi\in R(\bar{G})$ lives on a finite factor group. 
$W$ is uniquely determined, if it exists. This also proves that 
 the weak extensions $W$, if existing, are unique. 
\end{rem}

\subsection{Galois invariant, and abelian invariant}

Assume now that for any finite field extension $K/F$ we have $G_K=\mathrm{Gal}(\bar{F}/K)$, and map
\begin{center}
 $\omega_{K}:\widehat{G_K}\rightarrow \mathcal{A}$, $\varphi\mapsto \omega_{K}(\varphi)$
\end{center}
with values in an abelian group, which is \textbf{Galois invariant:}
\begin{center}
 $\omega_{sK}(^s\varphi)=\omega_{K}(\varphi)$
\end{center}
for all $s\in G_F$. We call such a family $\omega=(\omega_{K})_{K/F}$ an \textbf{abelian invariant} for $F$. 
The abelian 
invariants for $F$ are in $1-1$ correspondence with the maps
\begin{center}
 $\widetilde{\omega}_F:R_{+}(G_F)\rightarrow \mathcal{A}$, $[G_K,\varphi]\mapsto \omega_{K}(\varphi)$.
\end{center}
Because of the induction maps $\mathrm{Ind}_{+G_K}^{G_F}:R_{+}(G_K)\rightarrow R_{+}(G_F)$ the map
$\widetilde{\omega}_F$ will induce the maps
\begin{equation}
 \widetilde{\omega}_K:=
 \widetilde{\omega}_F\circ\mathrm{Ind}_{+G_K}^{G_F}:R_{+}(G_K)\rightarrow \mathcal{A} \label{eqn 2.18}
\end{equation}
such that 
\begin{center}
 $\widetilde{\omega}_{sK}([^s(H,\varphi)]_{sK})=\widetilde{\omega}_{K}([H,\varphi]_K)$,
\end{center}
for all $s\in G_F$, and $[H,\varphi]\in R_{+}(G_K)$, and $(H,\varphi)\in\mathcal{M}_{G_K}$.

\subsection{The canonical extension}

If $G_F$ is the absolute Galois group over $F$, then by definition, we have 
\begin{equation*}
 a_{G_F}=\varinjlim_{K/F}a_{G_{K/F}}:R(G_F)\rightarrow R_{+}(G_F)
\end{equation*}
the direct limit over the inflation maps $R(G_K)\rightarrow R(G_F)$ for finite factor groups 
$G_K\leftarrow G_F$.
Therefore, beginning from an abelian invariant $\omega=(\omega_K)_{K/F}$ over the base field $F$, we can use the 
compatible system $a_{G_K}:R(G_K)\rightarrow R_{+}(G_K)$ for all finite $K/F$ \textbf{to define}
\begin{equation}
 W_K:\widetilde{\omega}_K\circ a_{G_K}
\end{equation}
where $\widetilde{\omega}_K$ comes from $\widetilde{\omega}_F$ as in (\ref{eqn 2.18}). We call $W=(W_K)_{K/F}$ the 
\textbf{canonical extension} of $\omega$.

The Galois invariance of $\omega_K$, and $\widetilde{\omega}_K$ resp.(see equation (\ref{eqn 2.18})) implies that the system
$(W_K)_{K/F}$ is also Galois invariant. Therefore, the definition of $W$ gives 
an extension of $\omega$. Here, the canonical extension is defined for any arbitrary extendible  function.


\vspace{1cm}

\textbf{Acknowledgements.} I would like to thank Prof. E. W. Zink for his continuous suggestions, and comments for this work
and giving me a copy of his edited version of the Helmut Koch's paper on extendible function. 


\begin{thebibliography}{99}

  
\bibitem{SABFM}S. A. Biswas, Epsilon factors of symplectic type characters in the wild case.
Forum Math.33 (2021), no.2, 569-577.
\bibitem{SABTF} S. A. Biswas, Twisting formula of epsilon factors, Proc. Indian Acad. Sci. (Math. Sci.) Vol. {\bf 127}, No. 4,
September 2017, pp. 585-598. DOI 10.1007/s12044-017-0350-7.
\bibitem{SABLT} S.A. Biswas, Lamprecht-Tate formula, \url{https://arxiv.org/pdf/1702.04286.pdf}.

 
 
\bibitem{SAB1}S. A. Biswas, Computation of the Lambda function for a finite Galois extension, 
Journal of Number Theory, Volume {\bf 186}, May 2018, Pages 180-201.
\bibitem{SABLTRQ}S. A. Biswas, Langlands lambda function for quadratic tamely ramified extensions,
Journal of Algebra, and its Applications, 18 (2019), no. 7, 1950132, 10 pp,
DOI:\url{https://doi.org/10.1142/S0219498819501329}.
\bibitem{SABT}S. A. Biswas, Local constants for Galois representations - Some explicit results, Ph.D thesis, 
\url{https://arxiv.org/pdf/1603.06089.pdf}.
\bibitem{SABHTF} S. A. Biswas, An extension of Deligne-Henniart's twisting formula and its applications,
\url{https://arxiv.org/abs/2101.08118}.

 
\bibitem{RB} Robert Boltje, Canonical, and explicit Brauer induction in the character ring of finite group, and a 
 generalization for Mackey functors, Ph.D thesis, 1989, Universit\"{a}t Augsburg.
 \bibitem{RB2}Robert Boltje, A canonical Brauer induction formula, Ast\'{e}risque No. 181-182 (1990), 5, 31–59.


 
 
\bibitem{Br1} R. Brauer, On the Artin's $L$-series with general group characters, Ann. of Math. {\bf 48} (1947), 502-514. 
 



 
\bibitem{CB} Christophe Breuil, The emerging p-adic Langlands programme, Proceedings of the International Congress of 
Mathematicians, Hyderabad, India, 2010. 

\bibitem{BH} C.J. Bushnell, G. Henniart, The local Langlands conjecture for $GL(2)$, Springer-Verlag, 2006. 


 
 
 

 
 \bibitem{JWC}J. W. Cogdell, L-functions, and converse theorems for $GL_n$, IAS/Park City Mathematics Series, 
 Volume 12, 2002, pp. 97-177.
  \bibitem{CPS}J. W. Cogdell, I. I. Piatetski-Shapiro, Converse Theorems, Functoriality, and Applications to Number Theory,
  ICM 2002, Vol. III. 1-3, \url{http://arxiv.org/pdf/math/0304230.pdf}.


 
 \bibitem{CPS98}J. Cogdell, and I. Piatetski-Shapiro, Converse theorems for $GL_n$, 
Inst. Hautes Etudes Sci. Publ. Math. No. {\bf 79} (1994).






 

 
 

\bibitem{D1}P. Deligne, Les constantes des \'{e}quations fonctionnelle des fonctions L, in Modular functions of 
one variable II, Lecture Notes in Mathematics {\bf 349} (1972), 501-597, Springer-Verlag, Berlin-Heidelberg-New York. 
\bibitem{D2}\textemdash, Les constantes locale de l'equation fonctionnelle de la fonctions L d'Artin d'une 
representation orthogonale, Invent. Math. {\bf 35} (1976), 299-316.
 



\bibitem{DH} P. Deligne, and G. Henniart, 
Sur la variation, par torsion, des constantes locales d'équations fonctionnelles de fonctions L. (French) 
[On the variation through torsion of the local constants of functional equations of L functions] 
Invent. Math. {\bf 64} (1981), no. 1, 89-118.  
 
 
\bibitem{FV} I.B. Fesenko, S.V. Vostokov, Local fields, and their extensions, Second edition 2001.

\bibitem{AF} A. Fr\"{o}hlich, Artin root numbers, and normal integral bases for quaternion fields, Invent. Math. {\bf 17}, 143-166,
(1972).











 

 
 
 
 
 


 \bibitem{GP1} B. Gross, D. Prasad, On the decomposition of a representation of $SO_n$ when 
restricted to $SO_{n-1}$, Canad. J. Mathematics. {\bf 44} (1992), No. 5, 974-1002.
\bibitem{GP2} B. Gross, D. Prasad, On irreducible representations of 
$SO_{2n+1}\times SO_{2n}$, Canad. J. Math. {\bf 46} (1994), No. 5, 930-950.
 \bibitem{GGP}W.T. Gan, B. Gross, D. Prasad, Symplectic local root numbers, central critical L-values, and restriction
 problems in the representation theory of classical group, Astérisque {\bf 346} (2012), pp. 1-109.

 
\bibitem{HH} H. Hamburger, \"{U}ber die Riemannsche Funktionalgleichung der $\zeta$-funktion.
Math. Zeit. {\bf 10} (1921), 240-254; {\bf 11} (1922), 224-245; {\bf 13} (1922), 283-311.



\bibitem{EH} E. Hecke, \"{U}ber die Bistimmung Dirichletscher Reihen durch ihre Funktionalgleichung. 
Math. Ann. {\bf 112} (1936), 664-699.


\bibitem{VH} Volker Heiermann, Sur l'espace des repr\'{e}sentations irr\'{e}ducible du groupe de Galois d'un corps local,
C. R. Acad. Sci. Paris S\'{e}r. I Math. {\bf 323} (1996), no. 6. 571-576.




 




\bibitem{H00} G. Henniart, Une preuve simple des conjectures de Langlands pour $GL(n)$ sur un corps $p$-adique, Invent. Math.
{\bf 139} (2000), no. 2, 439-455.

\bibitem{GH84}G. Henniart, Galois $\epsilon$-factors modulo roots of unity, Invent. Math. {\bf 78}, 117-126 (1984).


 
 


\bibitem{DJ} D. Jiang, On local $\gamma$-factors.
Arithmetic geometry, and number theory, 1-28, Ser. Number Theory Appl., 1, World Sci. Publ., Hackensack, NJ, 2006. 

 

 



 
\bibitem{HK} H. Koch, Extendible  Functions, Centre Interuniversitaire en Calcul; Math\'{e}matique 
 Alge\'{e}brique, Concordia University in September 1990.
 
 \bibitem{HKEWZ} H. Koch, and E.-W. Zink, Extendible functions, and local root numbers: Remarks on a paper R. P. Langlands,
 Mathematische Nachrichten, DOI: 10.1002/mana.202200391.
 
\bibitem{KT}K. Kramer, and J. Tunnell, Elliptic curves, and local $\epsilon$-factors, Composition Mathematica,
{\bf 46} (1982), no. 3, p. 307-352. 
 
 
 \bibitem{VL}V. Lafforge, Shtukas for reductive groups, and Langlands correspondence for function fields, 
 \url{http://vlafforg.perso.math.cnrs.fr/files/cht-ICM-lafforgue.pdf}.
 
\bibitem{EL} E. Lamprecht, Allgemeine Theorie der Gauss'schen Summen in endlichen Kommutativen Ringen, 
Math. Nach. {\bf 9}(1953).

\bibitem{RL}R.P. Langlands, On the functional equation of the Artin $L$-functions, unpublished article (1970)
\url{https://publications.ias.edu/sites/default/files/a-ps.pdf}.
 
 
 
 


\bibitem{GL}G. Laumon, Correspondance de Langlands g\'{e}om\'{e}trique pour les corps de fonctions,
Duke Math. Jour., vol. {\bf 54} (1987), 309-359.


\bibitem{M}J. Martinet, Character theory of Artin $L$-functions, Algebraic Number Fields
 (L-functions, and Galois properties), Proceedings of Symposium, Edited by A. Fr\"{o}hlich, pp. 1-88.

 
   
   








\bibitem{DER} D. E. Rohrlich, Elliptic curves, and the Weil-Deligne group, Centre de Recherches Math\'{e}matiquea, 
CRM Proceedings, and Lecture Notes, Volume {\bf 4}, 1994, pp. 125-157.
 
   
\bibitem{TS}T. Saito, Local constant of $\text{Ind}_{K}^{L}1$, Comment. Math. Helvetici {\bf 70} (1995),
507-515.


\bibitem{JPS}J.-P. Serre, Local fields, Graduate Texts in Mathematics, {\bf 67}, Springer-Verlag, New York Inc, 1979.
\bibitem{JPSLR}J.-P. Serre, Linear representations of finite groups, Graduate Texts in Mathematics, {\bf 42}, Springer.




\bibitem{TT} Terence Tao, Tate's proof of the functional equation, 
\url{https://terrytao.wordpress.com/2008/07/27/tates-proof-of-the-functional-equation/#more-430}.

\bibitem{JT1}J. Tate, Local Constants, Algebraic number fields: L-functions, and Galois properties 
(Proc. Sympos., Univ. Durham, Durham, 1975), pp. 89-131.
\bibitem{JTPT}J. Tate, Fourier analysis in number fields, and Hecke's zeta-functions (Ph.D Thesis, 1950), Chapter XV,
Algebraic Number Theory (Proceedings of an  Instructional Conference, 
Edited by J.W.S. Cassels, and A. Fr\"{o}hlich, 1967), Academic Press.


\bibitem{JT2}J. Tate, Number theoretic background, Proc. Sympos. Pure Math., {\bf 33}, Automorphic forms, representations
and $L$-functions Part 2, pp. 3-26, Amer. Math. Soc., 1979.


 \bibitem{MJT} M.J. Taylor, On Fr\"{o}hlich's conjecture for rings of integers of tame extensions, Invent. Math. {\bf 63},
 41-79, (1981).

 
\bibitem{AW}A. Weil, \"{U}ber die Bistimmung Dirichletscher Reihen durch Funktionalgleichungen. Math. Ann. {\bf 168} (1967), 149-156.

 



 

  
 
 




 









 

 
 
 
 
  
 
 
 



















 
 
 































 








\end{thebibliography}
\end{document}